2
 1
 2
 1
 2
 1
\font\BibAuthorFont = cmbxsl10 scaled \magstep 0
 0 
\font\Smallcaps = cmcsc10 scaled \magstep 0 
 0 
\font\TitleFont = cmss17 scaled \magstep 0 
\font\SectionHeadingFont = cmss12 scaled \magstep 1 
\font\AuthorFont = cmss10 scaled \magstep 0 

\input amssym

{\TitleFont 
\noindent Weighted fourth moments  

\medskip 

\noindent of Hecke zeta functions with groessencharacters}

\vskip 4mm

\noindent{\AuthorFont Nigel Watt}

\vskip 8 mm



\noindent{\bf Abstract:}\ By using recently obtained bounds for 
certain sums of Kloosterman sums we prove new bounds for the 
sum $\sum_{-D\leq d\leq D}\int_{-D}^{D}|\zeta\left( 1/2 +it,\lambda^d\right)|^4
|\sum_{0<|\mu|^2\leq M}A(\mu) \lambda^d((\mu)) |\mu|^{-2it}|^2 {\rm d}t$, 
where $\lambda^d$ is the groessencharacter satisfying 
$\lambda^d((\alpha))=\lambda^d(\alpha{\Bbb Z}[i])=(\alpha /|\alpha|)^{4d}$,  for $0\neq\alpha\in{\Bbb Z}[i]$,  
and $\zeta(s,\lambda^d)$ is the Hecke zeta function that satisfies  
$\zeta(s,\lambda^d)=(1/4)\sum_{0\neq\alpha\in{\Bbb Z}[i]} \lambda^d((\alpha)) |\alpha|^{-2s}$ 
for $\Re(s)>1$,  while 
the numbers $D,M\in (0,\infty)$ and function $A : {\Bbb Z}[i]-\{ 0\}\rightarrow{\Bbb C}$ 
are arbitrary (though it is only in respect of cases in which  
$M$ is relatively small, compared to $D$, that our results are new and interesting). 
One of our new bounds may have an application in enabling 
a certain improvement of a result of P.~A. Lewis on the distribution of Gaussian primes. 

\smallskip


\smallskip 


\bigskip

\goodbreak 
\noindent{\SectionHeadingFont 1. Introduction}

\medskip

\noindent In [17,18] Hecke discovered a class of zeta functions with 
applications in the study of the multidimensional distribution of
prime ideals in algebraic number fields. Associated with the Gaussian number field
${\Bbb Q}(i)$, with ring of integers ${\frak O}={\Bbb Z}[i]$,
there is a family 
$(\zeta(s,\lambda^d))_{d\in{\Bbb Z}}$ of functions in Hecke's class 
that satisfy   
$$\zeta\left( s, \lambda^d\right)
=\sum_{{\frak a}\in I}\lambda^d ({\frak a}) \left( N({\frak a})\right)^{-2s}
={1\over 4}\,\sum_{0\neq\alpha\in\frak O}\Lambda^d (\alpha) |\alpha|^{-2s}
\qquad\quad\hbox{($s\in{\Bbb C}$, $\Re (s) >1$ and $d\in{\Bbb Z}$),}\eqno(1.1)$$
where $\Lambda^d$ is (for each $d\in{\Bbb Z}$) the  endomorphism 
of the group $({\Bbb Q}(i))^{*}$ given by 
$$\Lambda^d(\alpha)=\left( {\alpha}\over |\alpha|\right)^{\!\! 4d}
\qquad\quad\hbox{($0\neq\alpha\in{\Bbb Q}(i)$),}\eqno(1.2)$$
while the elements of the set $I$ are the non-zero ideals in ${\frak O}$, 
the `norm' $N$ satisfies $N((\alpha))=|{\frak O}/(\alpha)|=|\alpha|^2\,$ 
(when $(\alpha)=\alpha{\frak O}=\{\alpha\beta : \beta\in{\frak O}\}$), 
and each `groessencharacter' $\,\lambda^d : I\rightarrow{\Bbb C}^*$ is given by: 
$$\lambda^d\left( (\alpha)\right) =\Lambda^{d}(\alpha)\qquad\quad 
\hbox{($0\neq\alpha\in{\frak O}$).}$$  
For each $d\in{\Bbb Z}$ the function $\zeta\bigl( s, \lambda^d\bigr)$ has a meromorphic continuation
to all points $s\in{\Bbb C}$, with the only pole being that of the Dedekind zeta function
$\zeta\bigl( s, \lambda^0\bigr)$
at $s=1$. 

Our subject in this paper is the mean value
$$E(D;M,A)
=\ \ \sum_{\!\!\!\!-D\leq d\leq D}\ \int\limits_{-D}^{D}\left|\zeta\left( 1/2 +it,\lambda^d\right)\right|^4
\left| P_M\left( A;it,\lambda^d\right)\right|^2 {\rm d}t,\eqno(1.3)$$
where it is assumed that $A$ is a mapping from 
${\frak O}-\{ 0\}$ to ${\Bbb C}$, that $0<D,M<\infty$,  
and that 
$$P_M\left( A;s,\lambda^d\right)
=\sum_{0<|\mu|^2\leq M}A(\mu) \Lambda^d(\mu) |\mu|^{-2s}\qquad\qquad 
\hbox{for $d\in{\Bbb Z}$ and $s\in{\Bbb C}$}\eqno(1.4)$$
(note that here, and in all that follows,  
variables of summation designated by Greek letters, such as the `$\mu$' above,   
are ${\frak O}$-valued variables, whereas those designated by Latin letters
are ${\Bbb Z}$-valued variables).
The mean value $E(D;M,A)$ is a `Gaussian integer analogue' of the 
mean value 
$$S(T,N)=\int\limits_0^T \left|\zeta\left( 1/2 +it\right)\right|^4
\Biggl|\,\sum_{0<n\leq N} a_n n^{-it}\Biggr|^{\,2} {\rm d}t$$ 
which has been studied by Iwaniec [23], Deshouillers and Iwaniec 
[8] and Watt [43]. In [7] and [8] a connection was 
established between $S(T,N)$ and the spectral theory of the non-Euclidean 
Laplace-Beltrami operator $\Delta_2$ acting on certain spaces of 
functions defined on $2$-dimensional hyperbolic space ${\Bbb H}_2$ and automorphic 
with respect to the action of some Hecke congruence group 
$\Gamma\leq SL(2,{\Bbb Z})$ on ${\Bbb H}_2$;  
in particular, it was observed by 
Deshouillers and Iwaniec in [8] that, if 
one could assume Selberg's conjecture   
that, except for the eigenvalue $0$ associated with the space of constant functions, 
the relevant eigenvalues of 
$-\Delta_2$ lie in the interval $[1/4,\infty)$, then their method 
would yield the estimate 
$$S(T,N)\ll_{\varepsilon} \left( 1 + T^{-1/2} N^2\right) 
T^{1+\varepsilon}\sum_{0<n\leq N} \left| a_n\right|^2\;,\eqno(1.5)$$ 
for $\varepsilon >0$, $N\geq 1$, $T\geq 1$, and any complex sequence 
$(a_n)_{n\in{\Bbb N}}$. Selberg, in [40], had already shown that the 
eigenvalues in question do all lie in the interval $[3/16,\infty)$, and the  
lower bound $3/16$  in this result has since been improved to 
$975/4096=(1/4)-(7/64)^2$ by Kim and Sarnak [28]. 
\par 
By a non-trivial application of Selberg's lower bound for the 
above mentioned eigenvalues, Deshouillers and Iwaniec obtained, 
in [8, Theorem~1], the unconditional estimate   
$$S(T,N)\ll_{\varepsilon} \left( 1 + T^{-1/2} N^2 + T^{-1/4} N^{5/4}\right) 
T^{1+\varepsilon}\sum_{0<n\leq N} \left| a_n\right|^2\;,\eqno(1.6)$$ 
for $\varepsilon >0$, $N\geq 1$, $T\geq 1$, and any complex sequence 
$(a_n)_{n\in{\Bbb N}}$.    
If Kim and Sarnak's improvement of Selberg's lower bound for the eigenvalues 
had been available at the time, then Deshouillers and Iwaniec 
would have been able to replace the term $T^{-1/4} N^{5/4}\,$  (in brackets, in (1.6))  
by $T^{\alpha -1/2} N^{2-3\alpha}$, where $\alpha =7/64$. 

By refining one aspect of the method used by Deshouillers and Iwaniec 
in their proof of (1.6), we obtained, in [43, Theorem~1], the further 
unconditional estimate   
$$S(T,N)\ll_{\varepsilon} \left( 1 + T^{-1/2} N^2\right) 
T^{1+\varepsilon} N \max_{0<n\leq N} \left| a_n\right|^2\;,\eqno(1.7)$$ 
for $\varepsilon >0$, $N\geq 1$, $T\geq 1$, and any complex sequence 
$(a_n)_{n\in{\Bbb N}}$. For certain applications to the theory of the distribution of 
prime numbers (see, for example, [1], or [15, Chapters~7 and~9]) 
the estimate in (1.7) is 
more effective than that in (1.6), and is (moreover) just as effective as  
the conditional estimate in (1.5) would be (were an unconditional 
proof of (1.5) to be discovered). 
\par  
Our principal new result in this paper (Theorem~1, below) is an upper bound 
for $E(D;M,A)$ analogous to the bound (1.7) for $S(T,N)$. 
Our proof of this new result depends 
on estimates for a certain sum of Kloosterman sums; these 
estimates are supplied by Lemma~23 (below), which is a reformulation 
of a result from our paper [44]. Note that our work in [44] 
depends on the spectral large sieve inequality that we  
proved in [45]. In both [45] and [44], and in the present paper,  
we employ methods analogous to those pioneered by Iwaniec,
in [23], and by Deshouillers and Iwaniec, 
in [7] and [8]. 
In particular, with the aid of a slight extension of Lokvenec-Guleska's 
sum formula [34, Theorem~12.3.2] 
(which itself is a generalisation of 
the sum formula [5, Theorem~13.1] of  
Bruggeman and Motohashi),  we establish, in [44], 
a connection between 
the sum of Kloosterman sums occurring in the equation (7.7) (below)  
and the spectral theory of the 
non-Euclidean Laplace-Beltrami operator 
$$\Delta_3 
=r^2\left( {\partial^2\over\partial x^2} + {\partial^2\over\partial y^2} 
+ {\partial^2\over\partial r^2}\right) - r{\partial\over\partial r}\eqno(1.8)$$ 
acting on  certain spaces 
$${\frak D}={\frak D}_{\Gamma} 
=\left\{ f\in L^2\left(\Gamma\backslash{\Bbb H}_3\right)\cap C^2\left( {\Bbb H}_3\right) 
\,:\, \Delta_3 f\in L^2\left(\Gamma\backslash{\Bbb H}_3\right)\right\}\;,$$ 
the members of which are complex valued functions that are defined on  
${\Bbb H}_3 =\left\{ (x+iy , r)\,:\, x,y\in{\Bbb R} , r>0\right\}\,$  
(the upper half-space model for $3$-dimensional hyperbolic space) 
and that are automorphic with respect to 
the action on ${\Bbb H}_3$ of some Hecke congruence group 
$\Gamma =\Gamma_0(\omega)\leq SL(2,{\frak O})$, where $\omega$ is some non-zero 
Gaussian integer and 
$$\Gamma_0(\omega)=\left\{ \pmatrix{\alpha &\beta \cr \gamma & \delta}\in SL(2,{\frak O})\,:\, 
\gamma\in \omega{\frak O}\right\}\eqno(1.9)$$ 
(see [11, Chapters~1-4] for definitions of the space 
$L^2\left(\Gamma\backslash{\Bbb H}_3\right)$ and other terminology used here). 
\par 
A significant feature of the spectral theory just mentioned is that, 
even if one allows for repetitions (consistent with the relevant multiplicities), 
the eigenvalues of the operator 
$-\Delta_3 : {\frak D}_{\Gamma}\rightarrow L^2 (\Gamma\backslash{\Bbb H}_3)$ are 
the terms of an unbounded monotonically increasing sequence $(\lambda_n(\Gamma))_{n\in{\Bbb N}}$ 
such that $\lambda_1(\Gamma)>\lambda_0(\Gamma)=0\,$ (see, for example, [11, Theorem~4.1.8 and Chapter~8.9]), 
and so the operator $-\Delta_3$ with domain ${\frak D}_{\Gamma}$  
has a smallest positive eigenvalue, which is $\lambda_1(\Gamma)$.  
Preparatory to stating our new results, we define 
$$\Theta(\omega)=\sqrt{\max\left\{ 0 \,,\, 1-\lambda_1(\Gamma)\right\}} 
\in [0,1)\qquad\qquad 
\hbox{($0\neq \omega\in{\frak O}\,$ and $\,\Gamma =\Gamma_0(\omega)\leq SL(2,{\frak O})$)}  
\eqno(1.10)$$ 
and 
$$\vartheta =\sup_{0\neq \omega\in{\frak O}}\Theta(\omega)\;.\eqno(1.11)$$ 
Note that it is trivially the case that 
$$\lambda_1\left(\Gamma_0(\omega)\right)\geq 
1-\left(\Theta(\omega)\right)^2\geq 1-\vartheta^2\qquad\qquad 
\hbox{($0\neq \omega\in{\frak O}$).}\eqno(1.12)$$ 
Selberg's conjecture concerning the eigenvalues of $-\Delta_3$ for 
Hecke congruence subgroups of $SL(2,{\frak O})$ is that   
for all non-zero $\omega\in{\frak O}$ one has $\lambda_1(\Gamma_0(\omega))\geq 1$; 
an equivalent conjecture is that the constant $\vartheta$ is zero.  
Work of Kim and Shahidi, [29] and [30], has shown that 
$\Theta(\omega)<2/9$ for $0\neq \omega\in{\frak O}$, so that one has 
$$0\leq\vartheta\leq 2/9\eqno(1.13)$$ 
and, by (1.12), $\lambda_1(\Gamma_0(\omega))>77/81$ for $0\neq \omega\in{\frak O}$. 
\par 
We now state our new bounds for the mean value $E(D;M,A)$ 
that is defined in (1.3). 
\bigskip 

\goodbreak\proclaim{\Smallcaps Theorem 1}. Let $\vartheta$ be the real constant defined in (1.10) and (1.11). 
Let $\varepsilon > 0$. Then, 
for $D\geq 1$, $M\geq 1$ and all functions $A : {\frak O}-\{ 0\} \rightarrow {\Bbb C}$, 
one has both 
$$E(D;M,A)\ll_{\varepsilon} \left( D^{2+\varepsilon} 
+ \bigl( 1 + D M^{-3/2}\bigr)^{\vartheta} 
D^{1+\varepsilon} M^2\right)\sum_{0<|\mu|^2\leq M} \left| A(\mu )\right|^2\eqno(1.14)$$
and 
$$\eqalign{E(D;M,A) &\ll_{\varepsilon} D^{2+\varepsilon}
\biggl(\ \sum_{0<|\mu|^2\leq M} \left| A(\mu )\right|^2\biggr)\ +\cr 
 &\qquad +\left( 1 + D M^{-2}\right)^{\vartheta} 
D^{1+\varepsilon} M^3 \max\left\{\left| A(\mu )\right|^2\,:\, 
\mu\in{\frak O}\ {\rm and}\ 0<|\mu|^2\leq M\right\}\;.}\eqno(1.15)$$

\bigskip

In our proof of Theorem~1 we follow quite standard practice in 
utilising  an approximate functional equation for 
Hecke's zeta-functions $\zeta(s,\lambda^d)\,$ ($d\in{\Bbb Z}$). 
We remark that it might, perhaps, 
have been both interesting and profitable to have adopted 
a more novel approach, as Sarnak does in his proof 
(without the use of any approximate functional equation) of 
the sharp fourth power moment estimate 
$$\ \ \sum_{\!\!\!\!-D\leq d\leq D}\ \int\limits_{-D}^D \left|\zeta\left( 1/2 +it, \lambda^d\right)\right|^4
{\rm d}t
\ll D^2\log^4 D\qquad\qquad\hbox{($D\geq 2$),}\eqno(1.16)$$ 
which is [39, Theorem 1], 
and as Motohashi does in [35] 
(where he obtains an explicit formula for fourth moments of the Riemann zeta function); 
we however decided (at an early stage of this work) that we should   
simply concentrate on the approach in which we had the greatest confidence.   
We are similarly unadventurous in our use of the Poisson summation formula  
in proving Theorem~2 (below): 
we might instead have attempted, there, to emulate the approach  
taken by Blomer, Harcos and Michel in [3, Section~4.1],  
in which Jutila's variant of the circle method, from [25], is used  
to detect a condition of summation of the form $m_1 n_1 - m_2 n_2 =h$. 
\par 
A suitable approximate functional equation for $\zeta(s,\lambda^d)$ 
is proved in Section~3 of this paper: 
the proof we give is an implementation of a general method developed by 
Ivi\'{c}, in [21]. Although this approximate functional equation (which is Lemma~11, below) is 
very nearly contained in a more general theorem of Harcos 
[12, Theorem~2.5], it does have the merit of implying a 
slightly sharper bound for the relevant error term. 
\par 
By means of the approximate functional equation (Lemma~11), 
we show in Section~4 that Theorem~1 is 
a corollary of the asymptotic estimates contained in the results 
(1.30)-(1.32) of Theorem~2, below. Before coming to the statement of 
Theorem~2, we first define some of the notation used there (and subsequently). 
\par 
For $p\in [1,\infty]$, and for any function $b : X\rightarrow{\Bbb C}\,$ 
with domain $X\neq\emptyset$, we define  
$$\| b\|_p = 
\cases{\left(\sum\limits_{x\in X} |b(x)|^p\right)^{\!\!1/p}_{\matrix{\ }} &if $\,1\leq p<\infty\,$;\cr 
\sup\left\{ |b(x)| \,:\, x\in X\right\} &if $\,p=\infty\;$.}\eqno(1.17)$$ 
\par 
For $\alpha,\beta\in{\frak O}$, we shall (unless $\alpha =\beta =0$) 
take $(\alpha,\beta)$ to denote 
a highest common factor of $\alpha$ and $\beta$; at the same time 
$[\alpha,\beta]$ will denote a least common multiple of $\alpha$ and $\beta\,$ 
(so that if $\gamma =[\alpha,\beta]$ and $\delta =\alpha\beta /(\alpha,\beta)$ 
then $\gamma$ and $\delta$ are associates in the ring ${\frak O}$). 
\par 
Throughout this paper the notation ${\rm d}_{+}z$ denotes 
the standard Lebesgue measure on the set ${\Bbb C}\,$ 
(hence ${\rm d}_{+}z={\rm d}x\,{\rm d}y$, where $x$ and $y$ are the real and 
imaginary parts of the complex variable $z$). 

\bigskip 

\goodbreak\proclaim{\Smallcaps Theorem 2}. Let $\vartheta$ be the real constant defined in (1.10) and (1.11). 
Let $0<\varepsilon\leq 1/6$ and $0<\eta\leq(\log 2)/3$. 
Let $K_1 , K_2 , M_1\geq K_0=1$, and let $T>0$ satisfy 
$$T\gg\max\left\{ K_1^2 , K_2^2 , M_1^2\right\}\;.\eqno(1.18)$$
Let $w_0,w_1,w_2 : (0,\infty)\rightarrow{\Bbb C}$ be infinitely differentiable
functions satisfying, for $i\in\{ 0,1,2\}$, $j=0,1,2,\ldots\ $ and $x>0\;$,
$$w_i^{(j)}(x)=\cases{O_j\left( (\eta x)^{-j}\right) &if $e^{-\eta}K_i\leq x\leq e^{\eta}K_i\;$, \cr
0 &otherwise,}\eqno(1.19)$$
and let $a : {\frak O}-\{ 0\}\rightarrow{\Bbb C}$ be a function such that
$$|a(\mu)|>0\qquad{\rm only\  if}\qquad e^{-\eta}M_1\leq |\mu|^2 \leq e^{\eta}M_1\;.\eqno(1.20)$$
Define $C$ to be the mapping, from ${\frak O}-\{ 0\}$ to ${\Bbb C}$, given by  
$$C(\xi)= 
\sum\sum\sum_{\!\!\!\!\!\!\!\!\!\!\!\!\!\!\!\!\!\!\!\!\!\!{\scriptstyle 
\!\kappa_1\quad\!\!\kappa_2\quad\!\mu\atop\scriptstyle \kappa_1\kappa_2\mu=\xi}}
w_1\!\left(\left|\kappa_1\right|^2\right)
w_2\!\left(\left|\kappa_2\right|^2\right)
a(\mu)\qquad\qquad\hbox{($0\neq\xi\in{\frak O}$),}\eqno(1.21)$$ 
and put 
$$\eqalignno{ {\cal D}_0 &= 
\Biggl(\pi\int\limits_0^{\infty} w_0(x)\,{\rm d}x\Biggr)\,T\,\| C\|_2^2\;, &(1.22) \cr   
C^{\star}(\beta;z) &= 
\overline{w_1\left( |\beta z|^2\right)}
\,\sum\sum_{\!\!\!\!\!\!\!\!\!\!{\scriptstyle\mu\quad\,\kappa\atop\scriptstyle\mu\kappa =\beta}} 
\,a(\mu) w_2\left( |\kappa|^2\right) 
 &\hbox{($0\neq\beta\in{\frak O}$, $0\neq z\in{\Bbb C}$),}\qquad\qquad(1.23) \cr  
N &= 
T^{\varepsilon -1} K_1 K_2 M_1\;, &(1.24) \cr  
{\cal D}^{\star}_1 &= 
\Biggl(-4\pi\int\limits_0^{\infty} w_0(x)\,{\rm d}x\Biggr) 
\,T\int_{\Bbb C} 
\qquad\sum\!\!\!\!\!\!\!\sum_{\!\!\!\!\!\!\!\!\!\!\!\!\!\!\!\!{\scriptstyle 
\beta_1\neq 0\quad\beta_2\neq 0\atop\scriptstyle 
|(\beta_1,\beta_2)|^2<e^{-4\eta}T^{-\varepsilon}N}}\!\!C^{\star}\left(\beta_1;z\right) 
\overline{C^{\star}\left(\beta_2;z\right)}\,\left|\left(\beta_1,\beta_2\right)\right|^2 
{\rm d}_{+}z\;,\qquad\qquad\qquad &(1.25) \cr  
p(\alpha) &= 
\sum_{\nu\neq 0}|\nu|^{-4} \exp( 2\pi i \Re (\alpha\nu))  
 &\hbox{($\alpha\in{\Bbb C}$),}\qquad\qquad(1.26) \cr  
I_{w_0}(\gamma) &= 
{|\gamma |^4\over\pi^4}\int_{\Bbb C} 
\left(\left( x\,{{\rm d}^2\over{\rm d}x^2} + {{\rm d}\over{\rm d}x}\right)^{\!\!2} w_0(x)
\bigg|_{x=|s|^2}\right) 
p\!\left( {s\over\gamma}\right) 
{\rm d}_{+}s &\hbox{($0\neq\gamma\in{\Bbb C}$),}\qquad\qquad(1.27)}$$ 
and 
$$\qquad\ \,{\cal D}^{\star}_2 =4T\int_{\Bbb C} 
\qquad\ \sum\!\!\!\!\!\!\sum_{\!\!\!\!\!\!\!\!\!\!\!\!\!\!\!\!\!\!\!\!\!\!{\scriptstyle 
\beta_1\neq 0\quad\quad\!\beta_2\neq 0\atop\scriptstyle 
e^{-4\eta}T^{-\varepsilon}N\leq |(\beta_1,\beta_2)|^2<e^{\eta}N}} 
\!\!C^{\star}\left(\beta_1;z\right) 
\overline{C^{\star}\left(\beta_2;z\right)}\,\left|\left(\beta_1,\beta_2\right)\right|^2 
I_{w_0}\!\!\left( {\left[\beta_1,\beta_2\right] z\over T^{1/2}}\right) {\rm d}_{+}z\;.
\eqno(1.28)$$ 
For $d\in{\Bbb Z}$ and $t\in{\Bbb R}$, put 
$$c(d,it)= 
\left(\sum_{\kappa_1\neq 0}\!{w_1\bigl(\left|\kappa_1\right|^2\bigr)  
\Lambda^{d}(\kappa_1)\over\left|\kappa_1\right|^{2it}}\right)   
\left(\sum_{\kappa_2\neq 0}\!{w_2\bigl(\left|\kappa_2\right|^2\bigr) 
\Lambda^{d}(\kappa_2)\over\left|\kappa_2\right|^{2it}}\right)   
\left(\sum_{\mu\neq 0} 
{a(\mu)\Lambda^{d}(\mu)\over\left|\mu\right|^{2it}}\right) .\eqno(1.29)$$
Then one has  
$$\sum_{d=-\infty}^{\infty}\ \int\limits_{-\infty}^{\infty}
\left| c(d,it)\right|^2
w_0\!\left({|2d+it|^2\over\pi^2 T}\right) {\rm d}t
=2\pi{\cal D}_0 + (\pi /2)\left( {\cal D}^{\star}_1+{\cal D}^{\star}_2\right) 
+{\cal E}\;,\eqno(1.30)$$
where ${\cal E}$ is some complex number satisfying both   
$${{\cal E}\over K_1 K_2}\ll_{\eta,\varepsilon}\left( {M_1^2\over T^{1/2}}+ 
\left( {M_1^{2-(3/2)\vartheta}\over T^{(1-\vartheta)/2}}\right)
\!\left( {K_2\over T^{1/2}}\right)^{\!\!\vartheta}\right) 
T^{1+9\varepsilon} \| a\|_2^2\eqno(1.31)$$ 
and  
$${{\cal E}\over K_1 K_2}\ll_{\eta,\varepsilon}\left( {M_1^2\over T^{1/2}} +  
\left(\!\left( {K_1\over T^{1/2}}\right)^{\!\!\vartheta /2}
\!\!\left( {K_2\over T^{1/2}}\right)^{\!\!1-(\vartheta /2)} 
+\left( {K_2\over T^{1/2}M_1^{1/2}}\right)^{\!\!\vartheta}\right) 
\!\left( {M_1^2\over T^{1/2}}\right)^{\!\!1-\vartheta}\,\right) 
\!T^{1+11\varepsilon} M_1 \| a\|_{\infty}^2\;.\eqno(1.32)$$ 
The implicit constants in (1.31) and (1.32) are determined by those in (1.19) and (1.18), 
and by $\varepsilon$ and $\eta$. 

\bigskip

\goodbreak\proclaim{\Smallcaps Remarks~3}. For future reference, we observe here that the 
hypotheses (1.19), (1.20) and the definitions (1.22), (1.21) and (1.17) imply the 
upper bound 
$${\cal D}_0\ll 
T\sum_{|\kappa_1|^2\asymp K_1} 
\sum_{|\kappa_2|^2\asymp K_2}
\sum_{|\mu_1|^2\asymp M_1}
\sum_{|\kappa_3|^2\asymp K_1}
\sum_{|\kappa_4|^2\asymp K_2}
\sum_{\scriptstyle |\mu_2|^2\asymp M_1\atop\scriptstyle  
\kappa_3\kappa_4\mu_2 =\kappa_1\kappa_2\mu_1}
\left| a\left(\mu_1\right) a\left(\mu_2\right)\right|\;.$$
From this, the inequality 
$2|a(\mu_1)a(\mu_2)|\leq |a(\mu_1)|^2+|a(\mu_2)|^2$, and the bound 
(2.13) noted below, it follows that, when 
$K_1,K_2,M_1\geq 1$ and $T>0$ are such that (1.18) holds, one has: 
$$\eqalignno{ 
\left| {\cal D}_0\right| &\leq T 
\sum_{|\kappa|^2\asymp K_1} 
\sum_{|\lambda|^2\asymp K_2}
\sum_{|\mu|^2\asymp M_1} 
\left| a\left(\mu\right)\right|^2 
\sum_{\kappa'}
\sum_{\lambda'}
\sum_{\scriptstyle\mu'\atop\scriptstyle  
\!\!\!\!\!\!\!\!\!\!\!\!\!\!\!\!\!\!\!\!\!\!\!\!{\kappa'\lambda'\mu' =\kappa\lambda\mu}} O(1)  = \cr 
 &=T \sum_{|\kappa|^2\asymp K_1} 
\sum_{|\lambda|^2\asymp K_2}
\sum_{|\mu|^2\asymp M_1} 
\left| a\left(\mu\right)\right|^2 
O_{\varepsilon}\left( |\kappa\lambda\mu|^{\varepsilon}\right)   
\ll_{\varepsilon} T^{1+\varepsilon} K_1 K_2 \| a\|_2^2 
\qquad\qquad\hbox{($\varepsilon >0$).} &(1.33)}$$
Similar upper bounds may be obtained for the terms 
${\cal D}_1^{\star}$ and ${\cal D}_2^{\star}$ occurring, alongside ${\cal D}_0$,  
in the result (1.30). Indeed, by 
(1.23), (1.19) and the Cauchy-Schwarz inequality, it follows that
$$\left| C^{\star}(\beta ; z)\right|^2\ll \tau_2(\beta) 
\sum_{\scriptstyle\mu\mid\beta\atop\scriptstyle 
|\mu|^2\asymp M_1} |a(\mu)|^2\qquad\quad 
\hbox{for $\,0\neq\beta\in{\frak O}$, $\,z\in{\Bbb C}$}$$
(with $\tau_2(\beta)$ denoting, here and henceforth, the number of 
Gaussian integer divisors of $\beta$) 
and that it is, moreover, the case that $C^{\star}(\beta ; z)=0$ unless one has both
$|\beta|^2\asymp K_2 M_1$ and $|\beta z|^2\asymp K_1$. Hence, and by the
arithmetic-geometric mean inequality, one has, for some $z\in{\Bbb C}$, 
$$\eqalignno{
{\cal D}_1^{\star}
 &\ll T \sum_{\beta\neq 0} \left| C^{\star}(\beta ; z)\right|^2 
\!\!\!\sum_{\left|\beta'\right|^2\asymp K_2 M_1} 
\!\!\!\left|\left(\beta , \beta '\right)\right|^2 |\beta'|^{-2}K_1 \ll 
\qquad\qquad\qquad\qquad\qquad\qquad\qquad\qquad\qquad\qquad\cr 
 &\ll T K_1 \sum_{\beta\neq 0} \left| C^{\star}(\beta ; z)\right|^2 \tau_2(\beta) \ll\cr
 &\ll T K_1 \sum_{|\beta|^2\asymp K_2 M_1} 
\tau_2^2(\beta)\sum_{\scriptstyle\mu\mid\beta\atop\scriptstyle 
|\mu|^2\asymp M_1} |a(\mu)|^2 
\leq \Biggl( T K_1\sum_{|\kappa|^2\asymp K_2} {\tau_2^2(\kappa)\over 16}\Biggr) 
\Biggl(\;\sum_{|\mu|^2\asymp M_1} 
\tau_2^2(\mu) |a(\mu)|^2\Biggr)\,. &(1.34)}$$
Moreover, given the hypotheses (1.19), (1.20) and the definitions (1.23) and (1.24), 
it is implicit in (1.28) that 
$$\left|{\left[\beta_1,\beta_2\right] z\over\sqrt{T}}\right|^2
={\left|\beta_1\beta_2 z\right|^2\over\left|\left(\beta_1,\beta_2\right)\right|^2 T}
\ll{K_2 M_1 K_1\over e^{-4\eta}N T^{-\varepsilon +1}}
=e^{4\eta}\ll 1\;,$$
and so, as (1.26) and (1.27) trivially imply $I_{w_0}(\gamma)\ll |\gamma|^4$,
one may deduce from (1.23) and (1.28) essentially the same upper bound for $\bigl|{\cal D}_2^{\star}\bigr|$ 
as that found for $\bigl|{\cal D}_1^{\star}\bigr|$ in (1.34). Therefore it follows
by the bound (2.13) (below) that, when the hypotheses of Theorem~2  
are satisfied, the terms ${\cal D}_1^{\star}$ and ${\cal D}_2^{\star}$ in (1.30) 
satisfy 
$${\cal D}_j^{\star}\ll_{\varepsilon} T^{1+\varepsilon} K_1 K_2 \| a\|_2^2\qquad\quad 
\hbox{for $\,j=1,2\,$ and all $\,\varepsilon >0$.}\eqno(1.35)$$ 

\bigskip 

The details of our proof of Theorem~2 appear in Sections~5,~6 and~7 of this paper. 
\par 
The contents of Section~5 are several basic lemmas that are needed 
in Section~6: note that Lemma~17 is effectively the first step in 
our analysis of the sum occurring on the left hand side of Equation~(1.30). 
\par 
By means of the three (quite specialised) lemmas of Section~6, 
we transform the task of bounding the term ${\cal E}$ in Equation~(1.30) into 
a search for suitable bounds for the sum of Kloosterman sums occurring 
in the equation (6.47) of Lemma~22; aside from the use made of Lemma~17 
in the proof Lemma~20, it is fair to regard the proofs of 
Lemma~20, Lemma~21 and Lemma~22 as being exercises in the application of 
the Poisson Summation formulae that are contained in Lemma~14. 
\par 
In Section~7 we reformulate (as Lemma~23) the results of 
[44, Theorem~11], and we use most of the remainder of the section 
in applying those results so as to obtain (in Lemma~28) the sought for 
bounds on the sum of Kloosterman sums in (6.47); 
the supplementary bound (7.12) of Lemma~25 is needed in 
order to handle certain extreme cases (in which the condition (7.2) 
of Lemma~23 becomes an obstacle). At the end of Section~7 we complete 
our proof of Theorem~2, using only Lemma~20, Lemma~21, Lemma~22 and Lemma~28; 
our proof of Theorem~1 is, thereby, also completed (for we show in Section~4 that 
Theorem~2 implies Theorem~1). 
\par 
We conclude this introduction with a brief discussion of one 
likely application Theorem~1, followed by  
a remark in respect of one immediate implication of Theorem~1, 
and a remark on the possibility of  
generalising Theorem~1 in a non-trivial way.  
\par 
Our bound (1.7) for $S(T,N)$ played an essential part in 
work of Baker, Harman and Pintz on the distribution of rational 
primes; with its help, they showed, in [1], that there 
exists an $x_0\in(0,\infty)$ such that 
$$\left|\left\{ p\,: 
\,p\ {\rm is\ prime\ in}\ {\Bbb Z}\ {\rm and}\ x<p\leq x+x^{0.525}\right\}\right| 
\geq {9 x^{0.525}\over 100 \log x}\qquad\  
\hbox{for all $\,x\geq x_0$.}\eqno(1.36)$$ 
Similar progress on the distribution of Gaussian primes has 
been somewhat impeded by the the lack (until now) of a Gaussian integer 
analogue of (1.7). Nevertheless, in [16], Harman, Kumchev and Lewis 
have shown that there exist 
positive real numbers $c_1$ and $r_0$ such that, if one has $1\geq\alpha\geq 0.53$, then 
$$\left|\left\{ \pi_1\,: 
\,\pi_1\ {\rm is\ prime\ in}\ {\Bbb Z}[i]\ {\rm and} 
\ \left|\pi_1 - z\right|\leq |z|^{\alpha}\right\}\right| 
\geq {c_1 |z|^{2\alpha}\over\log |z|}\qquad\ 
\hbox{for all $z\in{\Bbb C}$ such that $|z|\geq r_0$.}\eqno(1.37)$$ 
Moreover, Lewis has improved on this by showing, in his thesis [33], 
that (1.37) holds (for some $c_1,r_0\in(0,\infty)$) whenever $1\geq\alpha >0.528$. 
Since our new estimate (1.15) implies the required Gaussian integer analogue of 
the bound (1.7), 
we expect that, by methods entirely analogous to those employed in the proof, in [1],  
of the result (1.36), it could now be proved that there exist 
positive real numbers $c_1$ and $r_0$ such that (1.37) 
holds whenever $1\geq\alpha\geq 0.525$. 
\par 
Our estimate (1.15) in Theorem~1 certainly does imply that, if    
$M\leq D^{1/2}$,  then one has 
$$\sum_{-D\leq d\leq D\,\hbox{\ }}\int\limits_{-D}^D \left|\zeta\left( 1/2 +it, \lambda^d\right)\right|^4 
\Biggl|\sum_{M/2<|\mu|^2\leq M} {A(\mu) \Lambda^d(\mu)\over |\mu|^{1+2it}}\Biggr|^2 
{\rm d}t \ll_{\varepsilon} D^{2+\varepsilon}\| A\|_{\infty}^2\qquad\quad   
\hbox{($\varepsilon >0$).}\eqno(1.38)$$ 
If it could be shown that (1.38) holds whenever $0\leq M\leq D$, 
then, by virtue of Corollary~12 (below), it would follow that 
$$\ \ \sum_{\!\!\!\!-D\leq d\leq D}\ \int\limits_{-D}^D \left|\zeta\left( 1/2 +it, \lambda^d\right)\right|^6  
{\rm d}t \ll_{\varepsilon} D^{2+\varepsilon} \qquad\   
\hbox{for all $\varepsilon >0$ and all $D\geq 1$.}$$ 
\par 
In [10, Theorem~2.2] it was shown by Duke that, if $F$ is a number field of 
degree $n$, if ${\frak q}$ is an ideal in ${\frak O}_F$ (the ring of 
integers of $F$), and if $\chi$ and $\{ \lambda_1,\ldots ,\lambda_{n-1}\}$ are
(respectively) a narrow class character $\bmod\,{\frak q}$ and a basis for 
the torsion free Hecke characters $\bmod\,{\frak q}$, then, for some $B\in(0,\infty)$, 
one has:  
$$\sum_{{\bf d}\in({\Bbb Z}\cap[-D,D])^{n-1}}\ \int\limits_{-D}^{D} 
\left|\zeta_{F}\!\left( 1/2+it ,\,\chi \lambda_{*}^{\bf d}\right) 
\right|^4 {\rm d}t\ \ll_{F,{\frak q}}\ D^n\log^B D\qquad\quad\ \hbox{($D\geq 2$),} 
\eqno(1.39)$$ 
where $\chi\lambda_{*}^{\bf d}$ is the groessencharacter 
$\chi\lambda_1^{d_1}\cdots\lambda_{n-1}^{d_{n-1}}$, while  
$\zeta_F(s,\,\chi\lambda_{*}^{\bf d})$ is the Hecke zeta function such that,  
for all $s_0\in{\Bbb C}$ such that $\Re (s_0)>1$, one has    
$\zeta_F(s_0,\,\chi\lambda_{*}^{\bf d}) 
=\sum_{{\frak a}\in I(F)}  
\,(\chi\lambda_{*}^{\bf d})({\frak a})\left( N({\frak a})\right)^{-s_0}$,  
with $I(F)$ being the set of non-zero 
ideals in ${\frak O}_F$, and with $N({\frak a})$ being the norm of ${\frak a}$. 
This result of Duke comes close to being 
a generalisation of Sarnak's bound in (1.16),  
falling short of that only by some power of $\log D$. 
\par\goodbreak  
In contrast with the considerable generality of the bounds in (1.39), 
our focus in the present paper 
is exclusively on Hecke zeta functions that are associated with ${\Bbb Q}(i)$.   
Similarly, in the papers [44] and [45] (upon which the present paper depends), 
we use only the case $F={\Bbb Q}(i)$ of the  
sum formulae [34, Theorem~11.3.3 and Theorem~12.3.2]. In [34] itself 
it is permitted that $F$ be any given imaginary quadratic field 
(while the relevant discrete group $\Gamma$ may be any Hecke congruence 
subgroup of $SL(2,{\frak O}_F)$). 
We expect that, 
by means of the sum formulae in [34] (or some slight extension thereof), 
it would be possible to generalise our work in [44,45] and the present 
paper so as to obtain new and useful upper bounds for 
sums of the form 
$$\ \ \sum_{\!\!\!\!-D\leq d\leq D\ }\ \int\limits_{-D}^{D} 
\left|\zeta_{F}\!\left( 1/2+it ,\,\chi \lambda_1^d\right) 
\right|^4 
\Biggl|\sum_{\scriptstyle {\frak a}\in I(F)\atop\scriptstyle N({\frak a})\leq M} 
{A({\frak a}) \left(\chi \lambda_1^d\right)\!({\frak a})\over  
\left( N({\frak a})\right)^{it}}\Biggr|^2
{\rm d}t\;,$$
where $F$ might be any imaginary quadratic field and  
$A$ might be any complex valued function defined on the set $I(F)$ of 
ideals in ${\frak O}_F$, while $\chi$, $\lambda_1$, 
$\zeta_F(s,\chi\lambda_1^d)$ and $N({\frak a})$ 
would all have the 
same meaning as in the above 
paragraph on Duke's result (1.39) (given that $F$ is now assumed 
to be a number field of degree $n=2$);  we conjecture that, by taking such an approach  
towards generalising Theorem~1, one might  
obtain mean value estimates for Hecke zeta functions 
capable of being used to obtain worthwhile improvements of certain 
results concerning the distribution of prime  ideals in 
imaginary quadratic fields ([16, Theorem~2], in particular). 

\medskip 

\goodbreak 
\noindent{\bf Notation and terminology that is fairly standard} 

\smallskip 

{\settabs\+\quad &$E_j^{\frak a}(q,P,K;N,{\bf b})\,$\  &--\quad &\cr 
\+&${\Bbb N}$ &--&the set $\{ n\in{\Bbb Z} : n\geq 1\}$;\cr  
\+&$|{\cal A}|$ &--&the cardinality of the set ${\cal A}\,$ 
(so that $|\{ x\in{\Bbb R} : x^2=1\}|=2$, for example); \cr 
\+&$\max{\cal X}$   &--&the greatest element of the set ${\cal X}\subset{\Bbb R}$ 
(where this exists);   \cr 
\+&$\min{\cal X}$   &--&the least element of ${\cal X}\subset{\Bbb R}$ (where this exists);   \cr 
\+&$\max_{A(x)} f(x)$   &--&equal to $\max\{ f(x) : A(x)\ {\rm is\ true}\}$,  
when $A(x)$ is some statement about $x$;  \cr 
\+&$\min_{A(x)} f(x)$   &--&equal to $\min\{ f(x) : A(x)\ {\rm is\ true}\}$,  
when $A(x)$ is some statement about $x$;  \cr 
\+&$\Re (z)$ and $\Im (z)$ &--&the real and
imaginary parts of the complex number $z$;  \cr 
\+&$|z|$ and $\overline{z}$   &--&the absolute value and complex conjugate 
of the complex number $z$;    \cr 
\+&${\rm Arg}(z)\in(-\pi,\pi]$  &--&the principal value of the argument of 
the non-zero complex number $z$;    \cr 
\+&$g\circ f$   &--&the function obtained by composing $f$ with $g\,$  
(so that $(g\circ f)(x)=g(f(x)));$   \cr 
\+&$f^{(j)}(s_0)$   &--&the 
$j$-th derivative of the function $f$ at the point $s=s_0\,$ (where this exists);   \cr 
\+&$f^{(0)}(s)$   &--&equal to $f(s)$;   \cr 
\+&$\pi$   &--&the ratio of the circumference of a circle to its diameter;    \cr 
\+&$e$   &--&the base of the natural logarithm function, $\ln(x)$;    \cr 
\+&$\log(z)$   &--&the principal value of the logarithm of $z$, equal to 
$\ln(|z|)+i{\rm Arg}(z)$;    \cr 
\+&$\exp(z)$   &--&equal to $e^z$, when $z\in{\Bbb C}$;   \cr 
\+&$i$   &--&most often denotes a square root of $-1$,  but sometimes  
is an integer variable;   \cr 
\+&$\Gamma(z)\ \ {\ \atop {\ \atop\ }}$   &--&Euler's Gamma function;   \cr 
\+&${\Gamma^{(1)}\over\Gamma}(z){{\ \atop\ } \atop\ }$   &--&the logarithmic derivative of $\Gamma(z)$, equal to 
both $\Gamma^{(1)}(z)/\Gamma(z)$ and ${{\rm d}\over{\rm d}z}\log\Gamma(z)$;   \cr 
\+&$\gamma$   &--&may denote either a variable or Euler's constant, 
$-\Gamma^{(1)}(1)=0.5772157\ldots\ $;   \cr 
\+&$n!$   &--&is `$n$-factorial' (equal to $\Gamma(n+1)$), when $n\in{\Bbb N}\cup\{ 0\}$;   \cr 
\+&$\bigl( {\scriptstyle n\atop\scriptstyle r}\bigr)$   &--&the `binomial coefficient' 
$(n!/(n-r)!)/r!$, when $n,r\in{\Bbb Z}$ satisfy $0\leq r\leq n$; \cr 
\+&$[x]$   &--&equal to $\max\{ m\in{\Bbb Z} : m\leq x\}$, 
when $x$ is a real number. \cr  
\+&${\bf u}\cdot{\bf v}$   &--&equal to $\sum_{j=1}^n u_j v_j$, 
when ${\bf u},{\bf v}\in{\Bbb C}^n$. \cr}  

\medskip 

\goodbreak 
\noindent $SL(2,R)$ denotes the group (under multiplication) of   
the $2\times 2$ matrices $M$ that have their elements in 
the integral domain $R$, and their 
determinants equal to $1$.

\smallskip 

\goodbreak 
\noindent When $D$ is an open subset of $\,{\Bbb C}$, a function $f : D\rightarrow{\Bbb C}$ 
may be termed `smooth' if and only if it is the case that, for all $j,k\in{\Bbb N}\cup\{ 0\}$, 
the partial derivative 
$({\partial^{j+k}/\partial x^j \partial y^k})f(x+iy)$ 
is defined and continuous at all points $(x,y)\in{\Bbb R}^2$ such that 
$x+iy\in D$.  

\smallskip 

\noindent The `Schwartz space' contains $F : {\Bbb R}^2\rightarrow{\Bbb C}$  if and only if,  
for all real $B\geq 0$ and all $j,k\in{\Bbb N}\cup\{ 0\}$, 
the partial derivative $({\partial^{j+k}/\partial x^j \partial y^k})F(x,y)$ 
is defined and continuous at all points $(x,y)\in{\Bbb R}^2$, and the mapping 
$(x,y)\mapsto (x^2+y^2)^{B}({\partial^{j+k}/\partial x^j \partial y^k})F(x,y)$ 
is a bounded function on ${\Bbb R}^2$. 

\smallskip 

\goodbreak 
\noindent{\bf Algebraic and number-theoretic notation} 

\smallskip 

{\settabs\+\quad &$E_j^{\frak a}(q,P,K;N,{\bf b})\,$\  &--\quad &\cr 
\+&$R^{*}$  &--&the group of units in $R$, when $R$ is a ring with an identity.       \cr  
\+&${\frak O}_F$  &--&the ring of integers of $F$, 
when $F$ is a number field; \cr  
\+&${\frak O}$  &--&the integral domain 
${\frak O}_{{\Bbb Q}(i)}={\Bbb Z}[i]=\{ m+n\sqrt{-1} : m,n\in{\Bbb Z}\}$;   \cr  
\+&Gaussian integer    &--&a number in the ring ${\frak O}$;       \cr  
\+&Gaussian prime  &--&a prime in the ring ${\frak O}$;    \cr  
\+&$(\alpha)$ or $\alpha{\frak O}$   &--&the ideal $\{\alpha\beta : \beta\in{\frak O}\}$ 
of ${\frak O}$, when $\alpha$ is a Gaussian integer;  \cr  
\+&$(\alpha)\mid(\delta)$    &--&signifies that the ideal $(\alpha)$ 
divides the ideal $(\delta)\,$  
(i.e. that ${\frak O}\supseteq \alpha{\frak O}\supseteq \delta{\frak O}$);  \cr  
\+&$\beta\equiv\alpha\bmod\gamma{\frak O}$    &--&signifies that 
$\beta$ is `congruent' to $\alpha$ $\bmod\ \gamma{\frak O}\,$ 
(i.e. that $(\gamma)\mid(\beta -\alpha)$);  \cr  
\+&${\frak O}/(\gamma{\frak O})$    &--&the ring of residue classes 
$\bmod\,\gamma{\frak O}\,$ (these being the cosets   
of $\gamma{\frak O}$ in ${\frak O}$);     \cr  
\+&$\Gamma_0(\omega)$ or $\Gamma$  &--&a Hecke congruence subgroup of 
$SL(2,{\frak O})\,$ (see Equation~(1.9) for the definition);  \cr 
\+&$\delta\mid\gamma$  &--&signifies that $\delta$ is a divisor of $\gamma\,$ 
(i.e. that one has both $\delta\in{\frak O}$ and $\gamma\in\delta{\frak O}$);  \cr 
\+&associate  &--&$\gamma\in{\frak O}$ and $\delta\in{\frak O}$ if and only if 
it one has both $\gamma\mid\delta$ and $\delta\mid\gamma$;  \cr 
\+&$\gamma\sim\delta$  &--&signifies that $\gamma,\delta\in{\frak O}$ 
are associates, so that one has $\gamma /\delta\in{\frak O}^{*}=\{ i,-1,-i,1\}$;  \cr 
\+&$\tau_n(\alpha)$    &--&the number of elements in 
the set $\{\;(\delta_1,\ldots ,\delta_n)\in{\frak O}^n : 
\,\delta_1\delta_2\cdots\delta_n=\alpha\}$;  \cr 
}

\smallskip 

\goodbreak 
\noindent $(\alpha , \beta)$  denotes a highest common factor of the 
Gaussian integers $\alpha$ and $\beta$; 
this notation is somewhat ambiguous, for if $\delta$ is a highest common factor of 
$\alpha$ and $\beta$, then so too are the three other associates of $\delta$ 
(i.e. $i\delta$, $-\delta$ and $-i\delta$); it does not, however, 
lead to any serious difficulties, since 
relations of the form $(\alpha ,\beta)\sim \delta$, or 
$|(\alpha ,\beta)|^2=n$, remain valid if 
the number $(\alpha ,\beta)$ is replaced by any one of its associates. 
In the case of ${\Bbb Z}$-valued variables or constants, $m$ and $n$ (say), 
we unambiguously put  
$(m , n)=\max\{ d\in{\Bbb N} : d\mid m\ {\rm and}\ d\mid n\}$.  

\smallskip 

\goodbreak 
\noindent $[\alpha,\beta]$ denotes 
a least common multiple of the Gaussian integers $\alpha$ and $\beta$;   
like the notation for highest common factors, this notation for 
least common multiples is harmlessly ambiguous. 

\smallskip

\goodbreak 
\noindent{\bf Specialised or customised notation} 

\smallskip  

{\settabs\+\quad &$E_j^{\frak a}(q,P,K;N,{\bf b})\,$\  &--\quad &\cr 
\+&$\zeta(s)$  &--&the Riemann zeta-function;  \cr  
\+&$L(s,\chi)$  &--&the Dirichlet $L$-function 
associated with the Dirichlet character $\chi$;  \cr  
\+&$\Lambda^d$  &--&$(\Lambda^d)_{d\in{\Bbb N}}$ 
is the family of endomorphisms  
of $({\Bbb Q}(i))^{*}$ that are given by (1.2);  \cr  
\+&$\lambda^d$  &--&the groessencharacter (or Hecke character) 
defined, for $d\in{\Bbb Z}$, just below (1.2);  \cr  
\+&$\zeta(s,\lambda^d)$  &--&a Hecke zeta function 
(meromorphic on ${\Bbb C}$ and, for $\Re (s)>1$, as stated in (1.1));  \cr  
\+&$P_M(A;s,\lambda^d)$  &--&given by (1.4), 
if $M\in(0,\infty)$, $A$ is a mapping from ${\frak O}-\{ 0\}$ to ${\Bbb C}$, 
$s\in{\Bbb C}$, $d\in{\Bbb Z}$;  \cr  
\+&$E(D;M,A)$  &--&the mean value defined in (1.1)-(1.4);  \cr  
\+&${\Bbb H}_3$    &--&a model for $3$-dimensional hyperbolic space (see 
below (1.8), and [11, Chapter~1]);  \cr  
\+&${\frak D}_{\Gamma}$  &--&a certain subspace of $L^2(\Gamma\backslash{\Bbb H}_3)\,$ 
(see below (1.8), and refer to [11, Chapters 1-4]);  \cr  
\+&$\Delta_3$  &--&the Laplacian operator, 
from ${\frak D}_{\Gamma}$ into $L^2(\Gamma\backslash{\Bbb H}_3)$, 
that is given by (1.8);  \cr  
\+&$\lambda_1(\Gamma)$  &--&the smallest positive eigenvalue of the operator $-\Delta_3$;  \cr  
\+&$\Theta(\omega)$   &--&the function, from ${\frak O}-\{ 0\}$ into $[0,2/9)$, that 
is defined in (1.10);  \cr  
\+&$\vartheta$  &--&the absolute constant defined in (1.11) 
(see also (1.12) and (1.13));  \cr  
\+&$\| b\|_p$    &--&defined in (1.17), it is, when $p\in [1,\infty]$, 
a norm of the complex valued function $b$;  \cr  
\+&${\rm d}_{+}z$    &--&a Lebesgue measure on ${\Bbb C}\,$ 
(${\rm d}_{+}z ={\rm d}x\,{\rm d}y$, where $x=\Re (z)$,  
$y=\Im (z)$);  \cr  
\+&$\int_{\Bbb C} f(z) {\rm d}_{+}z$   &--&
the integral of $f$, with respect to  ${\rm d}_{+}z$, 
over ${\Bbb C}\,$  
(equal to $\int_{-\infty}^{\infty}\int_{-\infty}^{\infty}f(x+iy) {\rm d}x\,{\rm d}y$); \cr  
\+&$X_d(s)$   &--&when $d\in{\Bbb Z}$, the mapping 
$s\mapsto X_d(s)$ is the meromorphic function given by (2.6);  \cr  
\+&$T(d,t)$    &--&when $d\in{\Bbb Z}$, the mapping 
$t\mapsto T(d,t)$ is the real function given by (2.6) and (2.8);  \cr  
\+&$\delta(\Lambda^d , n)$   &--&when $d\in{\Bbb Z}$, the real sequence 
$(\delta(\Lambda^d , n))_{n\in{\Bbb N}}$ is given by (2.10) and (1.2);  \cr  
\+&$z_{m,t}$ and $Z_{m,t}$  &--&the complex number in (3.2) and   
(in the proof of Lemma~7) its absolute value;  \cr  
\+&$\delta_{a,b}$  &--&is equal to $1$ if $a=b$, and is otherwise equal to zero 
(see (3.40));  \cr  
\+&${\rm e}(x)$    &--&denotes $\exp(2\pi i x)$, 
when $x\in{\Bbb R}$;  \cr  
\+&$\hat F : {\Bbb R}\times{\Bbb R}\rightarrow{\Bbb C}$&--
&defined as in Lemma~14, it is a Fourier transform of 
the function $F : {\Bbb R}\times{\Bbb R}\rightarrow{\Bbb C}$;  \cr  
\+&$\hat f : {\Bbb C}\rightarrow{\Bbb C}$    &-- 
&defined as in Lemma~14, it is a Fourier transform of the 
function $f : {\Bbb C}\rightarrow{\Bbb C}$;  \cr  
\+&$\Delta_{{\Bbb R}\times{\Bbb R}}$ and $\Delta_{\Bbb C}$  &-- 
&defined in Lemma~15, these are two Laplacian operators (of the Euclidean type);  \cr  
\+&$S_{\alpha}$  &--&defined in Remarks~16, it is 
(when $\alpha\in{\Bbb C}^{*}$) a certain `rotation-dilatation operator';  \cr  
\+&${\frak N}$  &--&the function ${\frak N} : {\Bbb C}\rightarrow[0,\infty)$ 
such that one has ${\frak N}(z)=|z|^2$ for all $z\in{\Bbb C}$;  \cr  
\+&$S(\alpha,\beta;\gamma)$             &--&defined in (5.5), a Kloosterman sum associated 
with the number field ${\Bbb Q}(i)$; \cr } 

\smallskip

\goodbreak 
\noindent{\bf Superscripts and subscripts}\quad  
\noindent The superscripts $\star$, $\flat$ and $\sharp$, the dash, $'$, and double-dash, $''$,
and the tilde and breve accents, $\ \widetilde{}\ $ and $\ \breve{}\ $, have no intrinsic 
meaning (we simply found it convenient to use them in devising names 
for certain variables);  
when used as subscripts the symbols $*$ and $\dagger$ are similarly devoid of meaning.   

\smallskip 

\goodbreak 
\noindent In relations such as $\alpha \delta^{*}\equiv\kappa\bmod\ \gamma{\frak O}$, 
and in expressions such as ${\rm e}(\Re (\alpha \delta^*/\gamma))$, 
or the highest common factor $(\alpha \delta^{*} , \gamma)$, 
it is to be understood that $\delta^*$ denotes a variable 
dependent upon a variable $\delta$, and that $\delta$ and $\delta^{*}$ 
take values in ${\frak O}$ such that 
$\delta^* \delta\equiv 1\bmod \gamma{\frak O}$; where such notation is used, 
the relevant ${\frak O}$-valued variables 
$\gamma$, $\delta$ and $\delta^{*}$ are 
implicitly constrained (by virtue of the congruence condition just mentioned) to take only 
combinations of 
values satisfying both $(\delta,\gamma)\sim 1$ and $(\delta^{*},\gamma)\sim 1$. 
One notable use of this `$*$-notation' is in our definition, in (5.5), of 
the Kloosterman sum $S(\alpha,\beta;\gamma)$: note, in particular, that 
the summand in (5.5) may be expressed as the product 
${\rm e}(\Re (\alpha\delta^{*}/\gamma))\,{\rm e}(\Re (\beta\delta /\gamma))$, 
and note also the indication given above as to how the relation 
$\delta\in({\frak O}/\gamma{\frak O})^{*}$ is to be interpreted when it occurs 
as a condition of summation. 

\smallskip 

\goodbreak 
\noindent{\bf Conventions concerning certain notation associated with sums and products}\quad 
With regard to products in which the factors are indexed by a  
variable non-zero ideal $(\alpha)$ 
(and in which the ideal itself determines the value of the corresponding factor), 
it should be understood that there is a
one-to-one correspondence between the relevant set of ideals and the 
factors in the product: this effectively means that one may view  
the factors in the product as being indexed by the Gaussian integer $\alpha$ 
that is implicit in the notation `$(\alpha)$', provided only that 
one constrains this $\alpha$  to satisfy the conditions 
$\Im (\alpha)\geq 0$ and $\Re (\alpha)>0$. 

\smallskip 

\goodbreak 
\noindent When $\gamma\in{\frak O}$ and $\gamma\neq 0$, the elements of 
the group of units $({\frak O}/(\gamma{\frak O}))^{*}$ 
are the `reduced residue classes $\bmod\ \gamma{\frak O}$' (and 
each of these is, of course, a coset of $\gamma{\frak O}$ in ${\frak O}$). 
However, where there occurs a condition of summation of the form 
$\delta\in ({\frak O}/(\gamma{\frak O}))^{*}$, 
the variable of summation $\delta$ should, by an abuse of notation,  
be understood to have a range that is a fixed 
set of representatives 
$\{ \delta_1,\delta_2\ldots 
,\delta_{|({\frak O}/(\gamma{\frak O}))^{*}|}\}\subset{\frak O}$ 
of the reduced residue classes $\bmod\ \gamma{\frak O}$. 
Similarly, any condition of summation of the form $\kappa\in{\frak O}/\gamma{\frak O}$ 
signifies that $\kappa$ has a range that is a fixed set of 
representatives 
$\{ \kappa_1,\kappa_2\ldots 
,\kappa_{|{\frak O}/(\gamma{\frak O})|}\}\subset{\frak O}$ 
of the residue classes $\bmod\ \gamma{\frak O}$. 

\smallskip 

\goodbreak 
\noindent Variables of summation designated by Latin letters are ${\Bbb Z}$-valued variables, 
whereas those designated by Greek letters are instead ${\frak O}$-valued variables (and may 
take any value in ${\frak O}$ permitted by the conditions attached to the summation). 

\smallskip 

\goodbreak 
\noindent Where there occur nested summations, such as sums of the form 
$\sum_{\alpha}f(\alpha)\sum_{\beta} g(\alpha,\beta)\,$ (for example), 
it should be understood that the variable $\alpha$ of the outer summation 
is constrained to satisfy the condition $f(\alpha)\neq 0$.

\smallskip 

\goodbreak 
\noindent{\bf Notation for bounds and asymptotic estimates}\quad  
The notation $O_{\alpha_1 ,\ldots,\alpha_n}(B)$, when it is the $k$-th instance 
of $O$-notation used in this paper, should be understood 
to denote a complex-valued variable $\iota_k$  
that satisfies a condition of the form $|\iota_k|\leq C(\alpha_1,\ldots,\alpha_n) B$, in which  
the `implicit constant' 
$C(\alpha_1,\ldots,\alpha_n)$ is positive and 
depends only on $\alpha_1,\ldots,\alpha_n$ and declared constants.   
Where this `$O$-notation' is used, the relevant variable or constant 
$B$ must necessarily satisfy $B\geq 0$.  

\smallskip  

\goodbreak 
\noindent We frequently employ Vinogradov's notation, as an alternative to 
the $O$-notation. Thus we may use either the notation 
$\xi\ll_{\alpha_1 ,\ldots,\alpha_n} B$, or the equivalent 
notation $B\gg_{\alpha_1 ,\ldots,\alpha_n} \xi$, to signify that 
one has $\xi =O_{\alpha_1 ,\ldots,\alpha_n}(B)$. 
Where $A\geq 0$ and $B\geq 0$, the notation 
$A\asymp_{\alpha_1 ,\ldots,\alpha_n}B$ signifies that  
$A\ll_{\alpha_1 ,\ldots,\alpha_n}B\ll_{\alpha_1 ,\ldots,\alpha_n}A\,$ 
(i.e. that one has both $A\ll_{\alpha_1 ,\ldots,\alpha_n}B$ 
and $B\ll_{\alpha_1 ,\ldots,\alpha_n}A$). 
There are a few places where, instead of attaching subscripts (to the $O$, $\ll$, $\gg$ or 
$\asymp$ sign), we explicitly state the parameters upon which 
the relevant implicit constant, or constants, may depend. 

\smallskip 

\goodbreak 
\noindent{\bf Epsilon and Eta}\quad   
In the stating the results of this paper we treat $\varepsilon$ and $\eta$ as 
positive valued variables. However, in any meaningful application of 
our results in which $\varepsilon$ had a part to play, 
it would be necessary either to assign to $\varepsilon$  a specific,   
fairly small, numerical value   
(such as $10^{-10}$, say), or else to have $\varepsilon$ 
function as an `arbitrarily small positive constant'; 
the same applies in respect of the variable $\eta$. 

\bigskip 

\goodbreak 
\noindent{\SectionHeadingFont 2. Essential properties of the zeta functions}

\medskip

\noindent By work of Hecke [18, Pages 34-35] it is known that if $d\in{\Bbb Z}-\{ 0\}$, 
then the function
$$\xi\left( s, \lambda^d\right) =\Gamma(s+2|d|)\pi^{-s}\zeta\left( s, \lambda^d\right)\eqno(2.1)$$
has an entire analytic continuation satisfying, for all $s\in{\Bbb C}$,
$$\xi\left( s, \lambda^d\right) =\xi\left( 1-s, \lambda^{-d}\right).\eqno(2.2)$$
By substituting $\overline{\alpha}$ for $\alpha$ in (1.1) one finds that
$$\zeta\left( s, \lambda^d\right) = \zeta\left( s, \lambda^{-d}\right)\quad\hbox{and}\quad 
\xi\left( s, \lambda^d\right) = \xi\left( s, \lambda^{-d}\right).\eqno(2.3)$$

The case $d=0$ of (2.1) defines a function $\xi\bigl( s, \lambda^0\bigr)$ with
a single valued analytic continuation to ${\Bbb C}-\{ 0,1\}$ satisfying
the case $d=0$ of (2.2). This function has simple poles at $s=0$ and $s=1$, 
where the residues are $-1/4$ and $1/4$, respectively. These facts, combined with 
(2.1) and (2.2), show that  $\zeta\bigl( 0 , \lambda^0\bigr) =-1/4$, so the pole of $\xi\bigl( s , \lambda^0\bigr)$ 
at $s=0$ is purely an effect of the `Gamma factor' in (2.1).
If $\Re (s)>1$, then it follows by (1.1), (1.2) and (2.3) that
$\overline{\zeta\bigl( s , \lambda^d\bigr)}=\zeta\bigl( \overline{s} , \lambda^{-d}\bigr)
=\zeta\bigl( \overline{s} , \lambda^d\bigr)$.
Therefore, by Schwarz's reflection principle,
$$\zeta\left(\overline{s} , \lambda^d\right) =\overline{\zeta\left( s , \lambda^d\right)}
\quad\hbox{for $s\in{\Bbb C}-\{ 1\}\;$.}\eqno(2.4)$$

To summarise some of the above, one may note firstly that
if $(0,1)\neq (d,s)\in{\Bbb Z}\times{\Bbb C}$, then
$$\zeta\left( s, \lambda^d\right) =X_d(s)\zeta\left( 1-s, \lambda^{-d}\right),\eqno(2.5)$$
where
$$X_d(s)={\pi^{2s-1}\Gamma(2|d|+1-s)\over\Gamma(2|d|+s)}\,,\eqno(2.6)$$
and secondly that
$$\pi^{-1} \zeta\left( s, \lambda^0\right) = 
{1/4\over s-1} + h_0 + O_{\rho}\left( |s-1|\right) \quad\hbox{for $\rho\geq |s-1|>0$,}\eqno(2.7)$$
where $h_0$ is some real absolute constant. Introducing further new notation, 
define $T(d,t)$ by
$$\log T(d,t)=-{X_d^{(1)}\left(\textstyle{1\over 2}+it\right)\over 
X_d\left(\textstyle{1\over 2}+it\right)} 
\quad\hbox{for $d\in{\Bbb Z}$ and $t\in{\Bbb R}\,$.}\eqno(2.8)$$
Due to the central part it plays in Lemma~11 below 
(an approximate functional equation for $\zeta\bigl( s,\lambda^d\bigr)\/$)
the number $T(d,\Im (s))$ figures significantly in much of what follows.

By the absolute convergence of the sums in (1.1), 
$$\zeta\left( s, \lambda^d\right)
=\sum_{n=1}^{\infty} {\delta\left(\Lambda^d,n\right)\over n^s}\quad\hbox{for $d\in{\Bbb Z}$ and 
$\Re (s)>1$,}\eqno(2.9)$$
where
$$\delta\left(\Lambda^d,n\right)
={1\over 4}\sum_{|\alpha|^2=n}
\Lambda^d (\alpha)\eqno(2.10)$$
(so that, by virtue of the definition (1.2), one has $\delta(\Lambda^d ,n)\in{\Bbb R}$). 
By (2.10), (1.2) and a theorem of Jacobi [14, Theorem 278] one has
$\delta\bigl(\Lambda^0,n\bigr) =\sum_{k\mid n}\chi_4(k)$
for $n\in{\Bbb N}$, where $\chi_4(n)$ is the non-principal Dirichlet character $\bmod\ 4$.
Therefore
$$\zeta\left( s, \lambda^0\right) = \zeta(s)L\left( s, \chi_4\right)\eqno(2.11)$$
and, for $\varepsilon >0$, $n\in{\Bbb N}$ and $d\in{\Bbb Z}$,
$$\left| \delta\left(\Lambda^d,n\right)\right|
\leq{1\over 4}\sum_{|\alpha|^2=n}
\left| \Lambda^d (\alpha)\right|
= \delta\left(\Lambda^0,n\right)
=\sum_{k\mid n}\chi_4(k)\leq\sum_{k\mid n} 1
=O_{\varepsilon}\left( n^{\varepsilon}\right) .\eqno(2.12)$$
By using (2.11) and information about $\zeta(s)$ and $L(s,\chi_4)$
(such as [42, Equation (2.1.16)] and 
[6, Equation (15) of Chapter 6]) one can determine that
(2.7) holds with $h_0=(1/4)\bigl(\gamma +{L^{(1)}\over L}\bigl( 1,\chi_4\bigr)\bigr)$, 
where $\gamma$ is Euler's constant.  
Note also that, by its definition in the section on notation, the `divisor function' $\tau_j(\alpha)$  
satisfies, for $j\in{\Bbb N}$ and $\alpha\in{\frak O}-\{ 0\}$, both 
$4\leq\tau_j(\alpha)\leq\left(\tau_2(\alpha)\right)^{j-1}$ and 
$\tau_2(\alpha)\leq\sum_{n\mid |\alpha|^2}4\delta\left(\Lambda^0 , n\right)$, 
so that by (2.12) one has 
$$\tau_j(\alpha)\ll_{\varepsilon ,j}\ |\alpha|^{2(j-1)\varepsilon}\quad 
\hbox{for $j\in{\Bbb N}$, $\varepsilon >0$ and $\alpha\in{\frak O}-\{ 0\}\;$.}\eqno(2.13)$$

As noted in [18, Equation (30)], the multiplicativity of the groessencharacters
and absolute convergence of the series in (1.1) together imply the Euler product formula:
$$\zeta\left( s , \lambda^d\right)
=\prod_{\scriptstyle{\frak p}\in I\atop\scriptstyle {\frak p}\ {\rm is\ prime}} 
{1\over 1-\lambda^{d}({\frak p}) (N{\frak p})^{-s}}
\qquad\quad\hbox{for $\Re (s)>1\;$,}\eqno(2.14)$$
where the product is taken over all prime ideals ${\frak p}$ in ${\frak O}$, 
while the other notation used is as defined in, and just below, Equation~(1.2).
This Euler product representation of $\zeta\bigl( s , \lambda^d\bigr)$
is of fundamental significance in respect of its applications 
concerning the distribution of Gaussian primes.
Here however it will serve merely to justify the
assertion that there is no integer $d$ for which  
the function $s\mapsto\zeta\bigl( s , \lambda^d\bigr)$ 
is identically equal to zero
(the convergence of the Euler product  guaranteeing that
$\zeta\bigl( s , \lambda^d\bigr)\neq 0$
in the half plane where $\Re (s)>0$\/).

\bigskip

\goodbreak 
\noindent{\SectionHeadingFont 3. An approximate functional equation}

\medskip

\noindent The principal result in this section, Lemma~11, 
is essentially a special case of the very general
class of approximate functional equations given by Ivi\'c's theorem in [21],
though additional work has been done to make the result more uniform in the
groessencharacter aspect. This uniformity issue has previously been addressed
by Harcos in [12, Theorem 2.5], which is an
application of the methods of [21] in the context of
$L$ functions corresponding to irreducible cuspidal automorphic representations
of general linear groups over number fields. For errata to [12] see [13]. 
\par 
Lemma~11 is applicable at any point in the 
strip $-1/3\leq \Re (s)\leq 4/3$, whereas 
[12, Theorem~2.5] applies only at points on the critical line $\Re (s)=1/2$.   
However, the case $\sigma =1/2$ is the only case of Lemma~11 that is used 
in this paper. In light of [12, Remark 2.7] it is clear that the case $\sigma =1/2$ of 
Lemma~11 is very nearly a direct corollary of [12, Theorem 2.5] 
(it surpasses what [12, Theorem~2.5] would 
imply only to the extent of providing, in (3.50)-(3.52),  bounds for 
the relevant error term that are better than the bound stated in [12, Equation~(2.12)]). 
The above mentioned results of Harcos would  
meet our needs in this paper, but we think that Lemma~11 has enough 
merits of its own to justify the inclusion of its proof in this section    
(which also serves to make this paper more self-contained). 
\par 
One can find, in the literature, other approximate functional equations 
for Hecke zeta functions with groessencharacters: one is 
Huxley's result [19, Theorem~2], which  
is applicable to Hecke zeta functions associated with number fields of arbitrary 
degree, but appears formidably complicated; another is Lavrik's theorem in [32], 
which is reasonably simple and involves only very small error terms, but is   
less convenient (for our purposes) than Lemma~11. 
See also [24, Theorem~5.3] for a very general approximate functional equation.

\bigskip 

\goodbreak
\proclaim{\Smallcaps Lemma~4}. Let $0<\delta <\pi$. Then, for $|{\rm Arg}(z)|\leq\pi -\delta$,
$${\Gamma^{(1)}\over\Gamma}\,(z)
=\log z -{1\over 2z}+O_{\delta}\left( {1\over |z|^2}\right) .\eqno(3.1)$$

\goodbreak
\noindent{\Smallcaps Proof.}\ 
See [20, Page 57], where it is noted that methods of complex analysis 
enable one to deduce the result (3.1) from the 
asymptotic estimate 
$$\log\Gamma(z)=\left( z-{1\over 2}\right)\log z -z +{1\over 2}\log(2\pi)
+O_{\delta}\left({1\over |z|}\right)\qquad\quad 
\hbox{($z\in{\Bbb C}-\{ 0\}$, $\,|{\rm Arg}(z)|\leq\pi -\delta<\pi$).}$$
A proof of the latter estimate is given in [46, Section 13.6]\quad$\square$

\bigskip 

\goodbreak\proclaim{\Smallcaps Lemma~5}.  The equations
(2.8) and (2.6) together define a function $T:{\Bbb Z}\times{\Bbb R}\rightarrow (0,\infty)$
such that, for $t,t_1\in{\Bbb R}$, $d\in{\Bbb Z}$, $C_0=4\pi e^{\gamma}$ and
$$z_{d,t}=2|d|+ 1/2 +it,\eqno(3.2)$$
one has:
$$\log T(d,t)=\log T\left( |d|,|t|\right)
=2\Re \left({\Gamma^{(1)}\over\Gamma}\left( z_{d,t}\right) -\log\pi\right);\eqno(3.3)$$
$$T\left( d,t_1\right) > T(d,t)\quad\hbox{if $\,\left| t_1\right| > |t|$}\;;\eqno(3.4)$$
$$\log T(d,0)=-2\log C_0+4\sum_{k=1}^{|d|}{1\over 2k-1}\;;\eqno(3.5)$$
$$\log T\left( 0,1/2\right) < -\log\left( 2 C_0\right)\;;\eqno(3.6)$$
$$T(d,t)={1\over\pi^2}\left( \left| z_{d,t}\right|^2 -\Re \left( z_{d,t}\right)\right) +O(1)
={|2d+it|^2\over\pi^2}+O(1).\eqno(3.7)$$

\goodbreak
\noindent{\Smallcaps Proof.}\  
Given the definition (2.8), the identity (3.3) follows by logarithmic differentiation
of (2.6) and an appeal to the reflection principle. By (3.3), $T(d,t)$ is positive valued.
One has, moreover, 
$${{\rm d}^2\over{\rm d}z^2}\log\Gamma(z)=\sum_{n=0}^{\infty}{1\over (n+z)^2}\quad
\hbox{ for $|{\rm Arg}(z)|<\pi$}$$
(see, for example, [46, Section 12.16]), and so it follows  
from (3.3) and (3.2) that, for $d\in{\Bbb Z}$, one has 
$${\partial\over\partial t}\log T(d, t)
=4 t\sum_{n=0}^{\infty}{\Re \left( n+z_{d,t}\right)\over\left| n+z_{d,t}\right|^4}$$ 
at all points $t\in{\Bbb R}$. 
By this equation one has
${\partial\over\partial t}\log T(d, -t) ={\partial\over\partial t}\log T(d, t)  >0$ for $t>0$, $d\in{\Bbb Z}$.
Therefore, when the integer $d$ is given, $\log T(d, t)$ is a strictly increasing function of $|t|$.
The result (3.4) follows, given the strict monotonicity of the function 
$\exp(x)$ on ${\Bbb R}$.

Since ${\Gamma^{(1)}\over\Gamma}(1/2)
={\Gamma^{(1)}\over\Gamma}(1)-2\log 2=-\gamma -\log 4\,$
(see [46, Sections 12.1 and 12.15]), one may deduce (3.5) from (3.3) and (3.2) by logarithmic differentiation of
the identity $\Gamma(z)=(z-1)(z-2)\ldots (z-2|d|)\Gamma(z-2|d|)$. Since
$${{\rm d}\over{\rm d}z}\,\log\Gamma(z)
=-\gamma -{1\over z} +z\sum_{n=1}^{\infty}{1\over n(n+z)}\quad\hbox{for $|{\rm Arg}(z)|<\pi$}$$
(see [46, Section 12.6]), one has
$$\Re \left({\Gamma^{(1)}\over\Gamma}\left({1+i\over 2}\right)\right)
=-\gamma -1+2\sum_{n=1}^{\infty}{n+1\over n\left( 4n(n+1)+2\right)}
<-\gamma -1+2\left({1\over 5}+{\zeta(2)-1\over 4}\right)\;.$$ 
By this, combined with the 
equality $\zeta(2)=\pi^2 /6$, the inequality 
$\pi^2 /6<\gamma+11/5-\log(8/\pi)$ and the case  $d=0$, $t=1/2$ of 
(3.2) and (3.3), one obtains the inequality  (3.6).

By (3.3) and (3.2),  the case $\delta =\pi /2$ of (3.1) of Lemma~4 shows that,
for $d\in{\Bbb Z}$ and $t\in{\Bbb R}$,
$$\log T(d,t) 
=2\log\left| z_{d,t}\right| -{\Re \left( z_{d,t}\right)\over\left| z_{d,t}\right|^2}
-2\log\pi +O\left({1\over\left| z_{d,t}\right|^2}\right) .$$
Here we have, by (3.2), $\bigl| z_{d,t}\bigr|\geq\Re \bigl( z_{d,t}\bigr)\geq 1/2$, 
so the estimate (3.7) follows\quad$\square$

\medskip

\goodbreak\proclaim{\Smallcaps Corollary~6}. With $C_0$ as in the above lemma, it follows by (3.2), (3.4), (3.5) and (3.7) that
$$C_0^{-2}\leq T(d,t)\asymp \left| 2|d|+ 1/2 +it\right|^2\quad\hbox{for $d\in{\Bbb Z}$, 
$t\in{\Bbb R}$.}\eqno(3.8)$$

\bigskip

\goodbreak\proclaim{\Smallcaps Lemma~7}. Let $0<\eta\leq 1/4$. Then, for $d\in{\Bbb Z}$
and $s=\sigma +it$ with
$$-{5\over 2}\leq\sigma\leq{5\over 2}\quad\hbox{and}\quad t\in{\Bbb R}\;,\eqno(3.9)$$
one has
$$X_d(s)\ll_{\eta} T^{1/2-\sigma}\quad
\hbox{if either $\;T=T(d,t)\geq T(0 , 1/2)\;$ or $\;\sigma\not\in\bigcup\limits_{n=1}^2 (n-\eta , n+\eta)\;$,}\eqno(3.10)$$
where $X_d(s)$ and $T(d,t)$ are as given by (2.6) and (2.8).

\goodbreak
\noindent{\Smallcaps Proof.}\ 
By (2.6) and [36, Lemma 3] one has
$$\left| X_d(s)\right| =\pi^{2\sigma -1}\left|{\Gamma\left( 2|d|+1-s\right)\over\Gamma\left( 2|d|+s\right)}\right|
\leq\left|{2|d|+1+s\over\pi}\right|^{1-2\sigma}\quad\hbox{if $\,-1/2\leq\sigma\leq 1/2\,$.}$$
Therefore, if one puts $Z_{d,t}=\left| 2|d|+\textstyle{1/2}+it\right|$, then
$$\left| X_d\left( s_1\right)\right|
\ll Z_{d,t}^{1-2\Re \left( s_1\right)}\quad\hbox{for $s_1\in{\Bbb C}$ such that $\Im \left( s_1\right) =t$ 
and $\left|\Re \left( s_1\right)\right|\leq 1/2\,$.}\eqno(3.11)$$
For cases where $|\sigma|>1/2$ an estimate similar to (3.11) is needed.

If $n\in{\Bbb N}$, then unless $s_1$ is an integer not satisfying $-2|d|-n<s_1<2|d|+1$,
it follows
by (2.6) and the functional equation of $\Gamma$ that
$${X_d\left( s_1\right)\over X_d\left( s_1+n\right)}
=\pi^{-2n}\prod_{k=0}^{n-1}\left( 2|d|-s_1-k\right)\left( 2|d|+s_1+k\right)\;.\eqno(3.12)$$
Therefore (3.12) holds for $\left|\Re \left( s_1\right)\right| <3$ (and for all $n\in{\Bbb N}$) 
when the integer $d$ is non-zero. If however $d=0$, then one may observe that
when $T\left( 0,\Im \left( s_1\right)\right)\geq T( 0 , 1/2 )$ it must follow by (3.3) and (3.4)
of Lemma~5 that $\left|\Im \left( s_1\right)\right|\geq 1/2>0$. The case $d=0$ of (3.12) 
therefore holds for $n\in{\Bbb N}$ and $s_1\in{\Bbb C}$ such that 
$T\left( 0,\Im \left( s_1\right)\right)\geq T( 0 , 1/2)$
(for $\left|\Im \left( s_1\right)\right|\geq 1/2$ certainly implies
$s_1\not\in{\Bbb Z}$). 

Suppose now that $d\in{\Bbb Z}$ and $s=\sigma +it$ with $\sigma ,t\in{\Bbb R}$
such that (3.9) holds and $T(d,t)\geq T( 0 , 1/2)$.
If $-5/2\leq\sigma<-1/2$, then, by the discussion of the previous paragraph, 
(3.12) may be applied with $s_1=s$ and
$n=-[\sigma +1/2]\in\{ 1,2\}$. Consequently one obtains the bound
$$X_d(s)\ll Z_{d,t}^{2n}\left| X_d (s+n)\right|\;,$$
and it follows by (3.11) that
$$X_d(s)\ll Z_{d,t}^{2n+1-2(\sigma +n)}
=Z_{d,t}^{1-2\sigma}\quad\hbox{if $-5/2\leq\sigma <1/2\;$.}\eqno(3.13)$$
If instead $1/2<\sigma\leq 5/2$, then (3.12) may be applied with
$s_1=s-n$ and $n=-[1/2-\sigma]\in\{ 1,2\}$. One has either $d=0$ and $|t|\geq 1/2$, or
$2|d|\geq 2\geq 1/2+\max\{ |\sigma -2| , |\sigma -1|\}$, so that (3.12) yields the
lower bound
$$\left| {X_d (s-n)\over X_d (s)}\right|\geq\left( {1\over\pi}\max\left\{ {|d|\over 2} , {1\over 2} , |t|\right\}\right)^{2n}
\gg Z_{d,t}^{2n}\;.$$
One has also $|\Re (s-n)|=|\sigma -n|\leq 1/2$ and $\Im (s-n)=t$, so
it follows by the above and (3.11) that
$$X_d(s)\ll Z_{d,t}^{-2n}\left| {X_d (s-n)\over X_d (s)}\right|
\ll Z_{d,t}^{-2n+1-2(\sigma -n)}=Z_{d,t}^{1-2\sigma}\quad\hbox{if $1/2 <\sigma\leq 5/2\;$.}\eqno(3.14)$$
By the result (3.8) of Corollary~6, 
one has $Z_{d,t}\asymp (T(d,t))^{1/2}$ in  (3.11), (3.13) and (3.14),
which therefore yield all those cases of the bound (3.10) in which
$T(d,t)\geq T(0,1/2)$. 

The only cases of (3.10) requiring further consideration are
therefore those in which $T(d,t)<T(0,1/2)$ 
and $\min\{ |\sigma -1|<\eta , |\sigma -2|\}\geq\eta$.
In such cases it follows by (3.3)-(3.6) of Lemma~5 (and the inequality $\exp(4-\gamma) >2\pi$)
that $d=0$ and $|t|<1/2$. 
Therefore it follows by (2.6) that in all these cases one has $X_d(s)\ll 1$,
for the set $\{ s\in{\Bbb C} : \Re (s)\in[-5/2 , 1-\eta]\cup [1+\eta , 2-\eta ]\cup [2+\eta , 5/2]\ {\rm and}\
|\Im (s)|\leq 1/2\}$ is a bounded closed region containing no pole of $X_0(s)$. By (3.8) one has $1\ll T^{1/2 -\sigma}$ in (3.10), so
the proof is now complete\quad$\square$

\bigskip

\goodbreak\proclaim{\Smallcaps Lemma~8}. Let $d\in{\Bbb Z}$.
Suppose moreover that the functions $X_d(s)$ and $T(d,t)$ are as in (2.8) and (2.6), and that
$s=\sigma+it$, where 
$$-{1\over 2}\leq\sigma\leq{3\over 2}\eqno(3.15)$$
and $t\in{\Bbb R}$ satisfies
$$T(0 , 1/2)\leq T(d,t) =T\quad\hbox{(say).}\;\eqno(3.16)$$
Then one may define  a meromorphic complex function $\tau \mapsto G_d(s,\tau)$ by:
$$G_d (s,\tau) = {X_d (s-\tau)\over X_d (s)}\,T^{-\tau} -1
\quad\hbox{for $\tau\in{\Bbb C}$ with $s-2|d|-\tau\not\in{\Bbb N}$}\;.\eqno(3.17)$$
This function  is analytic on the open disc 
$\{ \tau\in{\Bbb C} : |\tau |< R_{d,s}\}\,$, where
$$R_{d,s}=\left| 2|d|+1-s\right|\asymp T^{1/2}\;,\eqno(3.18)$$
and one has $G_d (s , 0)=0$.
There is a unique complex sequence $\left( a_k (d,s)\right)_{k\in{\Bbb N}}$ 
with the property that
$$\sum_{k=1}^{\infty}a_k (d,s)\tau^k = G_d (s,\tau)
\quad\hbox{for all $\tau\in{\Bbb C}$ 
such that $|\tau|<R_{d,s}\;$.}\eqno(3.19)$$
For some absolute constant $B\geq 1$, the above sequence $\left( a_k (d,s)\right)_{k\in{\Bbb N}}$ 
satisfies:
$$a_k(d,s)\ll \left( B|t|/T\right)^{k/2}+\left( B/T \right)^{k/3}\ll B^k T^{-k/4}
\quad\hbox{for $k\in{\Bbb N}\;$,}\eqno(3.20)$$
$$a_1(d,s)= - {2(\sigma -1/2) ti\over \left| z_{d,t}\right|^2}+O\left( T^{-1}\right)\;,
\qquad a_2(d,s)= {ti\over \left| z_{d,t}\right|^2}+O\left( T^{-1}\right)\;,\eqno(3.21)$$
$$a_3(d,s)\ll T^{-1}\quad\hbox{and}
\quad a_4(d,s)= - {(1/2) t^2\over \left| z_{d,t}\right|^4} +O\left( T^{-9/8}\right)\ll T^{-1}\;,\eqno(3.22)$$
where (as in Lemma~5) $\,z_{d,t}=2|d|+ 1/2 +it$.

\goodbreak
\noindent{\Smallcaps Proof.}\ 
By (3.17) and (2.6),
$$G_d(s,\tau) =g_{d,s}(\tau) -1,\eqno(3.23)$$
where
$$g_{d,s}(\tau) =\left( X_d(s)\right)^{-1} X_d(s-\tau) T^{-\tau} 
=\left({\Gamma(2|d|+1-s)\over\Gamma(2|d|+s)}\right)^{-1}
\left({\Gamma(2|d|+1-s+\tau)\over\Gamma(2|d|+s-\tau)}\right)
\left(\pi^2 T\right)^{-\tau}\;.\eqno(3.24)$$
Assuming that $s$ is neither a pole nor zero of $X_d(s)$, the function
$g_{d,s}(\tau)$ is analytic for
$$|\tau|<\left| 2|d|+1-s\right| =R_{d,s},$$
and is non-zero for
$$|\tau|<\left| 2|d|+s\right| =R_{d,s}^*\quad\hbox{(say).}$$
Therefore, and since $g_{d,s}(0)=1$, it follows by (3.23) and the theory of Taylor and Laurent series 
(for which see [42, Sections 2.43 and 2.9]) 
that (3.19) holds if and only if $a_k(d,s)=g_{d,s}^{(k)}(0)/(k!)$ for all $k\in{\Bbb N}$.

Completion of the proof requires verification of (3.18) and (3.20)-(3.22). For the latter part
it suffices to estimate of $g_{d,s}(\tau)$ in a subset of the region indicated in (3.19).
Cauchy's inequality (stated in [42, Section~2.5]) is useful in deducing (3.20).

If $d\neq 0$, then by (3.15)
$$\min\left\{ R_{d,s} , R_{d,s}^*\right\}\geq
\left| 2|d|- 1/2 +it\right|\geq (3/5)\left| z_{d,t}\right| ,$$
where $z_{d,t}=2|d|+ 1/2 +it$. If $d=0$, then by (3.3) and (3.4) of
Lemma~5 the lower bound on $T(d,t)$ in (3.16) implies that $|t|>1/2$, and so
$$\min\left\{ R_{d,s} , R_{d,s}^*\right\}\geq |t|\geq 2^{-1/2}\left| 1/2 +it\right| 
= 2^{-1/2}\left| z_{d,t}\right|$$
in this case. Since $2^{-1/2}>3/5$, one may deduce that
$$\min\left\{ \left| (2|d|+1)-s\right| , \left| (-2|d|)-s\right|\right\} \geq 3/10>0\;,$$
so that, by virtue of (3.15), one is sure to have $X_d(s)\in{\Bbb C}-\{ 0\}$.
Moreover, given (3.15), the facts already gathered are 
sufficient to justify the conclusions that
$$3\left| z_{d,t}\right|\geq |z|+\left|\sigma - 1/2 \right|
\geq\min\left\{ \left| 2|d|+1-s\right| , \left|2|d|+s\right|\right\}\geq (3/5)\left| z_{d,t}\right|\eqno(3.25)$$
and
$$\max\left\{ \left|{\rm Arg}(2|d|+1-s)\right| , \left|{\rm Arg}(2|d|+s)\right|\right\}\leq (3/4)\pi\;,$$
in all cases covered by the lemma. Consequently,
for $\tau\in{\Bbb C}$ such that
$$|\tau|\leq  (3/5)\left| z_{d,t}\right|\sin(\pi /6)=(3/10)\left| z_{d,t}\right|\;,\eqno(3.26)$$
one has
$$\max\left\{ \left|{\rm Arg}(2|d|+1-s+\tau)\right| , 
\left|{\rm Arg}(2|d|+s-\tau)\right|\right\}\leq (3\pi/4) + (\pi/6) =(11/12)\pi\;.\eqno(3.27)$$
Note that, by  (3.25) and 
the result (3.8) of Corollary~6, one has $R_{d,s}\asymp T^{1/2}$, as claimed
in (3.18).

Since one may define a single valued and analytic branch of $\log\Gamma(z)$ on 
${\Bbb C}-(-\infty ,0]$, it follows by (3.24), (3.27) and the result (3.3) of Lemma~5 that, 
for $\tau\in{\Bbb C}$ satisfying
(3.26), one has:
$$\eqalignno{{{\rm d}\over{\rm d}\tau}\,\log g_{d,s}(\tau)
 &={\Gamma^{(1)}\over\Gamma}( 2|d|+1-s+\tau)+{\Gamma^{(1)}\over\Gamma}( 2|d|+s-\tau)
-(\log T +2\log\pi) =\cr
 &=\left({\Gamma^{(1)}\over\Gamma}( 2|d|+1-s+\tau)-{\Gamma^{(1)}\over\Gamma}\left( z_{d,-t}\right)\right)
+\left({\Gamma^{(1)}\over\Gamma}( 2|d|+s-\tau)-{\Gamma^{(1)}\over\Gamma}\left( z_{d,t}\right)\right).&(3.28)}$$
If $\left| z_{d,t}\right|\leq (10/3)^3$, then for $\tau$ as in (3.26) it follows by 
(3.27), (3.28) and the case $\delta =\pi/12$ 
of the result (3.1) of Lemma~4 that
$$\eqalign{
\left| {{\rm d}\over{\rm d}\tau}\,\log g_{d,s}(\tau)\right| &\leq
\left|\log\left( 2|d|+1-s+\tau\right) -\log\left( 2|d|+ 1/2  -it\right)\right| + \cr
 &\quad +\left|\log\left( 2|d|+s-\tau\right) -\log\left( 2|d|+ 1/2  +it\right)\right|
+\bigl| O\bigl(\,\left| z_{d,t}\right|^{-1}\,\bigr)\bigr|  \leq \cr
 &\leq 2\left(\log(10/3)+17\pi/12\right)+O(1)\ll 1.}$$
From this it follows, since $\log g_{d,s}(0)=\log 1 =0$, that if 
$\left| z_{d,t}\right|\leq (10/3)^3$ then, for $\tau$ satisfying (3.26), one has
$\log g_{d,s}(\tau)\ll |\tau|\leq (10/3)^2$, and so $g_{d,s}(\tau)\ll 1$.
Consequently (given (3.18), (3.19) and (3.23)) it follows by Cauchy's inequality
that, in cases where $\left| z_{d,t}\right|\leq (10/3)^3$, one has
$$a_k (d,s)\ll\left( (3/10)\,\left| z_{d,t}\right|\right)^{-k}\quad\hbox{for $k\in{\Bbb N}\;$.}$$
This confirms (3.20) in those cases, 
since (3.8) shows $\left| z_{d,t}\right|\gg T^{1/2}\gg T^{1/3}$. One may
also  verify that in those same cases, in which $T\asymp 1$, 
the estimates in (3.21) and (3.22) are
no stronger than those provided (for $k=1,2,3,4$) by the confirmed bound (3.20).

The above completes the proof of the lemma for cases with $\left| z_{d,t}\right|\leq (10/3)^3$, so
henceforth it is to be supposed that $\left| z_{d,t}\right| >(10/3)^3$. This supposition is more than
sufficient to ensure that, for $\tau$ satisfying $(3.26)$ 
(and with $\sigma$ as in (3.15)), one has
$$(1/2)\left| z_{d,t}\right|\geq\left( (3/10)+(3/10)^3\right)\left| z_{d,t}\right|
\geq |\tau| +1\geq\left|\tau_{\sigma}\right|\;,\eqno(3.29)$$
where
$$\tau_{\sigma} = \tau - \left(\sigma - 1/2 \right)\;.\eqno(3.30)$$
Assuming (3.26), it follows by (3.28), the estimate (3.1) of Lemma~4, and (3.29)-(3.30), that
$$\eqalign{
{{\rm d}\over{\rm d}\tau}\log g_{d,s}(\tau)
 &=\left(\log\left( z_{d,-t}+\tau_{\sigma}\right)-\log\left( z_{d,-t}\right)\right)
+\left(\log\left( z_{d,t}-\tau_{\sigma}\right)-\log\left( z_{d,t}\right)\right) + \cr
 &\quad +{1\over 2}\left( {1\over z_{d,-t}} - {1\over z_{d,-t}+\tau_{\sigma}}\right)
+{1\over 2}\left( {1\over z_{d,t}} - {1\over z_{d,t}-\tau_{\sigma}}\right)
+O\left( \left| z_{d,t}\right|^{-2}\right) = \cr 
 &={\tau_{\sigma}\over z_{d,-t}}+{\left( -\tau_{\sigma}\right)\over z_{d,t}}
+O\left({\left|\tau_{\sigma}\right|^2\over \left| z_{d,t}\right|^2}\right)
+O\left( {\left| \tau_{\sigma}\right|+1\over \left| z_{d,t}\right|^2}\right) = \cr
 &={2it\tau_{\sigma}\over \left| z_{d,t}\right|^2}
+O\left( {1+\left| \tau_{\sigma}\right|^2\over\left| z_{d,t}\right|^2}\right) 
={2it\tau\over \left| z_{d,t}\right|^2} 
- {2it\left(\sigma - 1/2 \right)\over \left| z_{d,t}\right|^2}
+O\left( {1+\left| \tau\right|^2\over\left| z_{d,t}\right|^2}\right) .
}$$
Since $\log g_{d,s}(0)=0$, it follows by the above that if $\tau$ satisfies
(3.26) then
$$\log g_{d,s}(\tau)
=\beta_1\tau +\beta_2\tau^2+O\left( {|\tau|+|\tau|^3\over \left| z_{d,t}\right|^2}\right) ,\eqno(3.31)$$
where one has 
$$\beta_1 =\beta_1 (d,s)
=- {2i\left(\sigma - 1/2 \right) t\over \left| z_{d,t}\right|^2}
\quad\hbox{and}\quad
\beta_2 =\beta_2 (d,s)
={it\over \left| z_{d,t}\right|^2}\;,\eqno(3.32)$$
which, by (3.15), implies:
$$\left|\beta_j\right|\leq 2 |t|/\left| z_{d,t}\right|^2\leq 2/\left| z_{d,t}\right|
\quad\hbox{for $j=1,2\,$.}\eqno(3.33)$$
Now, if $\tau$ satisfies
$$|\tau|\leq {\left| z_{d,t}\right|\over |t|^{1/2}+\left| z_{d,t}\right|^{1/3} }\;,\eqno(3.34)$$
then, since we are to suppose that $\left| z_{d,t}\right|^{1/3}>10/3$,
it is certainly the case that $\tau$ also satisfies (3.26). Therefore it may be deduced from 
(3.31) and (3.33) that, for all
$\tau\in{\Bbb C}$ satisfying the condition (3.34), one has
$$g_{d,s}(\tau)
=\exp\left(  O\left( \left| z_{d,t}\right|^{-1/3}\right) 
+O(1) +O\left(\left| z_{d,t}\right|^{-4/3} +1\right)\right) =\exp\left( O(1)\right)\ll 1.
$$
By this last bound, together with (3.18), (3.19), (3.23) and Cauchy's inequality, one finds that
$$a_k (d,s)\ll \left( |t|^{1/2}\left| z_{d,t}\right|^{-1}+\left| z_{d,t}\right|^{-2/3}\right)^k
\quad\hbox{for $k\in{\Bbb N}\,$.}\eqno(3.35)$$
By (3.35) and (3.8) one obtains the bounds stated in (3.20).

One may next observe that, since $\left| z_{d,t}\right|\geq (10/3)^3$,
the condition (3.26) is certainly satisfied
whenever $|\tau |\leq\left| z_{d,t}\right|^{1/3}$ (say), 
so that in such cases (3.31)-(3.33) apply. Moreover it follows by
(3.31)-(3.33) that, for $\tau$ satisfying $1\leq |\tau |\leq\left| z_{d,t}\right|^{1/3}$, one has:
$$\eqalignno{
g_{d,s}(\tau)
 &=\exp\left(\beta_1 \tau\right) \exp\left(\beta_2 \tau^2\right) 
\left( 1+O\left( |\tau|^3\left| z_{d,t}\right|^{-2}\right)\right) = \cr
 &=\left( 1+\beta_1 \tau\right)\left( 1+\beta_2 \tau^2 +\beta_2^2 \tau^4/2\right)
+O\left( \left( |\tau|^2+|\tau|^3\right) 
\left| z_{d,t}\right|^{-2}+|\tau|^6 \left| z_{d,t}\right|^{-3}\right) = \cr
 &= 1+\beta_1\tau +\beta_2\tau^2 +\beta_2^2\tau^4 /2
+O\left( |\tau|^3\left| z_{d,t}\right|^{-2}\right) .&(3.36)}$$ 
Recalling now that $\left| z_{d,t}\right| >(10/3)^3>2^4$, it is implied
by (3.19), (3.23), (3.25) and (3.35) that 
$$\eqalignno{
g_{d,s}(\tau)
 &= 1+\sum_{k=1}^4 a_k (d,s)\tau^k
+\sum_{k=5}^{\infty} 
O\left( \left( {2|\tau|\over\left| z_{d,t}\right|^{1/2}}\right)^k\right) = \cr
 &=1+\sum_{k=1}^4 a_k (d,s)\tau^k
+O\left( |\tau|^3 \left| z_{d,t}\right|^{-2}\right)
\qquad\hbox{if $\ |\tau|\leq\left| z_{d,t}\right|^{1/4}\,$.}&(3.37)}$$  
Given (3.32), and given that $\left| z_{d,t}\right|^2\asymp T$ (by (3.8)),
a comparison of (3.37) with (3.36) at the four points $\tau =\pm 1 , \pm i$ (for example)
is sufficient to establish the first three estimates of (3.21) and (3.22).
The final estimate of (3.22) then follows by comparing (3.37) with (3.36) at
the point $\tau =\left| z_{d,t}\right|^{1/4}$\quad$\square$

\medskip

\goodbreak\proclaim{\Smallcaps Corollary~9}. Subject to the hypotheses of the
above lemma, and with the same 
absolute constant $B\geq 1$ as in (3.20),
one has, in (3.18) and (3.19), $R_{d,s}=| 2|d|+1-s|\geq B^{-1}T^{1/4}$ and, for all integers $K\geq 0$,
$$\sum_{k=K+1}^{\infty} a_k (d , s) \tau^k
\ll_K \left( \left( {|t|/T}\right)^{(K+1)/2} + T^{-(K+1)/3}\right) |\tau|^{K+1}
\quad\hbox{if $\,\tau\in{\Bbb C}\,$ and $\;|\tau|\leq (2B)^{-1} T^{1/4}$.}\eqno(3.38)$$

\bigskip

\goodbreak\proclaim{\Smallcaps Lemma 10 (convexity and sub-convexity estimates)}. 
Let $0<\varepsilon\leq 1/4$, $d\in{\Bbb Z}$ and $s=\sigma +it$, where
$t\in{\Bbb R}$ and $\sigma\geq -5/2$. 
Put $T=T(d,t)$ (where $T(d,t)$ is given by (2.8) and (2.6)).
Then
$$\zeta\left( s,\lambda^d\right) -{(\pi /4)\delta_{d,0}\over s-1}
=O\left( \varepsilon^{-2} T^{\max\{ 0 , (1-\sigma)/2+\varepsilon , 1/2-\sigma\} }\right)\;,\eqno(3.39)$$
where
$$\delta_{a,b}=\cases{1 &if $a=b\;$, \cr 0 &otherwise.}\eqno(3.40)$$
In the case $\sigma =1/2$ this may be strengthened to:
$$\zeta\left( 1/2 +it , \lambda^d\right)\ll_{\varepsilon} T^{1/6 +\varepsilon}\;.\eqno(3.41)$$

\goodbreak
\noindent{\Smallcaps Proof.}\ 
In cases where $\sigma >3/2$, the bound (3.39) follows since $\left| 1/(s-1)\right| <2$,
and since, by (2.9) and the case $\varepsilon =1/3$ of (2.12),
$$\left|\zeta\left( s , \lambda^d\right)\right|
\leq\sum_{n=1}^{\infty}
{\left|\delta\left(\Lambda^d ,n\right)\right|\over n^{\sigma}}
\ll\zeta\left(\sigma -1/3\right) <\zeta(7/6)$$
(this suffices, given that $\varepsilon^{-2}>1$ and
$\max\{ 0 , (1-\sigma)/2+\varepsilon , 1/2-\sigma\} =0$ in (3.39) ).

If instead $-5/2\leq\sigma\leq 3/2$ and $T(d,t)<T( 0 , 1/2 )$, then
by (3.4)-(3.6) of Lemma~5 (and the inequality $2\pi <\exp( 4-\gamma)$)
one must have $d=0$ and $|t|<1/2$, which implies that
$|s-1|=|\sigma -1+it |<|\sigma -1|+1/2\leq 4$.
Moreover, it follows by the case $\rho =4$ of (2.7) that (3.39) holds
whenever $|s-4|\leq 4$ (one need only note that, for such $s$, 
$\max\{ 0 , (1-\sigma)/2+\varepsilon , 1/2-\sigma\}\in [0,7/2]$,
and that $T^{7/2}\gg 1$, by (3.8)).

Given the above, the only cases left to consider are those in which
one has both
$-5/2\leq\sigma\leq 3/2$ and $|s-1|>4$, so that $T(d,t)\geq T(0 , 1/2)$
(as follows on account of the inference established in the first
sentence of the previous paragraph) . If, in such a case, one has also $\sigma <-1/2$,
so that $-5/2\leq\sigma <-1/2$, $|s-1|>4$ and $T(d,t)\geq T(0 , 1/2)$, then by Lemma~7 and the functional equation (2.5)
one finds that
$$\zeta\left( s , \lambda^d\right)
\ll T^{1/2-\sigma}\left|\zeta\left( 1-s , \lambda^{-d}\right)\right|\;.$$
where $\zeta\left( 1-s , \lambda^{-d}\right)\ll 1$, since $\Re (1-s)=1-\sigma >3/2$.
The bound (3.39) then follows, 
given that one has $\left| 1/(s-1)\right| <1/4$ and $5/4\geq -\sigma /2>1/4\geq\varepsilon$ and, by (3.8), $T\gg 1$.

In those cases of the lemma still to be considered
one may make use of Rademacher's bounds in [36, Theorem~4 and Theorem~5],
which imply that if $-1/2\leq\sigma\leq 3/2$ and $|s-1|>4$, then
$$\zeta\left( s , \lambda^d\right)
\ll\zeta^2\!\left( 1+\eta\right) \left| 2|d|+1+s\right|^{\eta +1-\sigma}\;,$$
for all $\eta\in(0,1/2]$ such that $-\eta\leq\sigma\leq 1+\eta$.
Now  $\zeta\left( 1+\eta\right)\ll\eta^{-1}$ 
for $\eta\in(0,1/2]$, 
while for $-1/2\leq\sigma\leq 3/2$ one has
$\left| 2|d|+1+s\right|\asymp\left| 2|d|+ 1/2 +it\right|$, and so, by taking (in the above)
$$\eta =\max\left\{ 2\varepsilon , \left|\sigma -1/2\right| -1/2\right\}
=\max\left\{ 2\varepsilon +1-\sigma , \left| 1/2 -\sigma\right| +1/2 -\sigma\right\}
- \left( 1- \sigma \right)\;,$$
one obtains:
$$\zeta\left( s , \lambda^d\right)
\ll \varepsilon^{-2}
\left| 2 |d|+ 1/2 +it\right|^{\max\left\{ 2\varepsilon +1-\sigma , 0 , 1-2\sigma\right\} }\;.
\eqno(3.42)$$
By the result (3.8) of Corollary~6, it follows 
that (3.42) implies (3.39) whenever 
$-1/2\leq\sigma\leq 3/2$ and $|s-1|>4$. 
This completes the proof of (3.39) in all the cases within the scope of the lemma.

The bound (3.41) follows from Kaufman's theorem in [26], which in fact
gives the sharper bound $\zeta\left( 1/2 +it , \lambda^d\right)
\ll T^{1/6}\log^C (T+2)\,$ (where $C$ is some positive absolute constant). 
See [27] and [41] for similar `sub-convexity' bounds on Hecke zeta functions
for general number fields.
\quad$\square$

\bigskip

\goodbreak\proclaim{\Smallcaps Lemma 11 (an approximate functional equation)}. 
Let $C_0=4\pi e^{\gamma}$. Let $\varepsilon\in(0,1/4)$ and $b\in(1,\infty)$, 
and suppose moreover that the 
infinitely differentiable function 
$\rho : (0,\infty)\rightarrow{\Bbb R}$ satisfies
$$\rho(u)+\rho\left({1/u}\right) =1\quad\hbox{for $u>0$}\eqno(3.43)$$
and
$$\rho(u)=0\quad\hbox{for $u\geq b$.}\eqno(3.44)$$
Let $d\in{\Bbb Z}$ and $\sigma+it=s\not\in\{ -2|d| , 2|d|+1\}$, where $\sigma, t\in{\Bbb R}$ and
$$-1/3\leq\sigma\leq 4/3,\eqno(3.45)$$
and take $X_d(s)\in{\Bbb C}$, the real 
sequence $\bigl( \delta\bigl( \Lambda^d , n\bigr)\bigr)_{n\in{\Bbb N}}$ 
and $T(d,t)\in(0,\infty)$ to be as in (2.6), (2.10) and (2.8), respectively. 
Then, for all non-negative integers $K$, and all $x,y,T\in(0,\infty)$ satisfying both
$$xy=T=T(d,t)\eqno(3.46)$$
and
$$1/\left( 2 C_0\right)\leq by\leq 2 C_0 T,\eqno(3.47)$$
one has
$$\eqalign{\zeta\left( s,\lambda^d\right)
 &=\sum_{n=1}^{\infty}\rho\left({n\over x}\right) \delta\left(\Lambda^d,n\right) n^{-s}
+X_d(s)\sum_{n=1}^{\infty}\widetilde\rho_K\left({n\over y}\right) 
\delta\left(\Lambda^{-d},n\right) n^{s-1} \ + \cr
 &\quad +{(\pi /4)\delta_{d,0}\over s-1}\,e^{-|s-1|}
+E_K\left( s , \lambda^d\right)\;,}\eqno(3.48)$$
where $\delta_{a,b}$ is as defined in (3.40), while 
$$\widetilde\rho_K (u)
=\rho(u)+\sum_{k=1}^K (-1)^k a_{k}(d,s)\left( u\,{{\rm d}\over{\rm d}u}\right)^{k}\rho(u)\eqno(3.49)$$
(with the same coefficients $a_{k}(d,s)$ as in (3.17)-(3.19) of Lemma~8) and
$$\eqalignno{
E_{K}\left( s, \lambda^d\right) 
&\ll_{\varepsilon ,K}  \,x^{1/2-\sigma}\left( (|t|/T)^{\alpha_K}+T^{-\beta_K}\right) 
\int\limits_{-T^{\varepsilon}}^{T^{\varepsilon}}
\left|\zeta\left( 1/2 +i(t+v) , \lambda^{d}\right)\right|
\,{{\rm d}v\over 1+v^2} \ll_{\varepsilon , K} &(3.50)\cr
&\ll_{\varepsilon , K} \,T^{1/6+\varepsilon} x^{1/2 -\sigma}\left( (|t|/T)^{\alpha_K}+T^{-\beta_K}\right) 
&(3.51)}$$
with exponents
$$\left(\alpha_K , \beta_K\right) 
=\cases{( 1 , 1 ) &if $K=0,1\;$,\cr ( 2 , 1 ) &if $K=2\;$,\cr $( (K+1)/2 , (K+1)/3 )$ &if $K\geq 3\;$,}
\eqno(3.52)$$
and implicit constants that depend only upon $\varepsilon$, $K$, the function $\rho(u)$ and $b$.
\medskip

\goodbreak
\noindent{\Smallcaps Proof.}\ 
Note firstly that in this proof we generally omit to indicate where there is
some dependence of an implicit constant upon either the function 
$\rho(u)$ or the related number $b$, so a statement 
such as `$U\ll_{\varepsilon} 1$' might mean only that $U\ll_{\varepsilon , \rho} 1$, or that $U\ll_{\varepsilon , b} 1$.

It helps to distinguish two cases: a case in which $T$, in (3.46), is 
`sufficiently large', and the complementary case, in which $T$ is `of bounded magnitude'. 
Taking the latter case first, 
one supposes that
$$T\leq (6B)^{4/(1-4\varepsilon)}\;,\eqno(3.53)$$
where $B\geq 1$ is the absolute constant of Lemma~8 and Corollary~9.
Then, by (3.8) and (3.53),
$$1\ll T\ll_{\varepsilon} 1,\quad d\ll_{\varepsilon} 1\quad\hbox{and}\quad t\ll_{\varepsilon} 1,\eqno(3.54)$$
$$1\ll x , y\ll_{\varepsilon} 1\eqno(3.55)$$
and
$$\int\limits_{-T^{\varepsilon}}^{T^{\varepsilon}}
\left|\zeta\left( 1/2 +i(t+v) , \lambda^{d}\right)\right|
\,{{\rm d}v\over 1+v^2}
\gg_{\varepsilon}
\,\int\limits_{t-1/C_0}^{t+1/C_0}
\left|\zeta\left( 1/2 +iw , \lambda^{d}\right)\right|
{{\rm d}w}\gg_{\varepsilon} 1,\eqno(3.56)$$
where the last lower bound follows since (3.54) restricts $d$ to a finite subset
of ${\Bbb Z}$ determined by $\varepsilon$, and restricts $t$ to a 
closed bounded interval $[-h,h]$ also determined by $\varepsilon$, while
the relevant integrand, 
$\left|\zeta\left( 1/2 +iw , \lambda^{d}\right)\right|$, 
is a continuous non-negative valued real function with 
only isolated zeros (each function $\zeta\bigl( z , \lambda^d\bigr)$ 
being analytic for $z\neq 1$, and, by the Euler product (2.14), non-zero for $\Re (z)>1$), 
so that the last integral in (3.56) is a continuous positive valued function of $t$ with an
infimum over $t\in[-h,h]$ that is also positive (this infimum necessarily being an
attained minimum of the function).

Therefore, for proof of the case of the lemma that we are currently considering, 
it suffices to check that:
$$\zeta\left( s , \lambda^d\right) -{(\pi /4)\delta_{d,0}\over s-1}\,e^{-|s-1|}
\ll_{\varepsilon} 1\;,\eqno(3.57)$$
$$\sum_{n=1}^{\infty}\rho\left( {n\over x}\right) \delta\left(\Lambda^d , n\right) n^{-s}
\ll_{\varepsilon} 1\eqno(3.58)$$
and
$$X_d(s)\sum_{n=1}^{\infty}\widetilde\rho_K\left( {n\over y}\right) \delta\left( \Lambda^{-d} , n\right)
n^{s-1}\ll_{\varepsilon , K} 1\;,\eqno(3.59)$$
where (both here and subsequently) the implicit constants may depend
on the choice of function $\rho(u)$ and associated constant $b>1$.

Since $\bigl( 1-\exp( -|s-1| )\bigr) / |s-1| < 1$, 
the bound (3.57) is a consequence  of (3.45), (3.54) and the bound (3.39) of Lemma~10.
The bound  (3.58) is straightforward to verify, given (3.44), (3.45),  (3.55) 
and the upper bound in (2.12). If $T\geq T(0 , 1/2)$, then it is similarly straightforward
to verify (3.59), since then (3.45), (3.54) and the case $\eta =1/4$ of the bound (3.10) in Lemma~7
show that $X_{d}(s)\ll T^{5/6}\ll_{\varepsilon} 1$, while (3.20) of Lemma~8 implies that
$a_{k}(d , s)\ll_k T^{-k/4}\ll_K 1$, for $k=1,2,\ \ldots\ ,K$.

If $T<T(0 , 1/2)$, then by (3.3)-(3.6) of Lemma~5 one has $d=0$ and $|t|<1/2$.
Though Lemma~8 itself does not cover this, one may in such a case nevertheless persist in
defining the sequence $\bigl( a_{k}(d,s)\bigr)_{k\in{\Bbb N}}$ through (3.17)-(3.19), 
given that (3.45) and the hypothesis that $s\not\in\{ -2|d| , 2|d|+1\} =\{ 0 , 1\}$
combine to imply that $G_{0}(s , \tau)$, in (3.17), is well defined and analytic
on the non-empty open set $\{ \tau\in{\Bbb C} : |\tau|< |1-s| \}$.
In particular, if $d=0$ then $X_{d}(s)=X_{0}(s)=\pi^{2s-1}\Gamma(1-s)/\Gamma(s)\in{\Bbb C}-\{ 0\}$
and, by (3.17), $G_{d}(s,0)=T^{0}-1=0$. Therefore, when $T<T(0 , 1/2)$ in (3.46), one may
verify (3.59) just by noting, firstly, that (3.47) and the inequality (3.6) of Lemma~5 then imply 
the inequality $by<1$, and, secondly, that by this inequality, (3.44) and (3.49), 
it follows that each summand in (3.59) equals zero.

The above completes the verification of (3.57)-(3.59), and so the case of the lemma 
in which (3.53) holds has been proved. The rest of this proof need only deal
with the case complementary to (3.53), so it is henceforth to be supposed that
$$T>(6B)^{4/(1-4\varepsilon)}\;.\eqno(3.60)$$

As in [21, Equation (21)], put
$$R(z)=\int\limits_{0}^{\infty}\rho(u) u^{z-1} {\rm d}u\quad\hbox{for $\Re (z)>0\,$.}\eqno(3.61)$$
Then (again as in [21]) it follows by (3.43), (3.44) and integration by parts that, where defined,
$$R(z)= - {1\over z}\,\int\limits_{1/b}^{b} \rho^{(1)}(u) u^z {\rm d}u\;,\eqno(3.62)$$
and that by means of this representation, (3.62), the function $R(z)$ has an analytic
continuation to ${\Bbb C}-\{ 0\}$ with, at $z=0$, a simple pole having
residue equal to $\rho(1/b) -\rho(b) = 1 - 0 = 1$.
Given (3.62), the condition (3.43) has the pleasing effect of ensuring that 
$$R(-z)=-R(z)\quad\hbox{for all $z\in{\Bbb C} - \{ 0\}\,$.}\eqno(3.63)$$
By (3.62) and the application of further integrations by parts, one obtains bounds
for $R(z)$ which, in combination with (3.63), imply the bounds: 
$$R(z)\ll_j {b^{|\Re (z)|}\over |z| (|z|+1)^{j-1}}
\quad\hbox{for $j\in{\Bbb N}$ and $z\in{\Bbb C}-\{ 0\}\,$,}\eqno(3.64)$$
with the relevant implicit constants being 
$b^{j-1}\int_{1/b}^{b}\bigl|\rho^{(j)}(u)\bigr| {\rm d}u\,$ ($j=1,2,\ldots\ $). 

Since $\rho(u)$ is infinitely differentiable, the hypotheses (3.43) and (3.44)
ensure that, 
given (3.61) and the case $j=2$ of (3.64), 
it follows from the Mellin inversion formula of [22, Equation (A.2)]  that one has 
$${1\over 2\pi i}\,\int\limits_{c-i\infty}^{c+i\infty}
R(z) u^{-z} {\rm d}z
=\rho(u)
\quad\hbox{for $c,u>0\,$.}\eqno(3.65)$$
Differentiating under the integral sign here one obtains:
$$- {1\over 2\pi i}\,\int\limits_{c-i\infty}^{c+i\infty} z R(z) u^{-z-1} {\rm d}z
=\rho^{(1)}(u)$$
(a result justified by virtue of the fact that the case $j=3$ of (3.64) 
implies that the improper integral on the left hand side of the equation 
is uniformly convergent for all $u$ in any given bounded closed
interval $\bigl[ u_0 , u_1\bigr]\subset(0,\infty)$). Upon multiplying both sides of the
above equation by $u$, one obtains the case $k=1$ of the equation
$${(-1)^k\over 2\pi i}\,\int\limits_{c-i\infty}^{c+i\infty}
z^k R(z) u^{-z} {\rm d}z
=\left( u\,{{\rm d}\over{\rm d}u}\right)^k \rho(u)
=\rho_k(u)\quad\hbox{(say),\ \  
for $c,u>0\,$.}\eqno(3.66)$$
The general case (all $k\in{\Bbb N}$) follows by induction (appealing, at the $k$-th step,
to the case $j=k+3$ of (3.64)).

Applying either (3.65) or (3.66), for $u=n/x$, and then multiplying the result
by $\delta\bigl(\Lambda^d ,n\bigr) n^{-s}$, before summing over $n\in{\Bbb N}$,
one finds by (3.44), (2.9) and (2.12) 
(which guarantee the required uniform absolute convergence of
$\sum_{n=1}^{\infty}\delta\bigl(\Lambda^d , n\bigr) n^{-s-z}$) that if $x>0$ 
and $c>\max\{ 0 , 1-\sigma\}$, then
$${(-1)^k\over 2\pi i}\,\int\limits_{c-i\infty}^{c+i\infty}
z^k R(z) x^{z} \zeta\left( s+z , \lambda^d\right) {\rm d}z
=\sum_{n=1}^{\infty}\rho_k\left( {n\over x}\right)
{\delta\left(\Lambda^d , n\right)\over n^{s}}\quad\hbox{for $k\in{\Bbb N}\cup\{ 0\}\,$,}\eqno(3.67)$$
where $\rho_k(u)=\rho(u)$ if $k=0$, and otherwise is as in (3.66).

A deduction may now be made from (3.67) which, though it is a digression
from the central theme of this proof, will later serve to keep
the error term estimate (3.50) as simple as it is. For this purpose it is
helpful to temporarily restrict  the discussion to the special case of
a function $\rho(u)$ such that (3.44) holds with $b=\sqrt{2}$
(there need be no question concerning the existence of such a function
$\rho(u)$, since the hypotheses of the lemma imply that if the given
$\rho(u)$ is not itself of the type sought, then for some choice of
constant $q>0$ the functions $\rho\bigl( u^{q}\bigr)$ will be).
In this special case, the application of (3.67) with $k=0$, $s=1/4+it$ and
$x=b=\sqrt{2}$ shows that
$${1\over 2\pi i}\,\int\limits_{1-i\infty}^{1+i\infty}
R(z) b^{z} \zeta\left( 1/4+it+z , \lambda^d\right) {\rm d}z= 1\;.$$
Indeed, by (3.43), (3.44) and (2.10),  one has
$\rho(n/x)=0$ for $n\geq 2$, while $\rho(1/x)=\rho(1/b)=1-\rho(b)=1$
and $\delta\bigl(\Lambda^d , 1\bigr)=1$.
By the estimate (3.39) of Lemma~10, together with (3.8) and the bounds of (3.64),
it is permissible  to shift the above line of integration to $\Re (z)=1/4$,
though if $d=0$ then the pole of $\zeta\bigl( 1/4+it+z , \lambda^d\bigr)$ at
$z=3/4-it$ makes it necessary to adjust the shifted integral by
addition of the appropriate residue. Then, by substituting $z=1/4+iv$, one obtains:
$${1\over 2\pi}\,\int\limits_{-\infty}^{\infty}
R\left(\textstyle{1\over 4}+iv\right) b^{1/4+iv} 
\zeta\left( 1/2 +i(t+v) , \lambda^d\right) {\rm d}v
+(\pi /4)\delta_{d,0} b^{3/4-it} R\left(\textstyle{3\over 4} -it\right) 
= 1\;,\eqno(3.68)$$
where $\delta_{m,m}=1$ and $\delta_{m,n}=0$ if $m\neq n$.
Using (3.8), the bound (3.39) of Lemma~10 and the case  $j=[3/(2\varepsilon)]+3$ (say)
of (3.64), one deduces from equation (3.68) that
$${1\over 2\pi}\,\int\limits_{-T^{\varepsilon}}^{T^{\varepsilon}}
R\left(\textstyle{1\over 4}+iv\right) b^{1/4+iv} 
\zeta\left( 1/2 +i(t+v) , \lambda^d\right) {\rm d}v
+O_{\varepsilon}\left( T^{-1}\right) = 1 $$
(any dependence of the implicit constant here upon $\rho(u)$
may be effectively nullified through a `one time only'
selection of $\rho(u)$ from amongst the set of all functions having the
required properties). Since it is evidently not possible
that both of the two terms on the left side of this last equation have
absolute value less than $1/2$, one therefore must have either
$$T\ll_{\varepsilon} 1,\eqno(3.69)$$
or else 
$$\int\limits_{-T^{\varepsilon}}^{T^{\varepsilon}}
\left| R\left(\textstyle{1\over 4}+iv\right)  
\zeta\left( 1/2 +i(t+v) , \lambda^d\right)\right| {\rm d}v
\geq\pi b^{-1/4}\;.\eqno(3.70)$$
As a corollary of the dichotomy just discerned one obtains a useful lower bound,
$$\int\limits_{-T^{\varepsilon}}^{T^{\varepsilon}}
\left|\zeta\left( 1/2 +i(t+v) , \lambda^d\right)\right|\,{{\rm d}v\over 1+v^2}
\gg_{\varepsilon} 1\;,\eqno(3.71)$$
for if $T$ satisfies (3.69) then (3.71) follows by reasoning similar to that
used to obtain (3.56) (subject to (3.53)), whereas if it is instead 
(3.70) which holds, then (3.71) follows from (3.70) by use of
the case $j=2$ of (3.64). It is certain that the lower bound (3.71) could be
greatly improved, and refined upon, by very generally applicable 
methods of Balasubramanian and Ramachandra
(for which see [37,2,38]).

With the bound (3.71) established, the digression on the special case $b=\sqrt{2}$ is
now ended. Returning to consideration of the 
general case (though with (3.60) still in force) suppose now that $x$ and $y$ are as in (3.46) and (3.47). 
By applying (3.67) for $k=0$ and $c=3/2$ (so that $c+\sigma\geq 7/6>1$)
and then shifting the contour of integration to the line $\Re (z)=-c=-3/2$ (say),
one arrives at the equation
$${1\over 2\pi i}\,\int\limits_{-c-i\infty}^{-c+i\infty}
R(z) x^{z} \zeta\left( s+z , \lambda^d\right) {\rm d}z
+\zeta\left( s , \lambda^d\right)
+\delta_{d,0} R(1-s) x^{1-s} \pi /4 
=\sum_{n=1}^{\infty}\rho\left( {n\over x}\right)
{\delta\left(\Lambda^d , n\right)\over n^{s}}\;,$$
in which appear the residue of the integrand at $z=0$ (the pole of $R(z)$) and,
if $d=0$, the residue at $z=1-s\neq 0$ (the pole of $\zeta\bigl( s+z , \lambda^0\bigr)$\/): 
the shifting of the contour may be justified
using the bounds (3.64), along with the result (3.8) of Corollary~6 
and the estimate for 
$\zeta\bigl( s , \lambda^d\bigr)$ implicit in the result (3.39) 
of Lemma~10, which
also contains the information needed for the residue calculation
at $z=1-s$. Using the functional equation (2.5), and the
substitution $z=-\tau$, one may rewrite the integral in the 
equation that was just arrived at, and so obtain a reformulation of that equation; 
by this reformulation, followed by two applications of (3.63), 
one finds that 
$$\zeta\left( s , \lambda^d\right)
=\sum_{n=1}^{\infty}\rho\left( {n\over x}\right)
{\delta\left(\Lambda^d , n\right)\over n^{s}}
+(\pi /4) \delta_{d,0}\,{R(s-1)  \over x^{s-1}} +I( d , s)\;,\eqno(3.72)$$
where
$$I( d , s)
={1\over 2\pi i}\,\int\limits_{c-i\infty}^{c+i\infty}
X_{d}(s-\tau) 
\zeta\left( 1-s+\tau , \lambda^{-d}\right){R(\tau)\over x^{\tau}}\,{\rm d}\tau\;,\eqno(3.73)$$
with (as before) $c=3/2$.

In the integral of (3.73) one has 
$\Re (1-s+\tau)=1-\sigma+3/2\geq 7/6>1$, by (3.45), and therefore $\zeta\bigl(1-s+\tau , \lambda^d\bigr)\ll 1$
(as follows, for example, by the case $\varepsilon =1/12$ of the result (3.39) of Lemma~10).  
By (3.64) and (3.45), one has also
$$R(\tau)\ll_j b^{3/2} |\tau|^{-j}\ll |\tau|^{-j}\quad\hbox{if $\sigma -1/2\leq\Re (\tau)\leq 3/2\,$, 
$\tau\neq 0$ and $j\in{\Bbb N}\,$,}$$
while, by (3.45) and (3.8), the  bound (3.10) of Lemma~7 shows that 
$$X_{d}(s-\tau)\ll\left| 2|d|+ 1/2 +i\left( t-\Im (\tau)\right)\right|^{2(2-\sigma)}
\ll \left( T |\tau|^{2}\right)^{2-\sigma}\ll T^3 |\tau|^5
\quad\,\hbox{when $\,\Re(\tau)=3/2\,$.}$$
Using the last two bounds to estimate parts of the integral in (3.73), one finds that
$$I( d , s)
={1\over 2\pi i}\,\int\limits_{c-i T^{\varepsilon}}^{c+i T^{\varepsilon}}
X_{d}(s-\tau) 
\zeta\left( 1-s+\tau , \lambda^{-d}\right){R(\tau)\over x^{\tau}}\,{\rm d}\tau
+O_j\left( T^{3-(j-6)\varepsilon} x^{-3/2}\right)\;,\eqno(3.74)$$
for each integer $j$ such that $j\geq 7$. 
By (3.45), $1/2-\sigma >-3/2$, while by (3.46) and (3.47) one has $T\gg x\gg 1$ here, 
so that on restricting to the case
$j=[(A+3)/\varepsilon]+7$ one obtains, in (3.74),
$$O_j\left( T^{3-(j-6)\varepsilon} x^{-3/2}\right)
\ll_{A,\varepsilon} x^{1/2 -\sigma} T^{-A}\quad\, 
\hbox{for $\,0\leq A<\infty$.}\eqno(3.75)$$

By (3.45), (3.46), (3.60) and the inequality (3.6) of Lemma~5, the
conditions (3.15) and (3.16) of Lemma~8 are satisfied, 
so by Lemma~8 and Corollary~9 one has:
$$X_{d}(s-\tau)
=X_{d}(s) T^{\tau}\left( 1+G_{d}(s , \tau)\right)\;,\eqno(3.76)$$
where $G_{d}(s , \tau)$ is analytic (as a function of $\tau$) on the
open disc $\{ \tau\in{\Bbb C} : |\tau|<R_{d,s}\}$ with radius
$R_{d,s}=| 2|d|+1-s|\geq T^{1/4}/B$.
Now, if $\tau$ satisfies
$$\sigma - 1/2 \leq \Re (\tau)\leq 3/2\quad\hbox{and}\quad
\left| \Im (\tau)\right|\leq T^{\varepsilon}\;,\eqno(3.77)$$
then, by (3.45) and (3.60), one has also
$$|\tau|<3/2+T^{\varepsilon}<T^{1/4}/(4B)+T^{1/4}/(6B)<T^{1/4}/(2B)\;,\eqno(3.78)$$
so that the bound (3.38) of Corollary~9 applies.
For later reference, observe moreover that the last two inequalities of (3.78)
certainly demonstrate that if $d=0$ then one has 
$$|1-\sigma|+|t|\geq |1-s|=R_{0,s}\geq T^{1/4}/B>3+2T^{\varepsilon}
>|1-\sigma|+2T^{\varepsilon}$$ 
(given (3.45)). Consequently 
$$|t|>2T^{\varepsilon}\quad\hbox{if $d=0\,$.}\eqno(3.79)$$

Suppose now that $K$ is an arbitrary non-negative integer. 
Since the conditions in (3.77) imply the inequalities in (3.78), 
and since $R_{d,s}\geq T^{1/4}/B$,  it follows by (3.76), 
and the identity (3.19) of Lemma~8, that one may rewrite  
the integrand in (3.74) by using: 
$$X_{d}(s-\tau)
=X_{d}(s) \left( T^{\tau} P_{K,d,s}^{+}(\tau) +T^{\tau} P_{K,d,s}^{-}(\tau) \right)\;,\eqno(3.80)$$
where
$$P_{K,d,s}^{+}(\tau)
=1+\sum_{k=1}^{K} a_{k}(d,s) \tau^k
=\sum_{k=0}^{K} a_{k}(d,s) \tau^k\quad\hbox{(with $a_{0}(d,s)=1$)}\eqno(3.81)$$
and
$$P_{K,d,s}^{-}(\tau)
=1+G_{d}(s,\tau) -P_{K,d,s}^{+}(\tau)
=\sum_{k=K+1}^{\infty} a_{k}(d,s) \tau^k\;.\eqno(3.82)$$
Using (3.80) and (3.75) in (3.74), one finds that, for arbitrary $A\in[0,\infty)$, one has 
$$I(d,s)
=X_{d}(s)\left( J_{\varepsilon ,K}^{+}(d,s) + J_{\varepsilon ,K}^{-}(d,s)\right)
+O_{A,\varepsilon}\left( x^{1/2 -\sigma} T^{-A}\right)\;,\eqno(3.83)$$
where (given (3.46) and (2.3)) 
$$J_{\varepsilon ,K}^{\pm}(d,s)
={1\over 2\pi i}\,\int\limits_{c-i T^{\varepsilon}}^{c+i T^{\varepsilon}}
P_{K,d,s}^{\pm}(\tau)
R(\tau) y^{\tau} \zeta\left( 1-s+\tau , \lambda^{\mp d}\right)\,{\rm d}\tau\;.\eqno(3.84)$$

By (3.81), the bounds (3.20) of Lemma~8, and the lower bound on $T(d,t)$ in (3.60),
one has
$$P_{K,d,s}^{+}(\tau)\ll 1+\sum_{k=1}^K B^k T^{-k/4} |\tau|^k
\ll_K |\tau|^K\quad\hbox{for $|\tau|\geq T^{\varepsilon}$.}$$
By this last bound, together with (3.64) and the fact (noted while obtaining (3.74)) that
$\zeta\bigl( z , \lambda^n\bigr)\ll 1$ when $\Re (z)\geq 7/6$, it follows 
from (3.84) that 
$$J_{\varepsilon ,K}^{+}(d,s)
={1\over 2\pi i}\,\int\limits_{c-i\infty}^{c+i\infty}
P_{K,d,s}^{+}(\tau)
R(\tau) y^{\tau} \zeta\left( 1-s+\tau , \lambda^{-d}\right)\,{\rm d}\tau
+O_{K,j}\left(  T^{-(j-K-1)\varepsilon} y^{3/2}\right)\;,\eqno(3.85)$$
for each integer $j$ that satisfies $j\geq K+2$.  
With regard to the $O$-term in (3.85), 
observe that if $A\in[0,\infty)$, and if one  
puts $j=[(A+3)/\varepsilon]+K+2$, then by (3.45), (3.47) and (3.60)  
it follows that 
$$ T^{-(j-K-1)\varepsilon} y^{3/2} 
\ll y^{\sigma -1/2} T^{3-(j-K-1)\varepsilon} 
\leq y^{\sigma -1/2} T^{-A}\;.\eqno(3.86)$$ 

By (3.85), (3.86), (3.81) and (3.67), one now has (when $0\leq A<\infty$)  
$$\eqalign{ 
J_{\varepsilon ,K}^{+}(d,s)
 &=O_{A, K, \varepsilon}\left( y^{\sigma -1/2} T^{-A}\right)
+\sum_{k=0}^{K}{a_{k}(d,s)\over 2\pi i}
\,\int\limits_{c-i\infty}^{c+i\infty} \tau^k
R(\tau) y^{\tau} \zeta\left( 1-s+\tau , \lambda^{-d}\right)\,{\rm d}\tau \cr 
 &=O_{A, K, \varepsilon}\left( y^{\sigma -1/2} T^{-A}\right)
+\sum_{k=0}^{K}(-1)^k a_{k}(d,s)
\sum_{n=1}^{\infty}\rho_k\left({n\over y}\right)
{\delta\left( \Lambda^{-d} , n\right)\over n^{1-s}}\;,}$$
where $a_{0}(d,s)=1$, while $\rho_k(u)$ is as in (3.66).
Therefore, with $\widetilde\rho_{K}(u)$ defined as in (3.49), 
one has  
$$J_{\varepsilon ,K}^{+}(d,s)
=O_{A, K, \varepsilon}\left( y^{\sigma -1/2} T^{-A}\right)
+\sum_{n=1}^{\infty}\widetilde\rho_K\left({n\over y}\right)
\delta\left( \Lambda^{-d} , n\right) n^{s-1}\quad\, 
\hbox{for all $\,A\in[0,\infty)$.}\eqno(3.87)$$

By (3.72), (3.83) and (3.87), the term $E_K\bigl( s , \lambda^d\bigr)$ in (3.48) satisfies
$$\eqalign{E_K\left( s , \lambda^d\right)
 &=(\pi /4)\delta_{d,0}\left(
{R(s-1)\over x^{s-1}} -{e^{-|s-1|}\over s-1}\right) + \cr
 &\quad +X_{d}(s)\left(
J_{\varepsilon , K}^{-}(d,s)
+O_{A,K,\varepsilon}\left(
y^{\sigma -1/2} T^{-A}\right)\right)
+O_{A,\varepsilon}\left(
x^{1/2-\sigma} T^{-A}\right)\;,}\eqno(3.88)$$
where the constant $A\in[0,\infty)$ is arbitrary.
If $d\neq 0$, then $\delta_{d,0}=0$, whereas if $d=0$ then 
$\delta_{d,0}=1$ and, by (3.64) and (3.79), and (3.45)-(3.47), one has 
$$\left| {R(s-1)\over x^{s-1}}\right| + \left| {e^{-|s-1|}\over s-1}\right|
\leq O_j\left( {b^{4/3}\over T^{j\varepsilon}}\,x^{1-\sigma}\right)
+{(j-1)!\over T^{j\varepsilon}}
\ll_{A,\varepsilon} x^{1/2-\sigma} T^{1-j\varepsilon} 
\leq x^{1/2-\sigma} T^{-A}\;,$$ 
for each pair $(j,A)\in{\Bbb Z}\times[0,\infty)$ such that 
$j=[(A+1)/\varepsilon]+1$.  
In either of the cases just mentioned it follows by (3.45)-(3.47), (3.60), and (3.6) of Lemma~5,
that the bound (3.10) of Lemma~7 gives
$$X_{d}(s)\ll T^{1/2-\sigma}=(xy)^{1/2-\sigma}\;.$$
Therefore, and by (3.88), one has:
$$E_K\left( s , \lambda^d\right)
\ll (xy)^{1/2 -\sigma} J_{\varepsilon , K}^{-}(d,s)
+O_{A,K,\varepsilon}\left( x^{1/2-\sigma} T^{-A}\right)\;,$$
for an arbitrary constant $A\in[0,\infty)$.
Since one may choose to put $A=\beta_K$ here (with $\beta_k$ as in (3.50)-(3.52)),
and since it was found that (3.71) holds, the bound (3.50) will now follow if
it can be proved that
$$J_{\varepsilon , K}^{-}(d,s)
\ll_{K,\varepsilon} y^{\sigma -1/2}
\left(\left( |t|/T\right)^{\alpha_K} +T^{-\beta_K}\right)
\int\limits_{-T^{\varepsilon}}^{T^{\varepsilon}}
\left|\zeta\left( 1/2 +i(t+v) , \lambda^d\right)\right|\,{{\rm d}v\over 1+v^2}
\;,\eqno(3.89)$$
with $\alpha_K$ and $\beta_K$ as in (3.52).

To obtain (3.89) consider first the analytic function $P_{K,d,s}^{-}(\tau)$
occurring in (3.82) and (3.84). For $\tau$ satisfying the conditions (3.77), 
the inequalities in (3.78) hold, so it follows by
(3.82) and Corollary~9 (the bound (3.38) in particular) that
$P_{K,d,s}^{-}(\tau)$ is analytic at all points of the rectangular subset of ${\Bbb C}$
given by (3.77), and on that subset satisfies:
$$P_{K,d,s}^{-}(\tau)
\ll_K \left( \left( |t|/T\right)^{(K+1)/2}+T^{-(K+1)/3}\right) |\tau|^{K+1}\eqno(3.90)$$
(having, in particular, a zero of order $K+1\geq 1$ at $\tau =0$). By this and (3.79), the
integrand in the `$-$'  case of (3.84) has no poles either on, or inside of, 
the rectangle given by (3.77). Therefore 
$$2\pi i J_{\varepsilon , K}^{-}(d,s)
={\cal J}\left( \overline{\kappa} , \overline{\mu}\right)
+{\cal J}\left( \overline{\mu} , \mu\right)
+{\cal J}\left( \mu , \kappa\right)
\quad\hbox{for $\kappa=c+iT^{\varepsilon}\,$,\ 
$\,\mu=\sigma -1/2+iT^{\varepsilon}\,$,}\eqno(3.91)$$
where $c=3/2$ and, for the relevant pairs of points, $\nu , \omega \in {\Bbb C}$,
$${\cal J}\left( \nu , \omega\right)
=\int\limits_{\nu}^{\omega} P_{K,d,s}^{-}(\tau) R(\tau) y^{\tau}
\zeta\left( 1-s+\tau , \lambda^{d}\right)
{\rm d}\tau\eqno(3.92)$$
with the contour of integration being the line segment between $\nu$ and $\omega$.

By (3.8) and (3.45), the bound (3.90) (valid in the rectangle (3.77)) 
certainly implies that one has 
$$P_{K,d,s}^{-}(\tau)\ll_K T^{-(K+1)/4}|\tau|^{K+1}\ll_K T^{-(1/4-\varepsilon)(K+1)}\ll_K 1$$
in the integrands of ${\cal J}\left( \overline{\kappa} , \overline{\mu}\right)$ and 
${\cal J}\left( \mu , \kappa\right)$ in (3.91)-(3.92).
Moreover, in those integrands 
$\Re (1-s+\tau)\geq 1/2$ and, by (3.79) and (3.8), 
$\delta_{d,0}T^{\varepsilon}<\left|\Im (1-s+\tau)\right|\leq |t|+T^{\varepsilon}\ll |t|+|d|$
(where $\delta_{d,0}$ is as in (3.68)), with it therefore following by the estimate
(3.39) of Lemma~10 (and (3.8), again) that one has there:
$$\zeta\left( 1-s+\tau , \lambda^{d}\right)\ll_{\varepsilon}
\left( 1+|t|+|d|\right)^{1/2+2\varepsilon}\asymp T^{1/4+\varepsilon}\;.$$
Given (3.45) and (3.47), the bound (3.64) (for $R(z)$\/) and points 
noted after (3.92) imply that 
$$\max\left\{
{\cal J}\left( \overline{\kappa} , \overline{\mu}\right) , 
{\cal J}\left( \mu , \kappa\right)\right\}
\ll_{K,\varepsilon ,j} (2-\sigma) (by)^{3/2}T^{ -j\varepsilon+1/4+\varepsilon}
\ll y^{3/2} T^{1/4 -(j-1)\varepsilon} \quad
\hbox{for $j\in{\Bbb N}\,$.}\eqno(3.93)$$

A suitable bound for the integral 
${\cal J}\bigl(\overline{\mu} ,\mu\bigr)$ in (3.91)-(3.92)  
may be obtained by applying the bound (3.90) for the factor $P_{K,d,s}^{-}(\tau)$
in (3.92), while using the case $j=K+3$ of (3.64) as a bound for $R(\tau)$ in (3.92);        
by following this with an application of the substitution $\tau =\sigma - 1/2 -iv$    
and an application of the identity  
$|\zeta( z , \lambda^d)| =|\zeta(\overline{z} , \lambda^d)|\,$ 
(which is a corollary of (2.4)), one finds, in particular, that     
$${\cal J}\left(\overline{\mu} , \mu\right)
\ll_K
y^{\sigma -1/2}\left( \left( {|t|\over T}\right)^{(K+1)/2} +T^{-(K+1)/3}\right)
\int\limits_{-T^{\varepsilon}}^{T^{\varepsilon}}
{\left|\zeta\left( 1/2 +i(t+v) , \lambda^{d}\right)\right| b^{|\sigma -1/2|} \over\left( 1+\left|\sigma - 1/2 +iv\right|\right)^2}\,
{\rm d}v\;.\eqno(3.94)$$

If $K\geq 3$, then, given (3.52) and (3.71), the desired bound for
$J_{\varepsilon , K}^{-}(d,s)$ (which is (3.89)) follows from (3.91), (3.93)
and (3.94)  on choosing $j$, in (3.93), to be large enough to ensure that
$$ y^{3/2} T^{1/4-(j-1)\varepsilon}\ll y^{\sigma -1/2} T^{-(K+1)/3}$$
(by (3.45) and (3.47), the choice
$j=[(K+9)/(3\varepsilon)]+2$ achieves this). Consequently (given  the discussion
which preceded (3.89)) the bound (3.50) of the lemma has now been shown
to hold when $K\geq 3$. 

In those cases in which the estimate (3.50) is, as yet, unproved 
(that is, the cases in which one has $K\in\{ 0 , 1 , 2\}$), 
one can improve on the bound (3.94) for ${\cal J}\bigl(\overline{\mu} , \mu\bigr)$
by combining the estimate
$$P_{K,d,s}^{-}(\tau)
=\sum_{k=K+1}^{3} a_{k}(d,s) \tau^{k}
+O\left( \left( \left( |t|/T\right)^2 +T^{-4/3}\right) |\tau|^4\right)$$
(from (3.82) and the case $K=3$ of (3.90)) with the estimates
for $a_1,a_2,a_3$ in (3.20)-(3.22) of Lemma~8.
These combined estimates show that if $\tau =\sigma - 1/2  +iv$
with $-T^{\varepsilon}\leq v\leq T^{\varepsilon}<T^{1/4}$, then, for
$K=0,1,2$,  one has 
$$P_{K,d,s}^{-}(\tau)
=O\left( {|\tau|^{K+1}+|\tau|^3\over T}+\left({|t|\over T}\right)^{\!\!2}\!|\tau|^4 \right)
+\cases{ 0 &if $K=2\,$, \cr it\bigl| z_{d,t}\bigr|^{-2}\tau^2 &if $K=1\,$, \cr
-it\bigl| z_{d,t}\bigr|^{-2}|\tau|^2 &if $K=0\,$,}$$ 
so that  
$$P_{K,d,s}^{-}(\tau)  
\ll\left( \left( |t|/T\right)^{\alpha_K}+T^{-\beta_K}\right) |\tau|^{K+1}
\left( 1+|\tau|\right)^2\qquad\qquad\hbox{($K=0,1,2$),}\eqno(3.95)$$
with $\alpha_K$, $\beta_K$ as in (3.52)
(note that $z_{d,t}=2|d|+ 1/2 +it$ here, so that by (3.8) one has
$\bigl| z_{d,t}\bigr|^2\asymp T$ and $|t|/\bigl| z_{d,t}\bigr|^2\asymp |t|/T\ll T^{-1/2}$, 
while $|\tau|^2 =(\sigma -1/2)^{2}+v^2<1+ T^{2\varepsilon}\ll T^{1/2}$\/).

By using the bound (3.95) and the case $j=K+5$ of (3.64) to estimate the
integral ${\cal J}\bigl(\overline{\mu} , \mu\bigr)$ in (3.91)-(3.92), one 
improves on the bound (3.94) to the extent that the exponents $(K+1)/2$ and $-(K+1)/3$ 
(in the factor $(|t|/T)^{(K+1)/2}+T^{-(K+1)/3}$ on the right-hand side of (3.94))
are sharpened  to $\alpha_K$ and $-\beta_K$, respectively. 
By reasoning similar to that employed in the paragraph below (3.94) it now follows
(given the improvement of (3.94) just obtained) that (3.89) holds for
each $K\in\{ 0 , 1 , 2\}$, and that the result (3.50) of the lemma has, as a consequence, now been proved
for these cases, which complement the case $K\geq 3$ established earlier.

The proof may now be completed with the help of   
Corollary~6 and the sub-convexity bound (3.41) of Kaufman.  
Indeed, by virtue of those results,  
the factor $|\zeta(1/2+i(t+v),\lambda^d)|$ occurring in the integrand in (3.50) is 
uniformly bounded above by $O_{\varepsilon}(T^{1/6+\varepsilon})$, and      
so, given that one has $\int_{-\infty}^{\infty}\bigl( 1+v^2\bigr)^{-1} {\rm d}v<\infty$,  
the estimate (3.51) follows from (3.50)\quad$\square$

\bigskip

\goodbreak\proclaim{\Smallcaps Corollary~12}. Let $\rho : (0,\infty)\rightarrow{\Bbb R}$ 
and $b\in(1,\infty)$ be as in
the above lemma. Suppose that $h\geq 1$, that $0<\eta\leq 1\leq T_{*}$ and that $\varepsilon >0$.
Then, provided that $T_{*}$ is sufficiently large (in terms of $b$ and $\eta$), one has
$$\zeta^{h}\left( 1/2+it , \lambda^d\right)
\ll {3^h\over\eta}\,\int\limits_{-2\eta}^{2\eta}\ 
\left|\sum_{n=1}^{\infty}\rho\left({n\over  e^{\theta} T_{*}^{1/2}}\right)\delta\left(\Lambda^d , n\right) n^{-1/2-it}\right|^{h}
{\rm d}\theta 
+O_{\rho,h,\varepsilon}\left( T_{*}^{(\varepsilon -1/3)h}\right) ,$$
for all $d\in{\Bbb Z}$ and all $t\in{\Bbb R}$ such that
$$e^{-\eta} T_{*}\leq \pi^{-2}|2d+it|^2\leq e^{\eta} T_{*}\;.\eqno(3.96)$$

\medskip

\goodbreak
\noindent{\Smallcaps Proof.}\ 
By (2.6), $\bigl| X_d(s)\bigr| =1$ if $\Re (s)=1/2$,
so, for $d,t$ as in (3.96), it follows from the case $\sigma =1/2$, $K=0$ of Lemma~11 that
$$\eqalign{
\left|\zeta\left( 1/2+it , \lambda^d\right)\right|
 &\leq\left|\sum_{n=1}^{\infty}\rho\left( {n\over x}\right)\delta\left( \Lambda^d ,n\right) n^{-1/2 -it}\right|
+\left|\sum_{n=1}^{\infty}\rho\left( {n\over y}\right) 
\delta\left( \Lambda^{-d} ,n\right) n^{-1/2 +it}\right| + \cr
 &\quad +O\left( T_{*}^{-1/2}\exp\left( -\pi\sqrt{T_{*} /e}\right)\right)
+O_{\rho,\varepsilon}\left( T^{1/6+\varepsilon}\left(\left( |t|/T\right)+T^{-1}\right)\right) ,}
\eqno(3.97)$$
were $T=T(d,t)$, and $x,y$ is any pair of real numbers satisfying (3.46) and (3.47).
By the condition (3.96), the result (3.7) of Lemma~5 and the hypothesis that $T_{*}$ is sufficiently large (in terms of
$\eta$), one may assume here that
$$e^{-2\eta} T_{*}\leq T\leq e^{2\eta} T_{*}\;.\eqno(3.98)$$
By this and (3.96), the $O$-terms in (3.97) are not greater than 
$O_{\rho,\varepsilon}\bigl( T_{*}^{\varepsilon -1/3}\bigr)$.
Moreover, using (3.98) one  can check that the condition (3.47) of Lemma~11 
is satisfied by all $y$ lying in the interval $\bigl[ e^{-\eta}\sqrt{T} , e^{\eta}\sqrt{T}\bigr]$
(given the hypotheses that $\eta\leq 1$ and that $T_{*}$ is sufficiently large in terms of $b$).
Since (3.46) implies $y^{-1}{\rm d}y=-x^{-1}{\rm d}x$, it follows by (3.98) and
the above discussion that, through an application of H\"{o}lder's inequality, followed
by integration with respect to $y$, one may deduce from (3.97) that
$$2\eta
\left|\zeta\left( 1/2+it , \lambda^d\right)\right|^h
\leq 3^{h-1}\sum_{w=\pm 1}\ \int\limits_{e^{-\eta}\sqrt{T}}^{e^{\eta}\sqrt{T}}\ 
\left|\sum_{n=1}^{\infty}\rho\left( {n\over x}\right)
\delta\left( \Lambda^{wd} ,n\right) n^{-1/2 -wit}\right|^{h} x^{-1}{\rm d}x
+O_{\rho,h,\varepsilon}\left( \eta T_{*}^{(\varepsilon -1/3) h}\right) .$$
Since $\rho(u)=\overline{\rho(u)}$ and 
$\delta\bigl(\Lambda^{-d} , n)=\overline{\delta\bigl(\Lambda^{d} , n)}$,
and since $|\overline{z}|=|z|\geq 0$ for $z\in{\Bbb C}$,
the corollary now follows by appealing to (3.98) and then making the 
substition $x=T_{*}^{1/2}e^{\theta}$\quad$\square$

\bigskip

\goodbreak\proclaim{\Smallcaps Remarks~13}. Suppose that $b>1$. Then one example of a function $\rho : (0,\infty)\rightarrow{\Bbb R}$ 
satisfying the hypotheses of Lemma~11 is that which is given by: 
$$\rho(u)=\left(\ \int\limits_{-\infty}^{\infty}\Phi(t) {\rm d}t\right)^{\!\!-1} 
\int\limits_{(\log u)/(\log b)}^{\infty}\Phi(t) {\rm d}t\qquad\qquad 
\hbox{($u>0$),}\eqno(3.99)$$ 
with 
$$\Phi(t)=\cases{\exp\left( -\left( 1 - t^2\right)^{-1}\right) 
 &if $\,-1<t<1$, \cr 0 &otherwise.}\eqno(3.100)$$ 
We leave it to the reader to verify that this  
$\rho$ is indeed an infinitely differentiable function on $(0,\infty)$ 
satisfying both of the conditions (3.43) and (3.44). 
Since the function $\Phi$ defined in (3.100) takes only non-negative real values, 
it follows from (3.99) and (3.100) that the range of the function $\rho$ is 
the interval $[0,1]$.  
It can moreover be shown (by induction) that this function $\rho$ 
is such that 
$$\rho^{(j)}(u)=u^{-j}\sum_{k=0}^{j-1} r(j,k) (\log b)^{-(k+1)} 
\Phi^{(k)}\!\left( {\log u\over\log b}\right)\qquad\qquad 
\hbox{($j\in{\Bbb N}$ and $u>0$),}$$ 
where each coefficient $r(j,k)$ is an integer that is completely determined 
by the pair $j,k$. 
Given the definition (3.100), one can furthermore establish   
that $\sup\{ |\Phi^{(k)}(t)| : t\in{\Bbb R}\} <\infty$ for all $k\in{\Bbb N}\cup\{ 0\}$. 
It therefore follows that if $1<b\leq e$, and if 
$\rho : (0,\infty)\rightarrow{\Bbb R}$ is given by (3.99) and (3.100), then one has 
$$\rho^{(j)}(u)\ll_j (u\log b)^{-j}\qquad\qquad 
\hbox{($j\in{\Bbb N}\cup\{ 0\}$ and $u>0$).}\eqno(3.101)$$

\bigskip

\goodbreak 
\noindent{\SectionHeadingFont 4. Proof that Theorem 2 implies Theorem 1}

\medskip

\noindent In this section we show that Theorem 1 is a corollary of Theorem 2:
the proof of the latter theorem is divided between Sections~6 and~7.
\par 
Throughout this section we assume that the hypotheses of Theorem~1 
concerning $\vartheta$ and $\varepsilon$ are satisfied. 
We assume also that  $\varepsilon$ satisfies $\varepsilon\leq 1$:  
no loss of generality results from this, for the 
cases of Theorem~1 in which one has $\varepsilon >1$ are a 
trivial corollary of the case in which one has $\varepsilon =1$. 
We suppose, furthermore, that 
$A$ is some complex valued function with domain 
${\frak O}-\{ 0\}$. 
For $M\geq 1$ and $D\geq M$, we put 
$$E^{\star}(D ; M , A) = E(D ; M , A) - E\!\left( M ; M , A\right)\;,\eqno(4.1)$$ 
where $E(D;M,A)$ is as defined in Section~1. 
\par 
Given (1.2) and (1.4), it follows by the Cauchy-Schwarz inequality and 
elementary estimates that 
$$\left| P_M\left( A ; it , \lambda^d\right)\right|^2 
\leq 8M\sum_{0<|\mu|^2\leq M}\left| A(\mu)\right|^2 
\leq 64 M^2 \max_{0<|\mu|^2\leq M} \left| A(\mu)\right|^2\qquad\quad 
\hbox{($t\in{\Bbb R}$, $\,d\in{\Bbb Z}$, $\,M\geq 1$).}\eqno(4.2)$$
By the definition (1.3)-(1.4) and the first part of (4.2), we find that 
$$E(D;M,A)\ll \biggl( M\sum_{0<|\mu|^2\leq M}\left| A(\mu)\right|^2 \biggr) 
E(D;1,U)\qquad\qquad\hbox{($D,M\geq 1$),}\eqno(4.3)$$ 
where $U$ is the complex function defined on ${\frak O}-\{ 0\}$ by: 
$$U(\mu) =\cases{1 &if $\,\mu =1$, \cr 0 &otherwise.}\eqno(4.4)$$ 
\par 
We are going to show that Theorem~2 implies that, for $M\geq 1$ and $D\geq M$, one has both  
$$E^{\star}(D;M,A)\ll_{\varepsilon} \left( D^{2+\varepsilon} 
+ \bigl( 1 + D M^{-3/2}\bigr)^{\vartheta} 
D^{1+\varepsilon} M^2\right)\sum_{0<|\mu|^2\leq M} \left| A(\mu )\right|^2\eqno(4.5)$$
and 
$$E^{\star}(D;M,A)\ll_{\varepsilon} D^{2+\varepsilon}
\biggl(\ \sum_{0<|\mu|^2\leq M} \left| A(\mu )\right|^2\biggr) 
+\left( 1 + D M^{-2}\right)^{\vartheta} 
D^{1+\varepsilon} M^3 \max_{0<|\mu|^2\leq M} | A(\mu )|^2\;.\eqno(4.6)$$
Before proceeding to the proof of (4.5) and (4.6), we consider 
the implications of these bounds. 
\par 
If (4.5) can be shown to hold for 
all $M\geq 1$ and all $D\geq M$, then (since we may substitute 
$U$, as defined in (4.4), in place of the function $A$)  
it will follow by (4.1), (4.5), (4.4), (1.3), (1.4), (1.2) and the bound 
(1.13) of Kim and Shahidi that, for $D\geq 1$, one has 
$$E\left( D; 1, U\right) 
=E^{\star}\!\left( D; 1, U\right) +E\left( 1 ; 1 , U\right) 
=\left( O_{\varepsilon}\left( D^{2+\varepsilon}\right) 
+ O(1)\right)  |U(1)|^2 
\ll_{\varepsilon} D^{2+\varepsilon}\eqno(4.7)$$ 
(note that the implied bound $E\left( D; 1, U\right) \ll_{\varepsilon} D^{2+\varepsilon}$ 
is not a new, for it may 
also be obtained directly from Sarnak's sharper estimate in (1.16)). 
Moreover, given (4.7), it follows by (4.1) and (4.3) that, for $M\geq 1$ and 
$D\geq M$, one has 
$$\eqalign{ 
E(D;M,A) &=E^{\star}(D;M,A) + E\left( M; M,A\right) = \cr 
 &= E^{\star}(D;M,A) + O\Biggr( \biggl( M\sum_{0<|\mu|^2\leq M} 
\left| A(\mu)\right|^2 \biggr) 
E\left( M ;1,U\right)\Biggr) = \cr 
 &= E^{\star}(D;M,A) + O_{\varepsilon}\Biggl( M^{3+\varepsilon} 
\sum_{0<|\mu|^2\leq M}\left| A(\mu)\right|^2 \Biggr) \;. 
}$$
Therefore, given that we also have the second inequality of (4.2), 
we may conclude that 
the bounds asserted in (4.5) and (4.6) would 
imply that the results (1.14) and (1.15) of Theorem~1 are valid 
whenever one has both $M\geq 1$ and $D\geq M$. 
We observe also that, by (4.3), the corollary of (4.5) noted in (4.7), 
and the second inequality of (4.2), it follows that 
$$E(D;M,A)\ll_{\varepsilon} D^{2+\varepsilon} M \sum_{0<|\mu|^2\leq M} 
|A(\mu)|^2
\ll D^{2+\varepsilon} M^2 \max_{0<|\mu|^2\leq M} |A(\mu)|^2
\qquad\qquad\hbox{($D,M\geq 1$).}$$ 
Consequently we may also conclude that the bound asserted in (4.5)  
would imply that the results (1.14) and (1.15) of Theorem~1 are 
valid whenever one has both $M\geq 1$ and $1\leq D<M$. 
This last conclusion, when combined with our earlier conclusion 
(in repect of the cases where $M\geq 1$ and $D\geq M$),  
is enough to show that Theorem~1 will follow if it can 
be proved that the bounds (4.5) and (4.6) hold whenever one has $M\geq 1$ and 
$D\geq M$. 
\par 
Our task now is to show, as was promised, that 
Theorem~2 implies the validity of (4.5) and (4.6) for 
all $M\geq 1$ and all $D\geq M$; 
it follows from our conclusions   
in the preceding paragraph that the succesful  
completion of this task will be enough to show   
that Theorem~1 is a corollary of Theorem~2.   
\par 
We assume henceforth that $M$ and $D$ are given real numbers satisfying   
$$1\leq M\leq D\;.\eqno(4.8)$$ 
Note that it is a consequence of the definitions given in (1.3), (1.4) and (4.1) 
that all terms occurring on either side of the bounds stated in 
(4.5) and (4.6) are independent of the mapping that is the restriction of $A$ 
to the set $\{\mu\in{\frak O} : |\mu|^2>M\}$. Therefore we may furthermore assume  
that $A(\mu)=0$ for all $\mu\in{\frak O}$ such that $|\mu|^2 >M$. 
This justifies our subsequent use of the convenient notation 
`$\| A\|_2^2$' and `$\| A\|_{\infty}^2$' to signify, respectively, the sum over $\mu$ and 
the maximum over $\mu$ occurring in (4.5)-(4.6) (i.e. our use of this notation 
will be consistent with our definition of it in (1.17)).  
\par 
By (4.8), (4.1) and (1.3), it follows that 
$$E^{\star}(D ; M , A)
\leq\sum_{d=-\infty}^{\infty}\ \int\limits_{-\infty}^{\infty}
\left|\zeta\left( 1/2 +it,\lambda^d\right)\right|^4
\left| P_M\left( A;it,\lambda^d\right)\right|^2 
\chi_{[M ,\sqrt{5}D]}\left( |2d+it|\right) 
{\rm d}t\;,\eqno(4.9)$$
where $\chi_{[a,b]}(x)=|\{x\}\cap [a,b]|\/$. 

We now put 
$$\eta =(\log 2)/3\qquad{\rm and}\qquad b=e^{\eta /3}=2^{1/9}\;,\eqno(4.10)$$ 
and we choose an infinitely differentiable function 
$W : [0,\infty)\rightarrow[0,1]$ such that
$$W(x)=\cases{1 &if $0\leq x\leq 1\/$,\cr 0 &if $x\geq e^{\eta}\/$,}\eqno(4.11)$$
and
$$W^{(j)}(x)\ll_j \ x^{-j}\qquad\hbox{for $\,x>0\,$ and $\,j=0,1,2,\ldots\ \;$}\eqno(4.12)$$ 
(an example of such a function is the mapping $x\mapsto\rho(x/b)$, in which 
we take $\rho$ to be as defined in (3.99) and~(3.100)). 
Then, on taking
$$L=1+\left[ 2\eta^{-1}\log\!\left( {\sqrt{5}D\over M}\right)\right]\eqno(4.13)$$
(which, by virtue of (4.8) and (4.10), makes $L\geq 7$), one has
$$W\left({x^2\over 5D^2}\right)
-W\left({x^2\over 5D^2 e^{-(L+1)\eta}}\right)
\geq\chi_{[M ,\sqrt{5}D]}\left( x\right) \quad\hbox{for $x\geq 0\/$,}$$
so it follows from (4.9) that
$$E^{\star}(D ; M , A)\leq\sum_{{\ell}=0}^L \widetilde{E}_{\ell} (D ; M , A)\;,\eqno(4.14)$$
where
$$\widetilde{E}_{\ell} (D ; M , A)
=\sum_{d=-\infty}^{\infty}\ \int\limits_{-\infty}^{\infty}
\left|\zeta\left( 1/2 +it,\lambda^d\right)\right|^4
\left| P_M\left( A;it,\lambda^d\right)\right|^2 
w_0\!\!\left({|2d+it|^2\over 5D^2 e^{-{\ell}\eta}}\right) {\rm d}t\eqno(4.15)$$
with
$$w_0(x)=W(x)-W\left( e^{\eta} x\right)\quad\hbox{for $x\geq 0\/$.}$$
Given (4.10), the conditions on $W(x)$, including (4.11) and (4.12), are enough to ensure that
$w_0(x)$ is an infinitely differentiable function from $[0,\infty)$ into $[0,1]$, 
and satisfies, for $x>0$ and $j=0,1,2,\ldots\ \;$,
$$w_0^{(j)}(x)=\cases{O_j\left( (\eta x)^{-j}\right) &for $e^{-\eta}\leq x\leq e^{\eta}\;$, \cr
0 &otherwise.}\eqno(4.16)$$

For the purpose of achieving our goals in this section, 
it will be enough that we obtain suitable bounds for 
each sum $\widetilde{E}_{\ell} (D ; M , A)$ occurring in (4.14).  
In accordance with this latter objective, 
we suppose now that that ${\ell}$ is one of the integers 
in the set $\{ 0,1,2,\ldots , L\}$, 
where $L$ is as defined in (4.13).
We also put 
$$T_{\ell} = 5\pi^{-2} D^2 e^{-{\ell}\eta}\;,\eqno(4.17)$$ 
and we define $\rho$ to be the mapping from $(0,\infty)$ to ${\Bbb R}$ 
that is given by (3.99) and (3.100) (with $b$ as specified in (4.10)).  
As we observed in Remarks~13, the function $\rho : (0,\infty)\rightarrow{\Bbb R}$ 
satisfies the conditions (3.43) and (3.44) of Lemma~11 and is infinitely 
differentiable on $(0,\infty)$. By the case $b=2^{1/9}$ of (3.101), one has, moreover,    
$$\rho^{(j)}(u)\ll_j u^{-j}\qquad\quad\hbox{ for all $u>0$ and all $j\in{\Bbb N}\cup\{ 0\}$.} 
\eqno(4.18)$$ 
Since ${\ell}\leq L$, it follows by (4.17), (4.13) and (4.10) 
that we have now 
$$T_{\ell}\geq 2^{-1/3} \pi^{-2} M^2\;.\eqno(4.19)$$ 
\par 
Given (4.10), (4.15)-(4.17) and our specific definition (in absolute terms) of 
the function $\rho$, 
it follows by the case $\varepsilon =1/12$, $h=4$, 
$T_{*}=T_{\ell}$ of Corollary~12 that,  
if $T_{\ell}$ is greater than a certain positive absolute constant, $B_0$ (say), then 
one has a bound of the form 
$$\widetilde{E}_{\ell}(D ; M , A)
\ll {1\over\eta}\,\int\limits_{-2\eta}^{2\eta}
\breve{E}\left( T_{\ell} ; M , A ; \theta\right) {\rm d}\theta
+ {1\over T_{\ell}}\,\sum_{d=-\infty}^{\infty}\ \int\limits_{-\infty}^{\infty}
\left| P_M\left( A;it,\lambda^d\right)\right|^2 
w_0\!\!\left({|2d+it|^2\over\pi^2 T_{\ell}}\right) {\rm d}t\;,\eqno(4.20)$$
where
$$\breve{E}\left( T_{\ell} ; M , A ; \theta\right)
=\sum_{d=-\infty}^{\infty}\ \int\limits_{-\infty}^{\infty}
\left|\sum_{n=1}^{\infty} {\rho\bigl(\,n e^{-\theta} T_{\ell}^{-1/2}\,\bigr)
\delta\left(\Lambda^d , n\right)\over n^{1/2+it}}\right|^{4}
\left| P_M\left( A;it,\lambda^d\right)\right|^2 
w_0\!\!\left({|2d+it|^2\over\pi^2 T_{\ell}}\right) {\rm d}t\;,\eqno(4.21)$$ 
while the implicit constant associated with the `$\ll$' notation   
in (4.20) is absolute.  
We observe (separately) that, since $w_0(x)\geq 0$ for all $x\in{\Bbb R}$, 
and since each mapping $t\mapsto\zeta(1/2+it,\lambda^d)$ is bounded on any 
bounded real interval, it follows from the definition (4.21), in combination with 
(4.15)-(4.17) and (4.10), that 
$$\breve E\left( T_{\ell} ; M , A ; \theta\right)\geq 0\qquad\qquad 
\hbox{($-2\eta\leq\theta\leq 2\eta$),}\eqno(4.22)$$ 
and that, for all $B\in(0,\infty)$, one has: 
$${1\over T_{\ell}}\,\sum_{d=-\infty}^{\infty}\ \int\limits_{-\infty}^{\infty}
\left| P_M\left( A;it,\lambda^d\right)\right|^2 
w_0\!\!\left({|2d+it|^2\over\pi^2 T_{\ell}}\right) {\rm d}t  
\gg_B\,\widetilde E_{\ell}(D;M,A)\qquad\quad 
\hbox{if $\,T_{\ell}\leq B$.}\eqno(4.23)$$ 
We note, in particular, that if $T_{\ell}$ is not greater than $B_0$  
(the positive absolute constant 
just referred to, above (4.20)), then one has $T_{\ell}\leq B_0$, 
and so, by virtue of (4.22) and the case 
$B=B_0$ of (4.23), the upper bound stated in (4.20)-(4.21) is valid.   
We may therefore conclude that, regardless of whether or not 
the condition $T_{\ell}>B_0$ is satisfied, one does have 
the upper bound stated in (4.20)-(4.21). 
\par 
Given (4.16), one may deduce from (4.20), (4.2) and (4.19) that
$$\widetilde{E}_{\ell} (D ; M , A)
\ll M \| A\|_2^2 
+ \max_{-2\eta\leq\theta\leq 2\eta} 
\breve{E}\left( T_{{\ell}}; M , A ; \theta\right)  
= O\left( T_{\ell}^{1/2} \| A\|_2^2\right) 
+ \max_{-2\eta\leq\theta\leq 2\eta} 
\breve{E}\left( T_{{\ell}}; M , A ; \theta\right) .\eqno(4.24)$$
By virtue of the conclusions reached in (4.14), (4.19) and (4.24), it will now
suffice (for the estimation of $E^{\star}(D;M,A)$) 
that, in the cases where (4.19) holds and $\theta$ satisfies 
$$-2\eta\leq\theta\leq 2\eta\;,\eqno(4.25)$$ 
we obtain suitable bounds for $\breve{E}(T_{\ell} ; M , A ; \theta)$. 
Given that the function $\rho$ satisfies the condition (3.44) of Lemma~11, 
and bearing in mind the definition (4.21) of $\breve E(T_{\ell};M,A;\theta)$, 
we either have 
$$e^{\theta} T_{\ell}^{1/2} > b^{-1}\;,\eqno(4.26)$$ 
or else have $\breve E(T_{\ell};M,A;\theta)=0$. 
Therefore, 
in our subsequent discussion of $\breve E(T_{\ell};M,A;\theta)$, 
we may assume that $\theta$ and $T_{\ell}$ satisfy the conditions (4.25) and (4.26).  
We remark that 
it is only when one has $M\ll 1$ that (4.26) is not implied by (4.10), (4.19) and (4.25). 
\par 
In order that we may make use of Theorem~2 in bounding 
$\breve{E}(T_{\ell} ; M ,A ; \theta)$, 
both the sum over $n$ occurring in (4.21) and the sum $P_M(A;it,\lambda^d)$ 
defined in (1.4) must first be split up into smaller subsums. 
In preparation for the first part of this splitting procedure, we put 
$$r(u)=\rho(u)-\rho(bu)\quad\hbox{for $u>0\;$,}\eqno(4.27)$$
and we set 
$$L_1=\left[\left(\theta+{1\over 2}\log T_{\ell}\right)\!/\log b\right] +1\eqno(4.28)$$
(so that, by (4.26) and (4.10), one has $L_1\geq 0$). 
Then, by (3.44),
$$\rho(u)=\sum_{{\ell}_1=0}^{L_1} r\left( b^{{\ell}_1} u\right)\quad
\hbox{if $u\geq e^{-\theta} T_{\ell}^{-1/2}\;$,}$$
and so it follows by H\"{o}lder's inequality that, 
for all $d\in{\Bbb Z}$ and all $t\in{\Bbb R}$, one has   
$$\left|\sum_{n=1}^{\infty}\rho\!\left({n\over  e^{\theta} T_{\ell}^{1/2}}\right)
\delta\left(\Lambda^d , n\right) n^{-1/2-it}\right|^{4}
\leq \left( 1 + L_1\right)^3\sum_{{\ell}_1=0}^{L_1}\left|
\sum_{n=1}^{\infty} r\!\left({n\over  e^{\theta} T_{\ell}^{1/2}b^{-{\ell}_1}}\right)
\delta\left(\Lambda^d , n\right) n^{-1/2-it}\right|^{4}\;.\eqno(4.29)$$
Likewise, by (1.4), (4.8), (4.10) and the Cauchy-Schwarz inequality, one has 
$$\left| P_M\left( A ; it , \lambda^d\right)\right|^2
\leq (1 + L_2)\sum_{{\ell}_2=0}^{L_2}
\left|\sum_{M b^{-{\ell}_2-1}<|\mu|^2\leq M b^{-{\ell}_2}} 
A(\mu)\Lambda^{d}(\mu) |\mu|^{-2it}\right|^2\qquad\quad\,  
\hbox{($d\in{\Bbb Z}$, $t\in{\Bbb R}$),}\eqno(4.30)$$
where
$$L_2 =\left[ (\log M)/\log b\right]\geq 0\;.\eqno(4.31)$$ 
\par 
Given (4.10), (4.19), 
(4.21), (4.25), (4.26) and the definition (2.10) of $\delta\bigl(\Lambda^d ,n\bigr)$, 
it follows by the bounds (4.28)-(4.29) and (4.30)-(4.31)  that  
for some pair $K,M'>0$ satisfying
$${1\over b} < K\leq e^{\theta} T_{\ell}^{1/2}\ll T_{\ell}^{1/2}\quad\hbox{and}\quad 
1\leq M'\leq M\ll T_{\ell}^{1/2}\;,\eqno(4.32)$$
one has
$$\eqalign{ 
\breve{E}\left( T_{\ell} ; M , A ; \theta\right)
 &\ll\left( 1+\left|\log T_{\ell}\right|\right)^6 K^{-2}
\sum_{d=-\infty}^{\infty}\ \int\limits_{-\infty}^{\infty}
\left|\,\sum_{\kappa\neq 0}w\!\left(|\kappa|^2\right) 
\Lambda^{d}(\kappa) |\kappa|^{-2it}\right|^{4} \times \cr 
 &\qquad\times\left|\sum_{M'/b<|\mu|^2\leq M'} A(\mu) 
\Lambda^{d}(\mu) |\mu|^{-2it}\right|^2
w_0\!\!\left({|2d+it|^2\over\pi^2 T_{\ell}}\right) {\rm d}t ,}\eqno(4.33)$$
with
$$w(x)=(x/K)^{-1/2} r(x/K)\quad\hbox{for $x>0\;$.}\eqno(4.34)$$ 
\par 
Since the function $\rho$ satisfies the conditions (3.43) and (3.44) 
of Lemma~11, it follows by (4.27) that 
the support of the function 
$r : (0,\infty)\rightarrow{\Bbb R}$ is contained in the interval 
$\bigl[ b^{-2} , b\bigr]$; 
we therefore find, by (4.34) and (4.10), that the support of 
the function $w : (0,\infty) \rightarrow{\Bbb R}$ is contained in 
the interval $[e^{-\eta}bK , bK]$. Moreover, by (4.34), (4.27), (4.18) and (4.10), 
the function $w$ is infinitely differentiable on $(0,\infty)$ and satisfies 
$w^{(j)}(x)\ll_j (K/x)^{1/2} x^{-j}$ for all $x>0$ and all $j\in{\Bbb N}\cup\{ 0\}$. 
Consequently (and since $\eta$ is, here, the positive absolute constant $(\log 2)/3$)   
it follows that, for all $x>0$ and all $j\in{\Bbb N}\cup\{ 0\}$, one has:  
$$w^{(j)}(x) 
=\cases{ O_j\left( (\eta x)^{-j}\right) &if $\,e^{-\eta} bK\leq x\leq e^{\eta} bK$, \cr 
0 &otherwise.}\eqno(4.35)$$ 
\par 
By (4.35), (4.32), (4.16) and (4.10), it follows that if we put 
$K_0=1$, $K_1=K_2=bK$, $w_1=w_2=w$, $M_1=M'$, $T=T_{\ell}$ and, 
for $\mu\in{\frak O}-\{ 0\}$, 
$$a(\mu)= \cases{A(\mu) &if $\,M'/b<|\mu|^2\leq M'$, \cr 0 &otherwise,}\eqno(4.36)$$ 
and if (at the same time) we substitute $\varepsilon /33\in(0,1/33]$ for $\varepsilon$, 
then the hypotheses of Theorem~2 
(up to, and including, (1.20)) will be satisfied.   
Therefore, by Theorem~2 and the elementary estimates (1.33) and (1.35) (all applied 
with $\varepsilon /33$ in place of $\varepsilon$), 
we are able to deduce from (4.33)-(4.34) that  
$$\breve E\left( T_{\ell} ; M , A ; \theta\right) 
\ll\left( 1+\left|\log T_{\ell}\right|\right)^6 K^{-2} 
\left( 2\pi {\cal D}_0 + (\pi /2){\cal D}_1^{\star}+(\pi /2){\cal D}_2^{\star}  
+{\cal E}\right) ,\eqno(4.37)$$ 
where, with $a : {\frak O}-\{ 0\}\rightarrow{\Bbb C}$ defined as in (4.36), 
the  terms ${\cal D}_0$, ${\cal D}_1^{\star}$, ${\cal D}_2^{\star}$ and ${\cal E}$ 
satisfy:  
$$\max\left\{ \left| {\cal D}_0\right|\,,\,\left| {\cal D}_1^{\star}\right|\,,  
\,\left| {\cal D}_2^{\star}\right|\right\}  
\ll_{\varepsilon} T_{\ell}^{1+(\varepsilon /33)} K^2 \| a\|_2^2 
\ll T_{\ell}^{1+(\varepsilon /3)} K^2 \| A\|_2^2 \;,\eqno(4.38)$$ 
$$\eqalignno{ 
K^{-2} {\cal E} &\ll_{\varepsilon} \left( 
{(M')^2\over T_{\ell}^{1/2}} 
+\left( {(M')^{2-(3/2)\vartheta}\over T_{\ell}^{(1-\vartheta)/2}}\right) 
\left( {K\over T_{\ell}^{1/2}}\right)^{\!\!\vartheta} 
\right) T_{\ell}^{1+(3/11)\varepsilon} \| a\|_2^2 \ll \cr 
 &\ll\left( 1 + {T_{\ell}\over M^3}\right)^{\!\!\vartheta /2}  
\!\!\left( {M^4\over T_{\ell}} \right)^{\!\!1/2}  
T_{\ell}^{1+(\varepsilon /3)} \| A\|_2^2 &(4.39)}$$ 
and  
$$\eqalignno{  
K^{-2} {\cal E} &\ll_{\varepsilon} \left( 
{(M')^2\over T_{\ell}^{1/2}} 
+\left( {K\over T_{\ell}^{1/2}} 
+\left( {K\over T_{\ell}^{1/2} (M')^{1/2}}\right)^{\!\!\vartheta}\right) 
\left( {(M')^2\over T_{\ell}^{1/2}}\right)^{\!\!1-\vartheta} 
\right) T_{\ell}^{1+(\varepsilon /3)} M'\| a\|_{\infty}^2 \ll \cr 
 &\ll\left( 1 + {T_{\ell}\over M^4}\right)^{\!\!\vartheta /2}  
\!\! \left( {M^4\over T_{\ell}} \right)^{\!\!1/2}  
T_{\ell}^{1+(\varepsilon /3)} M \| A\|_{\infty}^2 &(4.40)}$$ 
(note that, in each of (4.38), (4.39) and (4.40), the final upper bound follows by virtue 
of (1.13), (4.36), the hypothesis that $\varepsilon$ lies in $(0,1]$ and the 
bounds on $K$, $M$ and $T_{\ell}$ that are implied by (4.32)). 
\par 
By (4.37)-(4.40) and (4.32), we obtain both  
$$\breve E\left( T_{\ell};M,A;\theta\right) 
\ll_{\varepsilon} 
\left( 1 + \left( 1 + {T_{\ell}\over M^3}\right)^{\!\!\vartheta /2} 
\!\!\left( {M^4\over T_{\ell}} \right)^{\!\!1/2}\right)   
T_{\ell}^{1+(\varepsilon /2)} \| A\|_2^2$$ 
and 
$$\breve E\left( T_{\ell};M,A;\theta\right) 
\ll_{\varepsilon} 
T_{\ell}^{1+(\varepsilon /2)} \| A\|_2^2 
+\left( 1 + {T_{\ell}\over M^4}\right)^{\!\!\vartheta /2}  
\!\!\left( {M^4\over T_{\ell}} \right)^{\!\!1/2}  
T_{\ell}^{1+(\varepsilon /2)} M \| A\|_{\infty}^2\;.$$ 
Moreover, for the reasons that are mentioned below (4.8), 
we may assume here that the mapping $A$ is such that 
$\| A\|_2^2=\sum_{0<|\mu|^2\leq M} |A(\mu)|^2$ and 
$\| A\|_{\infty} =\max\{ |A(\mu)| : \mu\in{\frak O}\ {\rm and}\ 0<|\mu|^2\leq M\}$. 
Therefore, by  (4.24), (4.17), (4.8) and the bounds  
just obtained for $\breve E(T_{\ell};M,A;\theta)$,  
it follows that one has both  
$$\widetilde E_{\ell}(D;M,A) 
\ll_{\varepsilon} 
\biggl( D^{1+\varepsilon} + 
\left( 1 + {D^2\over M^3}\right)^{\!\!\vartheta /2} 
\!\!D^{\varepsilon} M^2\biggr)   
T_{\ell}^{1/2} \sum_{0<|\mu|^2\leq M} |A(\mu)|^2\eqno(4.41)$$ 
and 
$$\widetilde E_{\ell}(D;M,A) 
\ll_{\varepsilon} 
\Biggl(\;\sum_{0<|\mu|^2\leq M} |A(\mu)|^2\Biggr) D^{1+\varepsilon} T_{\ell}^{1/2}   
+ \biggl(\,\max_{0<|\mu|^2\leq M} |A(\mu)|^2\biggr) 
\!\left( 1 + {D^2\over M^4}\right)^{\!\!\vartheta /2} \!M^3 D^{\varepsilon} 
T_{\ell}^{1/2}\;.\eqno(4.42)$$ 
\par 
Since our conclusions in the above paragraph are 
valid for any $\ell\in\{ 0,1,2,\ldots ,L\}$, and 
since the definitions (4.17) and (4.10) imply that 
we have $\sum_{\ell =0}^{\infty} T_{\ell}^{1/2} 
=(\sqrt{5} /\pi) D\sum_{\ell =0}^{\infty} 2^{-\ell /6}=O(D)$, it follows that 
the upper bounds (4.14), (4.41) and (4.42) imply the 
bounds for $E^{\star}(D;M,A)$ that are stated in (4.5) and (4.6).   
Moreover, our only assumption concerning $M$ and $D$ has been that the condition 
(4.8) is satisfied, and so we have completed the set task 
of showing that Theorem~2 implies the validity of (4.5) and (4.6) for 
all $M\geq 1$ and all $D\geq M$; in view of our conclusions 
in the paragraphs preceding (4.8), we have thereby shown that 
Theorem~2 implies Theorem~1\quad$\square$ 

\bigskip

\goodbreak 
\noindent{\SectionHeadingFont 5. Some lemmas that we need for the proof of Theorem~2}

\medskip 

\noindent In this section we prepare for our proof of Theorem~2 (in 
in Sections~6 and~7) by stating some of the more basic lemmas that 
are used in that proof. Before proceeding to these lemmas (and their proofs), 
we define one further piece of notation by putting:  
$${\rm e}(x)=\exp(2\pi ix)\qquad 
\hbox{for all $\,x\in{\Bbb R}$.}$$ 
This convenient notation will be used freely in this section, and in those that follow.   

\bigskip 

\goodbreak\proclaim{\Smallcaps Lemma~14 (Poisson summation over ${\Bbb Z}^2$ 
and over ${\frak O}={\Bbb Z}[i]$)}. 
Let $f : {\Bbb C}\rightarrow{\Bbb C}$. Suppose that the function
$F : {\Bbb R}^2\rightarrow{\Bbb C}$ given by $F(x,y)=f(x+iy)$ lies in 
the Schwartz space: so that, for all real $A\geq 0$ and all integers $j,k\geq 0$,
the function $|x+iy|^A {\partial^{j+k}\over\partial x^j\partial y^k} F(x,y)$
is continuous and bounded on ${\Bbb R}^2$. Then the Fourier transforms
$$\hat F(u,v)
=\int\limits_{-\infty}^{\infty}\int\limits_{-\infty}^{\infty} 
F(x,y)\,{\rm e}(-ux-vy) {\rm d}x {\rm d}y$$
and
$$\hat f(w)=\int_{\Bbb C} f(z)\,{\rm e}\!\left( -\Re (wz)\right) {\rm d}_{+}z
=\hat F\left( \Re (w) , -\Im (w)\right)$$
are complex-valued functions defined on ${\Bbb R}^2$ and ${\Bbb C}$, respectively.
Moreover, one has
$$\sum_{x=-\infty}^{\infty}\ \sum_{y=-\infty}^{\infty}F(x,y)
=\sum_{u=-\infty}^{\infty}\ \sum_{v=-\infty}^{\infty}\hat F (u,v),\eqno(5.1)$$
and, for all $\tau\in{\Bbb C}$ and all $\alpha ,\gamma\in{\frak O}$ with $\gamma\neq 0$ :
$$\sum_{\nu\in{\frak O}}f(\nu){\rm e}\left( \Re (\tau\nu)\right)
=\sum_{\xi\in{\frak O}}\hat f(\xi-\tau) ,\eqno(5.2)$$
$$\sum_{\nu\equiv\alpha\bmod \gamma{\frak O}}f(\nu)
={1\over |\gamma|^2}\sum_{\xi\in{\frak O}}\hat f\left({\xi\over\gamma}\right)
{\rm e}\left( \Re \left( {\alpha\xi\over\gamma}\right)\right) ,\eqno(5.3)$$
and
$$\sum_{\nu\in{\frak O}}f(\nu){\rm e}\left(\Re \left(\alpha {\nu^{*}\over\gamma}\right)\right)
={1\over |\gamma|^2}\sum_{\xi\in{\frak O}}\hat f\left({\xi\over\gamma}\right) S(\alpha , \xi ; \gamma),\eqno(5.4)$$
where $S(\alpha , \beta ; \gamma)$ denotes the Kloosterman sum over ${\Bbb Q}(i)$ 
that is given by:  
$$S(\alpha , \beta ; \gamma)
=\sum_{\delta\in ({\frak O}/\gamma{\frak O})^{*}}
{\rm e}\left( \Re \left({\alpha\delta^{*} +\beta\delta\over\gamma}\right)\right) .\eqno(5.5)$$

\medskip

\goodbreak
\noindent{\Smallcaps Proof.}\ 
For (5.1) see, for example [31, Chapter 13, Section 6].
By replacing in (5.1) the function $F$ by $G(x,y)=f(x+iy){\rm e}(\Re ((x+iy)\tau))$
 (and $\hat F$ by $\hat G$),
one obtains (5.2): this is justified, since one has
$G(x,y)=F(x,y){\rm e}(\Re (\tau) x) {\rm e}(-\Im (\tau) y)$,
from which it may be deduced that $G$ lies in the Schwartz space if $F$ does.
Noting that
$$\sum_{\nu\equiv\alpha\bmod\gamma{\frak O}}f(\nu)
={1\over |\gamma|^2}\sum_{\beta\in{\frak O}/\gamma{\frak O}}
{\rm e}\left(\Re \left( {\alpha\beta\over\gamma}\right)\right)
\sum_{\nu}f(\nu){\rm e}\left(\Re \left(-{\beta\over\gamma}\,\nu\right)\right) ,$$
we see (5.3) now follows by applying (5.2) to the inner sum on the right here.
Finally:
$$\sum_{\nu}f(\nu){\rm e}\left( \Re \left(\alpha\,{\nu^*\over\gamma}\right)\right)
=\sum_{\delta\in ({\frak O}/\gamma{\frak O})^*}
{\rm e}\left(\Re \left(\alpha\,{\delta^*\over\gamma}\right)\right)
\sum_{\nu\equiv\delta\bmod\gamma{\frak O}}f(\nu) ,$$
so that (5.4)-(5.5) is a corollary of (5.3)\ $\square$

\bigskip

\goodbreak\proclaim{\Smallcaps Lemma~15}. Suppose that $f$ and $F$ are as in Lemma~14. 
For $z=x+iy$ with $x,y\in{\Bbb R}$, let
$$(\Delta_{\Bbb C} f)(z)=(\Delta_{{\Bbb R}\times{\Bbb R}} F)(x,y)
=\left({\partial^2\over\partial x^2}+{\partial^2\over\partial y^2}\right) f(x+iy)\;.$$
Then $\Delta_{{\Bbb R}\times{\Bbb R}} F$ (the Laplacian of $F$) is a member of the Schwartz space.
The functions $f$ and $\Delta_{\Bbb C} f$ have 
Fourier transforms $\hat{f},\widehat{\Delta_{\Bbb C} f} : {\Bbb C}\rightarrow{\Bbb C}$
(defined as in Lemma~14) that are related to one another by:
$$|2\pi w|^2 \hat{f}(w) = -\widehat{\Delta_{\Bbb C} f}(w)\quad\hbox{for $w\in{\Bbb C}$.}\eqno(5.6)$$
For all $w\in{\Bbb C}-\{ 0\}$ and all $j\in{\Bbb N}\cup\{ 0\}$ , one  has
$$\big|\hat{f}(w)\bigr| 
=\left( 2\pi |w|\right)^{-2j}\Bigl|\widehat{\Delta_{\Bbb C}^j f}(w)\Bigr|
\leq\left( 2\pi |w|\right)^{-2j}\widehat{\,\bigl|\!\Delta_{\Bbb C}^j f\bigr|}(0) 
=\left( 2\pi |w|\right)^{-2j}\int_{\Bbb C}\bigl| \bigl( \Delta_{\Bbb C}^j f\bigr)(z)\bigr| 
{\rm d}_{+}z\;.\eqno(5.7)$$

\smallskip 

\goodbreak
\noindent{\Smallcaps Proof.}\ 
This lemma is a restatement of [45, Lemma~2.8]: 
see there for a proof\quad$\square$ 

\bigskip 

\goodbreak\proclaim{\Smallcaps Remarks~16}. Let $\alpha$ be a non-zero complex constant, 
and let $S_{\alpha}$ be the
associated `rotation-dilatation operator' which maps any function $f$ of 
a complex variable $z$
to the function $g(z)$ such that $g(z)=f(\alpha z)$ for all $z\in{\Bbb C}$.
Then, when $f$ is as in Lemma~14, the function $G : {\Bbb R}^2\rightarrow{\Bbb C}$
given by $G(x,y)=\left( S_{\alpha}f\right)(x+iy)$ is a member of the Schwartz space, and a linear change of
variables of integration shows that one has 
$$\widehat{S_{\alpha}f}(w)=|\alpha|^{-2}\hat{f}(w/\alpha)\quad\hbox{for $w\in{\Bbb C}\;$,}\eqno(5.8)$$
where the Fourier transforms, $\widehat{S_{\alpha}f}$ and $\hat{f}$, are as defined in Lemma~14.
An immediate, and useful, consequence of (5.8) is that if
$f(z)=f(|z|)$ for all $z\in{\Bbb C}$, then 
$\hat{f}(w)=\hat{f}(|w|)$ for all $w\in{\Bbb C}$.
\bigskip
 
\goodbreak\proclaim{\Smallcaps Lemma~17}. Let ${\frak N} : {\Bbb C}\rightarrow[0,\infty)$ satisfy
${\frak N}(z)=|z|^2$ for $z\in{\Bbb C}$.
Let $\Upsilon : [0,\infty)\rightarrow{\Bbb C}$ be an infinitely
differentiable function on $[0,\infty)$, and let the support of 
$\Upsilon$ be contained in the interval $\left[ e^{-1} , e\right]$. 
Suppose moreover that $\eta_0\in(0,1]$, and that 
$$\Upsilon^{(j)}(x)\ll_j \,\eta_0^{-j}\qquad\quad 
\hbox{for all $\,x\in [e^{-1},e]$, $\,j\in{\Bbb N}\cup\{ 0\}$.}\eqno(5.9)$$ 
Let $X>0$, and let $C : {\frak O}-\{ 0\}\rightarrow{\Bbb C}$ satisfy both
$$C(\xi)=0
\quad\hbox{for all $\xi\in{\frak O}$ such that $|\xi|^2/X^2\in\left(0,1/2\right]\cup[2,\infty)$}\eqno(5.10)$$
and
$$C(i\xi)=C(\xi)\quad\hbox{for all $\xi\in{\frak O}-\{ 0\}\;$.}\eqno(5.11)$$
Let $D>0$ and put
$$E_{D,C}
={1\over 2 D^2}\,\sum_{d=-\infty}^{\infty}\ \int\limits_{-\infty}^{\infty} \Upsilon\!\left({|2d+it|^2\over D^2}\right)
\left|\sum_{\xi\neq 0}C(\xi)\Lambda^{d}(\xi) |\xi|^{-2it}\right|^2 {\rm d}t\;.$$
Then, for all $Q\in(0,1/2]$ and all $j\in{\Bbb N}$, one has
$$E_{D,C}
=\sum_{\xi_1\neq 0}\sum_{\scriptstyle\xi_2\neq 0
\atop\scriptstyle\!\!\!\!\!\!\!\!\!\!\!\!\!\!\!\!\!\left|\xi_1 -\xi_2\right|<QX}
C\left(\xi_1\right)\overline{C\left(\xi_2\right)}\ 
\widehat{\Upsilon\circ {\frak N}}\!\left({D\over\pi}\,\log\left({\xi_1\over\xi_2}\right)\right) 
+O_j\left( \left(\eta_0 QD\right)^{-2j} \| C\|_1^2\right) ,\eqno(5.12)$$
where the Fourier transform $\widehat{\Upsilon\circ {\frak N}}$ is defined as in Lemma~14. 

\medskip

\goodbreak
\noindent{\Smallcaps Proof.}\ 
Let $\sigma$ denote the sum over $\xi$ in the definition of $E_{D,C}$.
By first expanding $|\sigma|^2$ as $\sigma \overline{\sigma}$, one can then integrate and sum (over $d$) 
term by term, so as to obtain:
$$2 D^2 E_{D,C}=\sum_{\xi_1\neq 0}\ \sum_{\xi_2\neq 0}C\left(\xi_1\right)\overline{C\left(\xi_2\right)}\,
F\!\left({1\over\pi i}\log\left({\xi_1\over\xi_2}\right)\right) ,\eqno(5.13)$$
where (given that $\Lambda^{d}(\alpha)$ is as in (1.2)) one has
$$F(\alpha)
=\sum_{d=-\infty}^{\infty}\ \int\limits_{-\infty}^{\infty}\Upsilon\left({|2d+it|^2\over D^2}\right)
{\rm e}\left(\Re \left( (2d+it)\overline{\alpha}\right)\right) {\rm d}t\;.\eqno(5.14)$$

Suppose now that $\xi_1$ and $\xi_2$ are non-zero elements of ${\frak O}$.
Then, since $F(\alpha +1/2)=F(\alpha)$, one has
$$F\left({1\over\pi i}\log\left({\xi_1\over\xi_2}\right)\right) 
=F\left({1\over\pi i}\log\left({\epsilon\xi_1\over\xi_2}\right)\right) ,\eqno(5.15)$$
where $\epsilon$ is the unique unit of ${\frak O}$ for which
$$-{\pi /4}<{\rm Arg}\left(\epsilon\xi_1/\xi_2\right)\leq \pi /4\;.\eqno(5.16)$$
Therefore, supposing now that
$$\alpha ={1\over\pi i}\log\left({\epsilon\xi_1\over\xi_2}\right) ,\eqno(5.17)$$
one has, by (5.16),
$$-{1\over 4}<\Re (\alpha)\leq{1\over 4}\;.\eqno(5.18)$$
By (5.14),
$$\eqalignno{
F(\alpha) &=\sum_{d=-\infty}^{\infty}\ \sum_{n=-\infty}^{\infty}\ \int\limits_{-1}^1\Upsilon\left({|2d+2ni+ui|^2\over D^2}\right)
{\rm e}\left(\Re \left( (2d+2ni+ui)\overline{\alpha}\,\right)\right) {\rm d}t = \cr
 &=\int\limits_{-1}^{1} G(\alpha , u) {\rm e}\left( \Im (\alpha) u\right) {\rm d}u\;,&(5.19)}$$
where
$$G(\alpha , u)
=\sum_{\beta\in{\frak O}}\Upsilon\left({|2\beta +iu|^2\over D^2}\right)
{\rm e}\left(\Re \left( 2 \overline{\alpha} \beta\right)\right)\;.$$
For $u\in{\Bbb R}$, an application of (5.2) of Lemma~14 yields:
$$\eqalignno{G(\alpha , u)
 &=\sum_{\gamma\in{\frak O}}\int_{\Bbb C}\Upsilon\left({|2z +iu|^2\over D^2}\right)
{\rm e}\left( -\Re \left(\left(\gamma -2 \overline{\alpha}\right)\!z\right)\right) 
{\rm d}_{+}z = \cr
 &=(D/2)^2 {\rm e}\left( -\Im (\alpha) u\right)
\sum_{\gamma\in{\frak O}}\widehat{\Upsilon\circ {\frak N}}\!\left(\left(\gamma -2\overline{\alpha}\,\right)\!D/2\right)
{\rm e}\left( -\Im (\gamma)u/2\right)  &(5.20)}$$
(where the second line follows by a change of variables of integration, 
as in Remarks~16). 
\par 
By the result (5.7) of Lemma~15, we have the upper bound 
$$\left|\widehat{\Upsilon\circ {\frak N}}(w)\right|  
\leq\left( 2\pi |w|\right)^{-2j} \int\limits_{-\infty}^{\infty} \int\limits_{-\infty}^{\infty} 
\left| \left( {\partial^2\over\partial x^2}+{\partial^2\over\partial y^2}\right)^{\!\!j} 
\Upsilon\left( |x+iy|^2\right)\bigg|_{(x,y)=(x_1,y_1)}\right| {\rm d}x_1\,{\rm d}y_1\;,$$   
for all $w\in{\Bbb C}-\{ 0\}$ and all $j\in{\Bbb N}$:  
note that, by reasoning similar to that which is used in the 
proof of the results (7.23)-(7.24) of 
Lemma~26 (below), it can be shown to follow from our hypotheses 
concerning $\Upsilon$ and $\eta_0\,$ (including, in particular, (5.9))  
that the function 
$(x,y)\mapsto\Upsilon(|x+iy|^2)$ is in the Schwartz space, 
and that for all $k,\ell\in{\Bbb N}\cup\{ 0\}$, and  
all points $(x_1,y_1)\in{\Bbb R}^2$, one has  
$${\partial^{k+\ell}\over\partial x^k \partial y^{\ell}} 
\,\Upsilon\left( |x+iy|^2\right) \bigg|_{(x,y)=(x_1,y_1)} 
=\cases{ O_{k,\ell}\left( \eta_0^{-(k+\ell)}\right) &if $\,e^{-1}\leq |x_1+iy_1|^2\leq e$, \cr 
0 &otherwise.}$$ 
Consequently, for all $w\in{\Bbb C}-\{ 0\}$ and all $j\in{\Bbb N}$,  we have:   
$$\left|\widehat{\Upsilon\circ {\frak N}}(w)\right|  
\leq \left( 2\pi |w|\right)^{-2j} \int\limits_{-\sqrt{e}}^{\sqrt{e}} 
\ \,\int\limits_{-\sqrt{e}}^{\sqrt{e}}  
O_j\left(\eta_0^{-2j}\right) {\rm d}x_1\,{\rm d}y_1 
\ll_j \,\left(\eta_0 |w|\right)^{-2j}\;.\eqno(5.21)$$

The case $j=2$ of (5.21) 
ensures uniform absolute convergence (for all $u\in{\Bbb R}$) of the sum 
over $\gamma\in{\frak O}$ in (5.20), 
and so justifies both the substitution of the expression 
for $G(\alpha ,u)$ obtained in (5.20) into (5.19), and 
the term by term integration, with respect to $u\in[-1,1]$, of the  
series expansion of $G(\alpha,u) {\rm e}(\Im(\alpha) u)$  brought about 
by that substitution. This term by term integration  
has the effect of annihilating any term for which the index $\gamma$ is not real 
(this being so by virtue of the fact that 
the condition of summation $\gamma\in{\frak O}$ implies 
$\Im (\gamma)\in{\Bbb Z}$). 
Moreover, the terms with $\gamma\in{\frak O}\cap{\Bbb R}={\Bbb Z}$ are 
trivial to integrate (being independent of $u$), and so, by (5.19) and (5.20), one 
obtains:  
$$F(\alpha)={D^2\over 2}\sum_{g=-\infty}^{\infty}
\widehat{\Upsilon\circ {\frak N}}\!\left(\left( g/2 -\overline{\alpha}\,\right)\!D\right)\;.$$
Given (5.18), it follows by (5.21) that if $g\in{\Bbb Z}-\{ 0\}$, then 
$\widehat{\Upsilon\circ {\frak N}}\!\left(\left( g/2 -\overline{\alpha}\,\right)\!D\right) 
\ll_j (\eta_0 D|g|)^{-2j}$
for $j\in{\Bbb N}$. By applying this observation to the sum over $g\in{\Bbb Z}$ 
occurring in the expression just obtained for $F(\alpha)$, 
we find that 
$$F(\alpha) 
={D^2\over 2}\,\widehat{\Upsilon\circ {\frak N}}\!\left( -\overline{\alpha}\,D\right)
+O_j\!\left( \eta_0^{-2j} D^{2-2j}\right)\qquad\ \hbox{for all $\,j\in{\Bbb N}\;$.} 
\eqno(5.22)$$
Moreover, by (5.22) and (5.21) (again), it follows that if $\alpha\neq 0$ then
$$F(\alpha)\ll_j\,\left( |\alpha|^{-2j}+1\right)\eta_0^{-2j} D^{2-2j}\qquad\hbox{for all $\,j\in{\Bbb N}\;$.} 
\eqno(5.23)$$

The results (5.22) and (5.23) will be used to estimate the terms of the double sum in (5.13). 
Given (5.10), and (5.13), it is trivial that one needs these estimates only in the cases where
$$X/\sqrt{2}\leq \left|\xi_1\right| ,\left|\xi_2\right|\leq\sqrt{2} X\;.\eqno(5.24)$$
Choose now any $Q$ such that $0<Q\leq 1/2$. Taking, in turn, each pair
$\xi_1,\xi_2\in{\frak O}$ satisfying (5.24), suppose that
(5.22) is applied if and only if 
$$\left|\epsilon\xi_1 -\xi_2\right|<Q X\;,\eqno(5.25)$$
while (5.23) is applied if and only if the condition (5.25) is not satisfied.
Then, in the cases where the bound (5.23) is applied,  
one will have $|\alpha|>Q/6$: indeed, 
by (5.17), the condition  $Q\leq 1/2$ and the inequalitities  $|\exp(z) -1|\leq |z|\exp(|z|)$, 
$\pi \exp(\pi /12)<3\sqrt{2}$ and (5.24), it follows that 
if one were to have $|\alpha|\leq Q/6$ then 
the condition (5.25) would be satisfied, and so (5.23) would not be being applied.
Consequently, where it is applied, (5.23) shows:
$$F(\alpha)\ll_j\,\left( (Q/6)^{-2j} +1\right) \eta_0^{-2j} D^{2-2j}
\ll_j \left(\eta_0 Q\right)^{-2j} D^{2-2j}\;.$$
In the cases where it is instead (5.22) that is implied, 
one has (5.24) and (5.25), by which
(given that $0<Q\leq 1/2$) it follows that
$$\left|\epsilon\xi_1 -\xi_2\right| <\left|\xi_2\right| /\sqrt{2}\;.$$

Since the last inequality renders redundant the condition (5.16) that was
imposed on the choice of unit $\epsilon\in\{ 1 , i , -1 , -i\}$, it 
therefore follows from the conclusions of the above paragraph that, by means of
(5.13), (5.15), (5.17) and the estimates (5.22) and (5.23), one can deduce that 
$$E_{D,C}
={1\over 4}\,\sum_{\epsilon^4 =1}\ 
\sum_{\xi_1\neq 0}\sum_{\scriptstyle\xi_2\neq 0
\atop\scriptstyle\!\!\!\!\!\!\!\!\!\!\!\!\!\!\!\!\!\left|\epsilon\xi_1 -\xi_2\right|<QX}
C\left(\xi_1\right)\overline{C\left(\xi_2\right)}\ 
\widehat{\Upsilon\circ {\frak N}}\!\left( -\overline{{1\over\pi i}\,\log\left({\epsilon\xi_1\over\xi_2}\right)}\,D\right) 
+O_j\left( \| C\|_1^2\left(\eta_0 Q D\right)^{-2j}\right) ,$$
for all $j\in{\Bbb N}$. The  conclusion of the lemma now follows,
via the substitution $\xi_1=\epsilon^{-1} \xi$, the observation 
that the hypothesis (5.11) implies that 
$\sum_{\epsilon^4 =1} C\left(\epsilon^{-1} \xi\right) =4 C(\xi)$,  
for $\xi\in{\frak O}-\{ 0\}$, and the final point noted in Remarks~16\quad$\square$

\bigskip

\goodbreak\proclaim{\Smallcaps Remarks~18}. Let $0<\eta\leq 1$, and let  
$\Omega_{\eta}$ be an infinitely differentiable function 
from $[0,\infty)$ to $[0,\infty)$ that 
has as its support some subset of the interval $[1,e^{\eta}]$, and that satisfies both
$$\Omega_{\eta}^{(j)}(x)\ll_j \eta^{-j}\quad\hbox{(for $x>0$ and $j=0,1,2,\ldots\;$ )}
\qquad\hbox{and}\qquad{1\over\eta}\,\int\limits_{0}^{\infty} \Omega_{\eta}(x) {\rm d}x =1$$ 
(such functions do exist: one example is the mapping   
$x\mapsto\Phi(2 (e^{\eta}-1)^{-1} (x-1) - 1)$, where $\Phi : {\Bbb R}\rightarrow 
[0,e^{-1}]$ is given by the equation (3.100) of Remarks~13, above).  
Suppose moreover that $\delta$ is a real number satisfying $0<\delta\leq (4e)^{-1}$. 
Then, by multiplying both sides of the equation (5.12) by  
$(\delta\eta)^{-1}\Omega_{\eta}\bigl( Q^2 /\delta\bigr)$, applying 
the substitution $Q=\sqrt{\delta_1}$, and then  
integrating both sides of the resulting equation with respect to 
the (positive valued) variable $\delta_1$, one obtains: 
$$E_{D,C} = O_j\left(\delta^{-j} \left(\eta_0 D\right)^{-2j} \| C\|_1^2\right) 
+\sum_{\xi_1\neq 0}\ \sum_{\xi_2\neq 0}
W_{\eta}\!\!\left( {\left|\xi_1 -\xi_2\right|^2\over\delta X^2}\right)
C\!\left(\xi_1\right)\overline{C\!\left(\xi_2\right)}\ 
\widehat{\Upsilon\circ {\frak N}}\!\left({D\over\pi} 
\,\log\left({\xi_1\over\xi_2}\right)\right) ,\eqno(5.26)$$
where, for $u\geq 0$,
$$W_{\eta}(u)={1\over\eta}\,\int\limits_{u}^{\infty}\Omega_{\eta}(x) {\rm d}x\;.\eqno(5.27)$$
Moreover, given the stated properties of the function $\Omega_{\eta}(x)$, the function
$W_{\eta}(u)$ here is real valued, monotonic decreasing and infinitely differentiable on
$[0,\infty)$, and satisfies: $W_{\eta}(u)=1$ if $0\leq u\leq 1\;$;
$W_{\eta}(u)=0$ if $u\geq e^{\eta}$; 
$W_{\eta}^{(j)}(u)\ll_j \eta^{-j}$ for $u>0$ and $j\in{\Bbb N}\cup\{ 0\}$.

\bigskip

\goodbreak\proclaim{\Smallcaps Lemma~19 (estimates for Kloosterman sums)}. 
Let $\alpha,\beta,\gamma\in{\frak O}$, with $\gamma\neq 0$, and
let the Kloosterman sum $S=S(\alpha, \beta ; \gamma)$ be given by the equation (5.5) of 
Lemma~14. Then 
$$|S|\leq \phi_{\frak O}(\gamma) 
=|\gamma|^2\prod_{(\pi_1)\mid (\gamma)}\left( 1-{1\over |\pi_1|^2}\right)
\leq |\gamma|^2 ,\eqno(5.28)$$
where one has $\phi_{\frak O}(\gamma)=\bigl| ({\frak O}/\gamma{\frak O})^{*}\bigr|\,$
(Euler's function for Gaussian integers), and where the product is over distinct 
prime ideal factors $(\pi_1)=\pi_1{\frak O}$ 
of the ideal $(\gamma)=\gamma{\frak O}\subset{\frak O}$ (and so is 
equal to $1$ if and only if $\gamma$ is a unit of ${\frak O}$). One has also 
the (much deeper) `Weil-Estermann estimate': 
$$|S|^2\leq 2^3 \tau_2^2(\gamma)\,|(\alpha , \beta , \gamma) \gamma|^2\;,\eqno(5.29)$$ 
where $\tau_2(\gamma)=\sum_{\delta\mid\gamma} 1$ 
and $(\alpha,\beta,\gamma)$ is an arbitrary highest common factor of 
$\alpha$, $\beta$ and $\gamma$. 
Moreover, if $\gamma\mid\beta$ then one has the (elementary) `Ramanujan sum evaluation': 
$$S={1\over 4}\sum_{\nu\mid (\alpha,\gamma)}\mu_{\frak O}\!\left({\gamma\over\nu}\right) |\nu|^2 
=\mu_{\frak O}\!\!\left({\gamma\over (\alpha,\gamma)}\right)
|(\alpha,\gamma)|^2 
\!\!\!\!\prod_{\scriptstyle (\pi_1)\mid((\alpha,\gamma))\atop\scriptstyle 
(\pi_1)\,{\mid\!\!\!\not}\ \,(\gamma/(\gamma,\alpha))}
\!\!\left( 1-{1\over|\pi_1|^2}\right) 
={\mu_{\frak O}\!\left( \gamma /(\alpha,\gamma)\right) 
\phi_{\frak O}(\gamma)\over\phi_{\frak O}\!\left(\gamma /(\gamma,\alpha)\right)}\;,\eqno(5.30)$$
where, as in (5.28), the 
product is over distinct prime ideals $(\pi_1)\subset{\frak O}$, 
while $$\mu_{\frak O}(\kappa)=\cases{0 &if there exists a Gaussian prime $\pi_1$ such that 
$\pi_1^2\mid\kappa$, \cr 
(-1)^{\omega(\kappa)} &otherwise,}$$ 
with $\omega(\kappa)$ denoting the number of prime ideals in the 
ring ${\frak O}$ that contain $\kappa$. 
The result (5.30) implies that 
$$|S|\leq |(\alpha,\gamma)|^2\quad\ \hbox{if $\ \gamma\mid\beta\;$.}\eqno(5.31)$$

\medskip\par

\goodbreak
\noindent{\Smallcaps Proof.}\ 
The first two equalities of (5.30) follow from [44, Equations~(5.35)] 
and points noted in the last few lines of the proof of [44, Lemma~5.7]. 
The definition (5.5) implies that 
$|S|\leq S(0,0;\gamma)=\phi_{\frak O}(\gamma)$, 
and so we have both the first inequality in (5.28) and, by virtue of  
the case $\alpha =\beta =0$ of the first two equalities of (5.30), 
the equality in (5.28) also. As for the final inequality in (5.28),  
that is a trivial consequence of the fact that each factor $1-|\pi_1|^{-2}$ 
occurring in the product in (5.28) must certainly  
satisfy $0\leq 1-|\pi_1|^{-2}\leq 1$. 
Similarly, the result (5.31) 
follows from the first two equalities in (5.30).  
Two applications of the equality in (5.28) are enough to give the 
final equality in (5.30). 
\par 
The bound (5.29) is a corollary of the more general results 
[4, Proposition~9 and Theorem~10] of Bruggeman and Miatello\quad$\square$

\bigskip

\goodbreak 
\noindent{\SectionHeadingFont 6. Beginning the proof of Theorem~2}

\medskip 

\noindent In this section we show, in effect, that 
if ${\cal E}$ is the term given by Equation~(1.30) (in combination  
with the definitions in (1.21)-(1.29)) then 
$|{\cal E}|$ is either quite small (smaller, subject to the condition (1.18), 
than the bounds (1.31) and (1.32)  would imply), or else 
is bounded above by a certain sum of Kloosterman sums (a multiple of the sum 
$Y$ defined in Lemma~22).  
Note that the mean value on the left-hand side of Equation~(1.30) 
has an obvious `rational integer analogue', which was
studied by Deshouillers and Iwaniec: their work in [9, Section~2] and  
[7, Section~9.2] provided the model for the principal steps of this section. 
\par 
We assume throughout this section (and the next) that the numbers 
$\varepsilon$, $\eta$, $K_0$, $K_1$, $K_2$, $M_1$
are as in Theorem 2, 
that  the function $C : {\frak O}-\{ 0\}\rightarrow{\Bbb C}$ 
and terms $c(d,it)\,$ ($d\in{\Bbb Z}$, $t\in{\Bbb R}$) 
are as defined in (1.21) and (1.29),  
that $T$ satisfies (1.18), and that $N$ is given by the equation~(1.24).
The functions $w_0(x)$, $w_1(x)$, $w_2(x)$ and $a(\mu)$ are also
assumed to be as in Theorem 2, except that (to make conditions of the form $x\neq 0$, $\mu\neq 0$ implicit) 
their respective domains are enlarged to
include $0$, by defining $a(0)=0$ and $w_i(0)=0$ for $i=1,2,3$.
When $f : {\Bbb C}\rightarrow{\Bbb C}$, the Fourier transform $\hat f$ 
is as defined in Lemma~14; the Euclidean Laplace operator $\Delta_{\Bbb C}$ is 
defined in Lemma~15. As in Lemma~17, we take ${\frak N}(z)$ to denote the 
function on ${\Bbb C}$ satisfying ${\frak N}(z)=|z|^2$ for all $z\in{\Bbb C}$. 
We shall moreover assume that the functions 
$\Omega_{\eta}$ and $W_{\eta}$ are as described in Remarks~18 
(were it necessary, we could explicitly 
define just such a pair of functions).  
\par 
In cases where (in addition to (1.18) and (1.24)) one has $T=O(1)$, 
the results (1.30)-(1.32) of Theorem~2 can be verified by means 
of a direct (and quite trivial) upper bound estimate for 
the absolute value of the sum on the left-hand side of 
Equation~(1.30), in combination with similarly direct and trivial upper bound 
estimates for 
$|{\cal D}_0|$, $|{\cal D}^{*}_1|$ and $|{\cal D}^{*}_2|$ 
(where ${\cal D}_0$, ${\cal D}^{*}_1$ and ${\cal D}^{*}_2$ are 
as defined in (1.21)-(1.28)). Therefore, in completing our proof of 
Theorem~2, we may certainly suppose that 
$$T\geq 2^{10}\pi^8\;.\eqno(6.1)$$ 
See (7.52), below, for what motivates this particular choice of a lower bound for $T$.

\bigskip 

\goodbreak\proclaim{\Smallcaps Lemma~20}. Let 
$$\Xi=T^{\varepsilon} M_1 K_2 /K_1\eqno(6.2)$$ 
and, for $\psi,\psi',\nu,\xi\in{\frak O}$, put 
$$L\!\left(\psi,\psi';\nu,\xi\right) 
=\left|\psi\psi'\right|^2\!\!\!\int_{\Bbb C}\!\int_{\Bbb C} 
\!w_0\!\!\left({\left|\psi\psi' z_1 z_2\right|^2\over T}\right)
w_1\!\!\left(\left|\psi' z_2\right|^2\right)
\overline{w_1\!\!\left(\left|\psi z_2\right|^2\right)}\!\left| z_2\right|^2
{\rm e}\!\left(\Re \left(\nu z_1 +\xi z_2\right)\right)\!{\rm d}_{+}z_1\,{\rm d}_{+}z_2\;.
\eqno(6.3)$$
Then 
$$\eqalign{\sum_{d=-\infty}^{\infty}\ \int\limits_{-\infty}^{\infty}
\left| c(d,it)\right|^2
 &w_0\!\left({|2d+it|^2\over\pi^2 T}\right) {\rm d}t = \cr 
&={2\pi}\,\widehat{w_0\!\circ{\frak N}}(0)\,T \| C\|_2^2 
+ (\pi /2)\left( {\cal D}^{\star}+{\cal E}^{\star}\right) 
+O_{\varepsilon,\eta}\!\left(  T^{5\varepsilon +1/2} K_1 K_2 M_1 \| a\|_2^2\right)\;,}\eqno(6.4)$$ 
where
$${\cal D}^{\star}=\sum_{\varphi\neq 0}
\sum\sum\!\!\!\!\!\!\!\!
\sum_{\!\!\!\!\!\!\!\!\!\!\!{\scriptstyle\kappa_2\ \ \mu_1\ \ \kappa_4\ \ \mu_2
\atop\scriptstyle (\kappa_2\mu_1 , \kappa_4\mu_2) \sim\varphi}}
\!\!\!\!\!\!\!\!\sum w_2\!\left(\left|\kappa_2\right|^2\right) a\left(\mu_1\right)
\overline{w_2\!\left(\left|\kappa_4\right|^2\right) a\left(\mu_2\right)}
\ G_{\varphi}\!\left({\kappa_2\mu_1\over\varphi} , {\kappa_4\mu_2\over\varphi} ; 0\right) , 
\eqno(6.5)$$
and
$${\cal E}^{\star} 
 =\sum_{\varphi\neq 0}
\sum\sum\!\!\!\!\!\!\!\!
\sum_{\!\!\!\!\!\!\!\!\!\!\!{\scriptstyle\kappa_2\ \ \mu_1\ \ \kappa_4\ \ \mu_2
\atop\scriptstyle (\kappa_2\mu_1 , \kappa_4\mu_2) \sim\varphi}}
\!\!\!\!\!\!\!\!\sum w_2\!\left(\left|\kappa_2\right|^2\right) a\left(\mu_1\right)
\overline{w_2\!\left(\left|\kappa_4\right|^2\right) a\left(\mu_2\right)}
\ H_{\varphi}\!\left({\kappa_2\mu_1\over\varphi} , {\kappa_4\mu_2\over\varphi}\right) , 
\eqno(6.6)$$
with 
$$G_{\varphi}\left(\psi_1 , \psi_2 ; \xi\right)
 =\sum_{\nu\neq 0}W_{\eta}\!\left({|\varphi\nu|^2\over N}\right)
L\!\left(\psi_1,\psi_2;\nu,\xi\right) 
{\rm e}\left( \Re \left({\psi_1^{*}\nu\xi\over\psi_2}\right)\right)\eqno(6.7)$$
and 
$$H_{\varphi}\left(\psi_1 , \psi_2\right)
=\sum_{\xi\neq 0} W_{\eta}\biggl({\left|\varphi\xi\right|^2\over\Xi}\biggr)
G_{\varphi}\left(\psi_1 , \psi_2 ; \xi\right) . \eqno(6.8)$$

\medskip

\goodbreak
\noindent{\Smallcaps Proof.}\  By (1.2), (1.29) and (1.21), one has
$$c(d , it)=\sum_{\xi\neq 0} C(\xi) \Lambda^{d}(\xi) |\xi|^{-2it}\;.\eqno(6.9)$$
Since $|i\kappa|^2=|\kappa|^2$ for $\kappa\in{\frak O}$, and since $i{\frak O}={\frak O}$,
it is evident from (1.21) that $C(\xi)$ satisfies the condition (5.11) of Lemma~17.
Moreover, by (1.19) and (1.20) (in which $0<\eta\leq (\log 2)/3$\/), this $C(\xi)$ also
satisfies, for
$X=( K_1 K_2 M_1 )^{1/2}$, the condition (5.10) of Lemma~17.
Therefore Lemma~17 may be applied with $\Upsilon(x)=w_0(x)$,  
and with $C(\xi)$ as in (1.21). In particular,  
by (6.9) and the case $\eta_0=\eta$, $X=( K_1 K_2 M_1 )^{1/2}$, 
$D=\pi T^{1/2}$, $\delta =T^{\varepsilon -1}$
of (5.26) (the corollary of Lemma~17 noted in Remarks~18), and by 
(1.18)-(1.21), (1.24) and the Cauchy-Schwarz inequality,  
one has  
$$\eqalign{\sum_{d=-\infty}^{\infty}\ &\int\limits_{-\infty}^{\infty}
\left| c(d,it)\right|^2
w_0\!\left({|2d+it|^2\over\pi^2 T}\right) {\rm d}t = \cr 
 &=O_{j,\eta}\left( T^{7/2-j\varepsilon}\| a\|_2^2\right) 
+2\pi T\sum_{\xi_1}\sum_{\xi_2}
W_{\eta}\!\!\left( {\left|\xi_1 -\xi_2\right|^2\over N}\right)
C\!\left(\xi_1\right)\overline{C\!\left(\xi_2\right)}\ 
\widehat{w_{0}\circ {\frak N}} 
\!\left(\!\sqrt{T}\,\log\!\left({\xi_1\over\xi_2}\right)\right) ,} \eqno(6.10)$$
for all $j\in{\Bbb N}$. 

Taking $j=[(A+4)/\varepsilon]+1$ (for an arbitrary $A\geq 0$\/) 
one may replace
the last term in (6.10) by 
$O_{\varepsilon,\eta,A}(T^{-A}\| a\|_2^2)$.
The terms of the first sum on the right in (6.10) are `diagonal' when $\xi_1=\xi_2$,
and otherwise are `off-diagonal'. Let $C\bigl(\xi_1\bigr)$ and $C\bigl(\xi_2\bigr)$
in (6.10) be expressed, through (1.21), in terms of variables of summation
$\kappa_1,\kappa_2,\mu_1$ and $\kappa_3,\kappa_4,\mu_2$ satisfying
$\kappa_1\kappa_2\mu_1=\xi_1$ and $\kappa_3\kappa_4\mu_2=\xi_2$.
Then, by separating diagonal from off-diagonal terms, and, for each pair $\varphi,\nu\in{\frak O}-\{ 0\}$, 
grouping together off-diagonal terms with both  
$$\left( \kappa_2\mu_1 , \kappa_4\mu_2\right)\sim\varphi$$
and 
$$\kappa_1\kappa_2\mu_1-\kappa_3\kappa_4\mu_2=\xi_1 -\xi_2 =\varphi\nu\;,$$
one may rewrite the case $j=[(A+4)/\varepsilon]+1$ of (6.10) as 
$$\sum_{d=-\infty}^{\infty}\ \int\limits_{-\infty}^{\infty}
\left| c(d,it)\right|^2
w_0\!\left({|2d+it|^2\over\pi^2 T}\right) {\rm d}t
={2\pi}\,\widehat{w_0\!\circ{\frak N}}(0)\,T \| C\|_2^2 
+ (\pi /2)\,{\cal D}' + O_{\varepsilon,\eta,A}\left(T^{-A}\| a\|_2^2\right)\;, 
\eqno(6.11)$$
where
$${\cal D}'
=\sum_{\varphi\neq 0}\sum_{\nu\neq 0}W_{\eta}\!\left({|\varphi\nu|^2\over N}\right)
\sum\sum\!\!\!\!\!\!\!\!
\sum_{\!\!\!\!\!\!\!\!\!\!\!{\scriptstyle\kappa_2\ \ \mu_1\ \ \kappa_4\ \ \mu_2
\atop\scriptstyle (\kappa_2\mu_1 , \kappa_4\mu_2) \sim\varphi}}
\!\!\!\!\!\!\!\!\sum w_2\!\left(\left|\kappa_2\right|^2\right) a\left(\mu_1\right)
\overline{w_2\!\left(\left|\kappa_4\right|^2\right) a\left(\mu_2\right)}
\,U_{\nu}\!\left({\kappa_2\mu_1\over\varphi} , {\kappa_4\mu_2\over\varphi}\right) ,
\eqno(6.12)$$
with
$$U_{\nu}\!\left(\psi_1 , \psi_2\right)
=\ T\ \sum\!\!\!\!\!\!\!\!\sum_{\!\!\!\!\!\!\!\!\!\!{\scriptstyle\kappa_1\quad\kappa_3\atop\scriptstyle\psi_1\kappa_1-\psi_2\kappa_3=\nu}}
\widehat{w_0\circ{\frak N}}\!\left(\!\sqrt{T}\,\log\!\left({\psi_1\kappa_1\over\psi_2\kappa_3}\!\right)\right) 
w_1\!\left(\left|\kappa_1\right|^2\right)
\overline{w_1\!\left(\left|\kappa_3\right|^2\right)  }
\eqno(6.13)$$
(note that we have here simplified the first term on the right-hand side 
Equation~(6.11) by making use of the equality $W_{\eta}(0)=1$, which is 
one of the properties of the function $W_{\eta}$ that are mentioned below (5.27)).  
Within the sum in (6.13) one has 
$$\kappa_3= {\psi_1\kappa_1-\nu\over\psi_2}=\left( 1-{\nu\over\psi_1\kappa_1}\right) {\psi_1\kappa_1\over\psi_2}\;,$$
where, given how (1.19), (1.20) and the properties of $W_{\eta}(u)$ noted below (5.27) 
effectively  restrict the ranges of each of $|\varphi\nu|$, $|\kappa_2\mu_1|$,  
$|\kappa_4\mu_2|$ and $|\kappa_1|$ in the sums occurring in (6.12) and (6.13), 
it may be supposed that 
$$\psi_2\neq 0\qquad{\rm and}\qquad\left|{\nu\over\psi_1\kappa_1}\right|^2\leq{e^{4\eta} N\over K_1 K_2 M_1} 
<3T^{\varepsilon-1}<{1\over 5}\;.$$
By this, the Taylor expansions of $\log z$ and $\exp(w)$ about $z=1$, $w=0$, 
the definition of the Fourier transform, and (1.19), for $i\in\{ 0,1\}$, 
one finds that, for the relevant (constrained) values of 
$\varphi,\nu,\psi_1,\psi_2,\kappa_1\in{\frak O}$, 
$$\widehat{w_0\circ{\frak N}}\!\left(\!\sqrt{T}\,\log\!\left({\psi_1\kappa_1\over\psi_2\kappa_3}\!\right)\right) 
\overline{w_1\!\left(\left|\kappa_3\right|^2\right)}
= \widehat{w_0\circ{\frak N}}\!\left(\sqrt{T}\,{\nu\over\psi_1\kappa_1}\right) 
\overline{w_1\!\left(\left|{\psi_1\over\psi_2}\,\kappa_1\right|^2\right)}
+O\!\left( \eta^{-1} T^{(\varepsilon -1)/2}\right) .$$ 
Therefore, by (6.13) and the case $i=1$ of (1.19), one has in (6.12):
$$U_{\nu}\!\left(\psi_1 , \psi_2\right)
=\widetilde{U}_{\nu}\!\left(\psi_1 , \psi_2\right)
+V_{\nu}\!\left(\psi_1 , \psi_2\right)\;,\eqno(6.14)$$
where
$$V_{\nu}\!\left(\psi_1 , \psi_2\right)
\ll \eta^{-1} T^{(\varepsilon +1)/2}\sum_{\scriptstyle \kappa_1\equiv\psi_1^{*}\nu\bmod \psi_2 {\frak O} 
\atop\scriptstyle e^{-\eta}\leq |\kappa_1|^2 /K_1\leq e^{\eta}} 1  
\eqno(6.15)$$
and 
$$\widetilde{U}_{\nu}\!\left(\psi_1 , \psi_2\right)
=\  T\sum_{\kappa_1\equiv\psi_1^{*}\nu\bmod \psi_2 {\frak O}}
\widehat{w_0\circ{\frak N}}\!\left({\sqrt{T}\,\nu\over\psi_1\kappa_1}\right)w_1\!\left(\left|\kappa_1\right|^2\right)
\,\overline{w_1\!\left(\left|{\psi_1\over\psi_2}\,\kappa_1\right|^2\right)} ,\eqno(6.16)$$
with $\psi_1^{*}\in{\frak O}$ being an arbitrary solution of the congruence
$\psi_1\psi_1^{*}\equiv 1\bmod \psi_2 {\frak O}$. 

Using the definition of $\widehat{w_0\circ{\frak N}}(s)$ as an integral, it follows
by a change of variable of integration, and a change in  the order of summation and
integration, that (6.16) may be rewritten as:
$$\widetilde{U}_{\nu}\!\left(\psi_1 , \psi_2\right)
=\hat{F}_{\psi_1,\psi_2,\nu}(\nu)\;,\eqno(6.17)$$
where
$$F_{\psi_1,\psi_2,\nu}(z)=\sum_{\kappa\equiv\psi_1^{*}\nu\bmod \psi_2 {\frak O}}
f_{\psi_1,\psi_2,z}(\kappa)\;,\eqno(6.18)$$
with
$$f_{\psi_1,\psi_2,z}(\kappa)
= \left|\psi_1\kappa\right|^2 w_1\!\left(\left|\kappa\right|^2\right)
\overline{w_1\!\!\left(\left|{\psi_1\over\psi_2}\,\kappa\right|^2\right)}
\ w_0\!\!\left({\left|z\psi_1\kappa\right|^2\over T}\right)\;.\eqno(6.19)
$$

By (6.19), (1.19) and (for example) [44, Lemma~9.4], 
the functions $w_0(x)$ and $w_1(x)$ (extended to have $w_0(0)=w_1(0)=0$) are
such that the function $G : {\Bbb R}^2\rightarrow{\Bbb C}$ given by
$G(x,y)=f_{\psi_1,\psi_2,z}(x+iy)$ lies in the Schwartz space. Therefore
it follows from (6.18) and the result (5.3) of Lemma~14 that
$$F_{\psi_1,\psi_2,\nu}(z)
={1\over\left|\psi_2\right|^2}\sum_{\xi\in{\frak O}}
\hat{f}_{\psi_1,\psi_2,z}\left({\xi\over\psi_2}\right)
{\rm e}\left(\Re \left({\psi_1^{*}\nu\xi\over\psi_2}\right)\right) \qquad\quad\  
\hbox{($z\in{\Bbb C}$).}\eqno(6.20)$$
Note that when $\psi_1,\psi_2,z$ are known (and so to be treated as constants)
it then follows by (6.19) and (1.19) that one  has
$f_{\psi_1,\psi_2,z}(\kappa)=h\bigl(|\kappa|^2\bigr)$, where the function
$h : [0,\infty)\rightarrow{\Bbb C}$ is infinitely differentiable, 
satisfies
$$h^{(k)}(x)=\cases{O_k\left( K_1 \left|\psi_1\right|^2 (\eta x)^{-k}\right) &if $e^{-\eta}K_1\leq x\leq e^{\eta}K_1\;$, \cr
0 &otherwise,}\eqno(6.21)$$
for all $x>0$ and all $k\in{\Bbb N}\cup\{ 0\}$, and is, moreover, identically zero
unless it is the case that 
$$e^{-2\eta} T/K_1\leq \left| z\psi_1\right|^2\leq e^{2\eta} T/K_1 .\eqno(6.22)$$

As an operator on smooth functions of a complex variable $s=u+iv$ the Laplacian $\Delta_{\Bbb C}$ satisfies:
$$\Delta_{\Bbb C}={\partial^2\over\partial u^2}+{\partial^2\over\partial v^2}= 4{\partial\over\partial s}{\partial\over\partial\overline{s}}$$
where
$${\partial\over\partial s}={1\over 2}\left({\partial\over\partial u}-i{\partial\over\partial v}\right)\quad
\hbox{and}\quad
{\partial\over\partial\overline{s}}={1\over 2}\left({\partial\over\partial u}+i{\partial\over\partial v}\right) .$$
Therefore, and since (as was found above) $f_{\psi_1,\psi_2,z}(s)=h\bigl( |s|^2\bigr)=h\bigl(s\overline{s}\bigr)$,
one may determine that
$$\left(\Delta_{\Bbb C}^j f_{\psi_1,\psi_2,z}\right)(s)
=\sum_{r=0}^j {4^j (j!)^2  |s|^{2j-2r}h^{(2j-r)}\!\left( |s|^2\right)\over (r!) ((j-r)!)^2}  
\qquad\ \hbox{for $\,s\in{\Bbb C}$, $j\in{\Bbb N}\cup\{ 0\}$.}$$
By this, (6.21) and the bound (5.7) of Lemma~15, it follows that if $\xi\neq 0$, 
then, for all $j\in{\Bbb N}\cup\{ 0\}$, 
$$\eqalign{
\hat{f}_{\psi_1,\psi_2,z}\left({\xi\over\psi_2}\right) 
\ll \left|{\xi\over\psi_2}\right|^{-2j}\int_{\Bbb C}\left|\left(\Delta_{\Bbb C}^j f_{\psi_1,\psi_2,z}\right)(s)\right|
{\rm d}_{+}s 
 &=\left|\psi_2\right|^{2j}\left|\xi\right|^{-2j} O_j\!\left( K_1^2\left|\psi_1\right|^2
\left(\eta K_1^{1/2}\right)^{-2j}\right) \ll_{\eta , j} \cr
&\ll_{\eta , j} K_1^{2-j}\left|\psi_1\right|^2\left|\psi_2\right|^{2j} |\xi|^{-2j}\;.
}$$
Given (1.19), and given how $\psi_1$ and $\psi_2$ arise in (6.12)-(6.13), it is 
here sufficient to treat cases in which one has $\left|\psi_i\right|^2\asymp |\varphi|^{-2}K_2 M_1$ 
for $i=1,2$. In such a case the above bound 
for $\hat{f}_{\psi_1,\psi_2,z}(\xi /\psi_2)$ 
implies that if $X_1\asymp |\varphi|^{-2}\Xi$, where
$\Xi$ is given by (6.2), then, for $j\geq 2$, one has 
$${1\over\left|\psi_2\right|^2}\sum_{|\xi|^2\geq X_1}\left|\hat{f}_{\psi_1,\psi_2,z}\left({\xi\over\psi_2}\right)\right|
\ll_{\eta , j} K_1^{2-j}\left({K_2 M_1\over |\varphi|^2}\right)^{\!j} X_1^{1-j}
\asymp_j T^{-j\varepsilon}K_1^2 X_1\;.$$

By the last estimates, combined with (6.20) and (1.18), one obtains:
$$F_{\psi_1,\psi_2,\nu}(z)
={1\over\left|\psi_2\right|^2}\sum_{\xi\in{\frak O}}W_{\eta}\left({|\varphi\xi|^2\over\Xi}\right) 
\hat{f}_{\psi_1,\psi_2,z}\left({\xi\over\psi_2}\right)
{\rm e}\left(\Re \left({\psi_1^{*}\nu\xi\over\psi_2}\right)\right)
+O_{\eta , j}\left( T^{1-j\varepsilon}|\varphi|^{-2}\Xi\right) ,\eqno(6.23)$$
with the function $W_{\eta}(u)$ being as described below (5.27).

Since $f_{\psi_1,\psi_2,z}(s)$ (and hence also its Fourier transform) is identically zero
unless $\psi_1$ and $z$ satisfy (6.22), it follows from (6.17), (6.20), (6.2) and (6.23) that, 
when $\bigl|\psi_i\bigr|^2\asymp |\varphi|^{-2} K_2 M_1$ (for $i=1,2$), one will have:
$$\widetilde{U}_{\nu}\left(\psi_1 ,\psi_2\right)
=\sum_{\xi} W_{\eta}\left({|\varphi\xi|^2\over\Xi}\right) 
\hat{g}_{\psi_1,\psi_2,\xi}(\nu)
\,{\rm e}\!\left(\Re \left({\psi_1^{*}\nu\xi\over\psi_2}\right)\right)
+O_{\eta , j}\left( K_1^{-2}T^{2-(j-1)\varepsilon}\right)\quad\hbox{for $j\geq 2\;$,}\eqno(6.24)$$
where   
$$g_{\psi_1,\psi_2,\xi}(z)
={1\over\left|\psi_2\right|^2}\,\hat{f}_{\psi_1,\psi_2,z}\left({\xi\over\psi_2}\right) 
=\int_{\Bbb C} f_{\psi_1,\psi_2,z}\left(\psi_2 z_2\right) 
{\rm e}\left( -\Re \left(\xi z_2\right)\right) {\rm d}_{+}z_2\;,$$ 
so that, by (6.19), 
$$\hat{g}_{\psi_1,\psi_2,\xi}(\nu)
=\int_{\Bbb C}\int_{\Bbb C} f_{\psi_1,\psi_2,z_1}\left(\psi_2 z_2\right)
{\rm e}\left( -\Re \left(\xi z_2 +\nu z_1\right)\right) {\rm d}_{+}z_2\,{\rm d}_{+}z_1 
=L\left(\psi_1,\psi_2;\nu,\xi\right) ,\eqno(6.25)$$
with $L(\psi,\psi';\nu,\xi)$ defined as in (6.3). 

By (6.12), (6.14), 
(6.24) (with $j=[2/\varepsilon]+2$) and (6.25), it follows that 
$${\cal D}'={\cal D}^{\star}+{\cal E}^{\star} 
+{\cal E}'+O_{\eta,\varepsilon}\left( {\cal E}''\right) ,\eqno(6.26)$$
where ${\cal D}^{\star}$ and ${\cal E}^{\star}$ are as defined in (6.5)-(6.8),   
while, by (1.18)-(1.20), in combination with the properties of $W_{\eta}(u)$ 
(described below (5.27)), (1.24), (6.15), (6.1), the bound (2.13) and 
the hypotheses that $\varepsilon\in(0,1/6]$ and $\eta\in(0,(\log 2)/3]$, one has: 
$$\eqalignno{
{\cal E}'
 &=\sum_{\varphi\neq 0}\sum_{\nu\neq 0}W_{\eta}\!\left({|\varphi\nu|^2\over N}\right)
\sum\sum\!\!\!\!\!\!\!\!
\sum_{\!\!\!\!\!\!\!\!\!\!\!{\scriptstyle\kappa_2\ \ \mu_1\ \ \kappa_4\ \ \mu_2
\atop\scriptstyle (\kappa_2\mu_1 , \kappa_4\mu_2) \sim\varphi}}
\!\!\!\!\!\!\!\!\sum w_2\!\left(\left|\kappa_2\right|^2\right) a\left(\mu_1\right)
\overline{w_2\!\left(\left|\kappa_4\right|^2\right) a\left(\mu_2\right)}
\,V_{\nu}\!\left({\kappa_2\mu_1\over\varphi} , {\kappa_4\mu_2\over\varphi}\right) \ll \cr 
 &\ll\sum_{\varphi\neq 0}\sum_{\scriptstyle \nu\neq 0\atop\scriptstyle |\varphi\nu|^2\leq e^{\eta}N}
\qquad\quad\ \,\sum\!\!\!\!\!\!\!\!\!\! 
\sum_{\!\!\!\!\!\!\!\!\!\!\!\!\!\!\!\!\!\!\!\!\!\!K_2 e^{-\eta}\leq\left|\kappa_2\right|^2, 
\,\left|\kappa_4\right|^2\leq e^{\eta} K_2} 
\qquad\quad\ \sum\!\!\!\!\!\!\!\!\!\!\!\!\sum_{\!\!\!\!\!\!\!\!\!\!\!\!\!\!\!\!\!\!\!{\scriptstyle 
M_1 e^{-\eta}\leq\left|\mu_1\right|^2,\,\left|\mu_2\right|^2\leq e^{\eta} M_1\atop\scriptstyle 
(\kappa_2\mu_1 , \kappa_4\mu_2) \sim\varphi}} 
\left|a\left(\mu_1\right) a\left(\mu_2\right)\right| \times \cr 
 &\quad\ \times \sum_{\scriptstyle K_1 e^{-\eta}\leq |\kappa_1|^2\leq e^{\eta} K_1 
\atop\scriptstyle \kappa_1\kappa_2\mu_1\equiv\varphi\nu\bmod \kappa_4\mu_2 {\frak O}}
\!\!\!\!\!\!\eta^{-1}T^{\varepsilon +1/2} \ll \cr
 &\ll \eta^{-1} T^{\varepsilon +1/2}
\ \,\sum\quad\ \sum\!\!\!\!\!\sum_{\!\!\!\!\!\!\!\!\!\!\!\!\!\!\!\!\!\!\!\!\!\!\!\!\!\!\!\!\!\!\!\!\!\!\!\!\!\!{\!\!\!\!\scriptstyle 
\kappa_1\qquad\,\kappa_2\qquad\,\mu_1\atop\scriptstyle
1/2\leq\left|\kappa_1\kappa_2\mu_1\right|^2/\left( K_1 K_2 M_1\right)\leq 2}}
\quad\ \ \,\sum\ \sum\!\!\!\!\!\!\!\!\!\!\!\!\!\!\!\!
\sum_{\!\!\!\!\!\!\!\!\!\!\!\!\!\!\!\!\!\!{\!\!\!\!\!\!\!\!\!\scriptstyle\kappa_3\quad\,\kappa_4\quad\,\mu_2\atop\scriptstyle 
0<\left|\kappa_1\kappa_2\mu_1 -\kappa_3\kappa_4\mu_2\right|^2\leq  2N
}} \left(\left|a\left(\mu_1\right)\right|^2+\left|a\left(\mu_2\right)\right|^2\right) \ll \cr 
 &\ll \eta^{-1} T^{\varepsilon +1/2}\sum_{e^{-\eta}M_1\leq |\mu|^2\leq e^{\eta}M_1} 
\left| a(\mu)\right|^2
\sum_{|\alpha|^2\asymp K_1 K_2} \tau_2(\alpha)
\sum_{\scriptstyle\beta\neq 0\atop\scriptstyle 0<|\beta -\mu\alpha|^2\ll N} \tau_3(\beta) \ll \cr 
 &\ll  \eta^{-1} T^{\varepsilon +1/2}\,\sum_{\mu}\left| a(\mu)\right|^2
O_{\varepsilon}\left( \left( K_1 K_2\right)^{1+\varepsilon}\right)
\sum_{0<|\delta|^2\ll N} O_{\varepsilon}\left( 
\left( K_1 K_2 M_1\right)^{\varepsilon}\right) = \cr 
 &=O_{\varepsilon}\!\!\left( \eta^{-1} T^{\varepsilon +1/2}\left( K_1 K_2\right)^{1+2\varepsilon}M_1^{\varepsilon}N\right)
\sum_{\mu} |a(\mu)|^2
\ll_{\varepsilon , \eta} T^{5\varepsilon -1/2}K_1^2 K_2^2 M_1 \| a\|_2^2 &(6.27)}$$  
and
$$\eqalignno{
{\cal E}''
 &=\sum_{\varphi\neq 0}\sum_{\nu\neq 0}W_{\eta}\!\left({|\varphi\nu|^2\over N}\right)
\sum\sum\!\!\!\!\!\!\!\!
\sum_{\!\!\!\!\!\!\!\!\!\!\!{\scriptstyle\kappa_2\ \ \mu_1\ \ \kappa_4\ \ \mu_2
\atop\scriptstyle (\kappa_2\mu_1 , \kappa_4\mu_2) \sim\varphi}}
\!\!\!\!\!\!\!\!\sum \left| w_2\!\left(\left|\kappa_2\right|^2\right) a\left(\mu_1\right)
w_2\!\left(\left|\kappa_4\right|^2\right) a\left(\mu_2\right)\right| K_1^{-2} \ll \cr
 &\ll K_1^{-2}\sum_{0<|\nu|^2\ll N} 
\ \sum_{\left|\kappa_2\right|^2\asymp K_2} 
\ \sum_{\left|\kappa_4\right|^2\asymp K_2}
\ \sum_{\left|\mu_1\right|^2\asymp M_1} 
\ \sum_{\left|\mu_2\right|^2\asymp M_1}
\left(\left|a\left(\mu_1\right)\right|^2+\left|a\left(\mu_2\right)\right|^2\right) \ll \cr
 &\ll K_1^{-2} N K_2^2 M_1\sum_{\mu}\left| a(\mu)\right|^2
= T^{\varepsilon -1} K_1^{-1} K_2^3 M_1^2 \| a\|_2^2\;. &(6.28) }$$
Given that (1.18) holds, the result (6.4) of the lemma now follows  
by virtue of the equations (6.11) and (6.26), 
and the bounds noted in (6.27) and (6.28) 
\quad$\square$ 

\bigskip 

\goodbreak
\proclaim{\Smallcaps Lemma~21}. Let ${\cal D}^{\star}$ be given by the equations 
(6.5), (6.7) and (6.3) of Lemma~20. Then, for all $A>0$, one has 
$${\cal D}^{\star}
={\cal D}_1^{\star}+{\cal D}_2^{\star}
+O_{\eta,\varepsilon,A}\left( T^{-A} K_2^2 M_1 \| a\|_2^2\right) ,
\eqno(6.29)$$
where ${\cal D}_1^{\star}$ and ${\cal D}_2^{\star}$ are as defined in 
the equations (1.23)-(1.28) of Theorem~2. 

\medskip 

\goodbreak
\noindent{\Smallcaps Proof.}\  
We begin by reformulating (at the cost of having to introduce a small error term)  
the factor of form $G_{\varphi}\bigl(\psi_1,\psi_2;0)$ occurring in the summand 
on the right-hand side of Equation~(6.5). 
In doing so, we may assume that 
$$e^{-2\eta}K_2 M_1\leq\left|\varphi\psi_i\right|^2\leq e^{2\eta}K_2 M_1\qquad\hbox{for $\,i=1,2$}
\eqno(6.30)$$
(for, by (1.19) and (1.20), the coefficient of 
$G_{\varphi}\bigl(\kappa_2\mu_1/\varphi , \kappa_4\mu_2/\varphi ; 0)$ 
in (6.5) will otherwise be zero). 
Now observe that the properties of $W_{\eta}(u)$
(outlined after (5.27)) ensure that (6.7) implies
$$\left| G_{\varphi}\left(\psi_1 ,\psi_2 ; 0\right)
-\sum_{\nu\neq 0}L\left(\psi_1,\psi_2;\nu,0\right)\right|
\leq\sum_{\scriptstyle\nu\atop\scriptstyle N<|\varphi\nu|^2}
\left|L\left(\psi_1,\psi_2;\nu,0\right)\right|\;.\eqno(6.31)$$
Here, upon using a change of variable of integration to rewrite 
the case $\xi =0$ of (6.3), one finds that
$$L\left(\psi_1,\psi_2;\nu,0\right) 
=T\int_{\Bbb C} \widehat{w_0\circ{\frak N}}\left({\sqrt{T}\nu\over\psi_1\psi_2 z_2}\right)
w_1\left(\left|\psi_2 z_2\right|^2\right)
\overline{w_1\left(\left|\psi_1 z_2\right|^2\right)}
{\rm d}_{+}z_2\qquad\qquad\hbox{($\nu\in{\frak O}$).}\eqno(6.32)$$
By (1.19), one has $|w_1(|\psi_2 z_2|^2)|>0$   
only when $|\psi_2 z_2|^2\asymp K_1$; by this, in combination with 
the result (5.7) of Lemma~15, the case $i=0$ of (1.19), 
the hypothesis that $K_0=1$, the inequalities (6.30) and 
the definition (1.24), it follows that when  $|w_1(|\psi_2 z_2|^2)|>0$ 
one has   
$$\widehat{w_0\circ{\frak N}}\left({\sqrt{T}\nu\over\psi_1\psi_2 z_2}\right)
\ll_j\,\left|{\eta\sqrt{T}\nu\over\psi_1\psi_2 z_2}\right|^{-2j}
\ll_{\eta,j} \left({T |\varphi\nu|^2\over K_1 K_2 M_1}\right)^{\!\!-j} 
=\left({T^{\varepsilon}|\varphi\nu|^2\over N}\right)^{\!\!-j} 
\qquad\hbox{($0\neq\nu\in{\frak O}$, $j\in{\Bbb N}$).}\eqno(6.33)$$
Therefore, given (6.32) and (1.19) (and subject to (6.30) holding), 
one has:
$$L\left(\psi_1,\psi_2;\nu,0\right) 
\ll_{\eta,j} {|\varphi|^2 T K_1\over K_2 M_1}\left({T^{\varepsilon}|\varphi\nu|^2\over N}\right)^{-j}
\qquad\ \hbox{for all $\,j\in{\Bbb N}$.}$$
By applying these estimates to the right-hand side of (6.31) one finds that, 
for $j\geq 2$, 
$$\left| G_{\varphi}\left(\psi_1 ,\psi_2 ; 0\right)
-\sum_{\nu\neq 0}L\left(\psi_1,\psi_2;\nu,0\right)\right|
\ll_{\eta,j}
{|\varphi|^2 T K_1\over K_2 M_1}\left({T^{\varepsilon}|\varphi|^2\over N}\right)^{\!\!-j}
\left({N\over |\varphi|^2}\right)^{\!\!1-j}
={K_1 N T^{1-j\varepsilon}\over K_2 M_1}
=K_1^2 T^{-(j-1)\varepsilon}\;.$$
By this one obtains (given (1.18), and subject to (6.30) holding):
$$G_{\varphi}\left(\psi_1 ,\psi_2 ; 0\right)
=\sum_{\nu\neq 0}L\left(\psi_1,\psi_2;\nu,0\right) 
+O_{\eta,\varepsilon,A}\left( T^{-A}\right)\qquad\hbox{for $\,A>0\;$.}
\eqno(6.34)$$

Since (6.33) and (1.19) imply that the sum 
$\sum_{\nu\neq 0}\widehat{w_0\circ{\frak N}}(T^{1/2}\nu\psi_1^{-1}\psi_2^{-1} z_2^{-1}) 
\,w_1(|\psi_2 z_2|^2)\,\overline{w_1(|\psi_1 z_2|^2)}$ 
is uniformly convergent for all $z_2\in{\Bbb C}$ such that 
$|z_2|^2\asymp |\psi_2|^{-2} K_1$, 
one can deduce from (6.32) and (6.34) (via a change in 
the order of summation and integration) that, for $A>0$, 
$$\eqalign{ 
G_{\varphi}\left(\psi_1,\psi_2;0\right) 
+L\left(\psi_1,\psi_2;0\right) 
 &=T\int_{\Bbb C} 
\left(\sum_{\nu}\widehat{w_0\circ{\frak N}}\biggl( 
{\sqrt{T}\nu\over\psi_1\psi_2 z_2}\biggr)\right) 
w_1\!\!\left(\left|\psi_2 z_2\right|^2\right)
\overline{w_1\!\!\left(\left|\psi_1 z_2\right|^2\right)}{\rm d}_{+}z_2 \ + \cr 
 &\quad +O_{\eta,\varepsilon,A}\left( T^{-A}\right) .}\eqno(6.35)$$ 
Here, by means of the identity noted in (5.8), we are able to apply 
Poission summation (i.e. the case $\tau =0$ of the result (5.2) of Lemma~14) 
so as to obtain: 
$$\sum_{\nu}\widehat{w_0\circ{\frak N}}\biggl( 
{\sqrt{T}\nu\over\psi_1\psi_2 z_2}\biggr) 
=\left| {\psi_1\psi_2 z_2\over\sqrt{T}}\right|^2 
\sum_{\alpha}\left( w_0\circ{\frak N}\right) 
\left( {\psi_1\psi_2 z_2 \alpha\over\sqrt{T}}\right) 
={\left|\psi_1\psi_2 z_2\right|^2\over T}\sum_{\alpha} 
w_0\left( {\left|\psi_1\psi_2 z_2 \alpha\right|^2\over T}\right) .$$ 
By this, (1.19) and the hypothesis $K_0 =1$, it follows that the 
integral on the right-hand side of Equation~(6.35) is trivially equal to zero if 
$\max\{ |\psi_1|^2 , |\psi_2|^2\} > e^{2\eta} T K_1^{-1}$. 
Therefore, given the inequalities in (6.30), it is certainly 
the case that the integral in (6.35) is zero whenever $\varphi$ 
satisfies $|\varphi|^2<e^{-4\eta} T^{-1} K_1 K_2 M_1$. 
This means that when $0<A<\infty$ one has  
$$G_{\varphi}\left(\psi_1,\psi_2;0\right) 
= - L\left(\psi_1,\psi_2;0\right) 
+O_{\eta,\varepsilon,A}\left( T^{-A}\right)\qquad\  
\hbox{if $\ 0<|\varphi|^2<e^{-4\eta} T^{-\varepsilon} N$.}\eqno(6.36)$$ 
Note now that, due to the presence of the 
factor $W_{\eta}\bigl( |\varphi\nu|^2 /N\bigr)$
in (6.7), it is effectively a condition of the summation over $\varphi$ 
in (6.5) that $\varphi$ must satisfy 
$$|\varphi|^2\leq e^{\eta} N\;.\eqno(6.37)$$
Consequently it is only when $|\varphi|^2$ lies in the interval  
$[e^{-4\eta}T^{-\varepsilon}N , e^{\eta} N]$
that one has to refrain from using the result (6.36), 
and make do with (6.34) instead. The combination of 
(6.34) and (6.36) completes our 
reformulation of $G_{\varphi}(\psi_1,\psi_2;0)$. 

Following a change of the variable of integration, 
the case $\nu =0$ of (6.32) gives:
$$L\left(\psi_1,\psi_2;0,0\right) 
=\widehat{w_0\circ{\frak N}}(0)\,T\,|\varphi|^2\int_{\Bbb C}
w_1\left(\left|\varphi\psi_2 z\right|^2\right)
\overline{w_1\left(\left|\varphi\psi_1 z\right|^2\right)}
{\rm d}_{+}z\qquad\qquad\hbox{($0\neq\varphi\in{\frak O}$).}$$ 
Given that (subject to (6.30) holding) one has (6.34) and (6.36), 
and given the last equation (above), the observation (6.37)  and the hypothesis (1.19), 
it follows (by virtue of the arithmetic-geometric mean inequality) that the 
sum ${\cal D}^{\star}$ defined in (6.5) satisfies 
$${\cal D}^{\star}
={\cal D}_1''+{\cal D}_2'' 
+O_{\eta,\varepsilon,A}\left( T^{-A} K_2^2 M_1\sum_{\mu} |a(\mu)|^2\right) \qquad\qquad 
\hbox{($0\leq A<\infty$),}\eqno(6.38)$$
with 
$$\eqalign{ 
{\cal D}_1''
 &=4\quad\ \,\sum\!\!\!\!\!\!\sum_{\!\!\!\!\!\!\!\!\!\!\!\!\!\!\!\!\!\!\!{\scriptstyle 
\beta_1\neq 0\quad\ \beta_2\neq 0\atop\scriptstyle 
\left|\left(\beta_1,\beta_2\right)\right|^2<e^{-4\eta}T^{-\varepsilon}N}}
\Biggl( 
\sum\sum_{\!\!\!\!\!\!\!\!\!\!\!{\scriptstyle\kappa_2\quad\mu_1\atop\scriptstyle 
\kappa_2\mu_1 =\beta_1}}
w_2\left( \left|\kappa_2\right|^2\right) a\left(\mu_1\right)\Biggr)   
\overline{\Biggl(\sum\sum_{\!\!\!\!\!\!\!\!\!\!\!{\scriptstyle 
\kappa_4\quad\mu_2\atop\scriptstyle \kappa_4\mu_2 =\beta_2}}
w_2\left(\left|\kappa_4\right|^2\right) a\left(\mu_2\right)\Biggr) } \ \times \cr 
 &\qquad\times\left( - \widehat{w_0\circ{\frak N}}(0)\right) T  
\left|\left(\beta_1,\beta_2\right)\right|^2 
\int_{\Bbb C} w_1\left( \left|\beta_2 z\right|^2\right) 
\,\overline{w_1\left(\left|\beta_1 z\right|^2\right)}\,{\rm d}_{+}z}$$
and 
$${\cal D}_2''
=\qquad\,\sum_{\!\!\!\!\!\!\!\!\!\!\!\!\!\!e^{-4\eta}T^{-\varepsilon}N\leq|\varphi|^2\leq e^{\eta}N}
\sum\sum\!\!\!\!\!\!\!\!
\sum_{\!\!\!\!\!\!\!\!\!\!\!{\scriptstyle\kappa_2\ \ \mu_1\ \ \kappa_4\ \ \mu_2
\atop\scriptstyle (\kappa_2\mu_1 , \kappa_4\mu_2) \sim\varphi}}
\!\!\!\!\!\!\!\!\sum w_2\!\left(\left|\kappa_2\right|^2\right) a\left(\mu_1\right)
\overline{w_2\!\left(\left|\kappa_4\right|^2\right) a\left(\mu_2\right)}
\ \lambda\!\left({\kappa_2\mu_1\over\varphi} , {\kappa_4\mu_2\over\varphi} \right) , 
\eqno(6.39)$$
where 
$$\lambda\!\left(\psi_1 , \psi_2 \right) 
=\sum_{\nu\neq 0} L\left(\psi_1,\psi_2;\nu,0\right)\;.\eqno(6.40)$$

Given the equality 
$\widehat{w_0\circ{\frak N}}(0)=\pi\int_0^{\infty} w_0(x) {\rm d}x\,$ 
(for which see (7.72), below), it may be deduced from the above that in (6.38) one has 
$${\cal D}_1''={\cal D}_1^{\star}\;,\eqno(6.41)$$ 
with ${\cal D}_1^{\star}$ given by the equations (1.25) and (1.23) of Theorem~2. 

We now have only to attend to the term ${\cal D}_2''$ defined in (6.39)-(6.40). 
Each term appearing in the sum on the right of (6.40) 
may be expressed, via the equation (6.32), as an integral; 
upon applying the result (5.6) of Lemma~15 (twice) to 
the factor $\widehat{w_0\circ{\frak N}}(T^{1/2}\nu\psi_1^{-1}\psi_2^{-1}z_2^{-1})$ in
the integrand in (6.32), and then recalling the definition of the Fourier transform, 
one finds (thereby) that
$$L\left(\psi_1,\psi_2;\nu,0\right) 
={\left|\psi_1\psi_2\right|^4\over |2\pi\nu|^4 T}
\int_{\Bbb C}\int_{\Bbb C}\left(\Delta_{\Bbb C}^2\!\left( w_0\circ{\frak N}\right)\right)\!(s)
\,{\rm e}\!\left(\! -\Re \!\left({\sqrt{T}\nu s\over\psi_1\psi_2 z_2}\right)\!\right)\!{\rm d}_{+}s 
\ w_1\!\!\left(\left|\psi_2 z_2\right|^2\right)
\overline{w_1\!\!\left(\left|\psi_1 z_2\right|^2\right)}
\left| z_2\right|^4 {\rm d}_{+}z_2 .$$
Substituting here $z_2 = -\varphi z$, and then summing over non-zero $\nu\in{\frak O}$,
one obtains (given (6.40)):   
$$\lambda\!\left(\psi_1 , \psi_2\right)
={|\varphi|^6\left|\psi_1\psi_2\right|^4\over (2\pi)^4 T}
\int_{\Bbb C}\int_{\Bbb C}
\left(\Delta_{\Bbb C}^2\!\left( w_0\circ{\frak N}\right)\right)\!(s)
\ w_1\!\left(\left|\varphi\psi_2 z\right|^2\right)
\overline{w_1\!\left(\left|\varphi\psi_1 z\right|^2\right)}
|z|^4\,p\!\left(\!{\sqrt{T} s\over \varphi\psi_1\psi_2 z}\!\right) 
{\rm d}_{+}s\,{\rm d}_{+}z\;,$$
where $p(\alpha)$ is as defined in Theorem~2, Equation~(1.26). 
Since $p(i\alpha)=p(\alpha)$, and since 
$$\left(\Delta_{\Bbb C}^2\left( w_0\circ{\frak N}\right)\right)(s) 
=\left(\Delta_{\Bbb C}^2\left( w_0\circ{\frak N}\right)\right)(|s|) 
=\left( {\partial^2\over\partial r^2}+r^{-1}\,{\partial\over\partial r}\right)^{\!\!2}  
w_0\!\left( r^2\right)\bigg|_{r=|s|} 
=\left( 4x\,{{\rm d}^2\over{\rm d}x^2}+4\,{{\rm d}\over{\rm d}x}\right)^{\!\!2} w_0(x) 
\bigg|_{x=|s|^2}\;,$$ 
it follows by (6.39) and the expression just obtained for $\lambda(\psi_1,\psi_2)$ 
that one has 
$$\eqalignno{ 
{\cal D}_2''
 &={1\over 4 \pi^2 T}\qquad\quad\sum\!\!\!\!\!\!\!\!\!\!\sum_{\!\!\!\!\!\!\!\!\!\!\!\!\!\!\!\!\!\!\!\!\!{\scriptstyle 
\beta_1\neq 0\quad\ \beta_2\neq 0\atop\scriptstyle 
e^{-4\eta}T^{-\varepsilon}N\leq\left|\left(\beta_1,\beta_2\right)\right|^2\leq e^{\eta}N}}
\left|\left(\beta_1,\beta_2\right)\right|^2 
\left|\left[\beta_1 , \beta_2\right]\right|^4 \ \times \cr 
 &\quad\ \times\Biggl( 
\sum\sum_{\!\!\!\!\!\!\!\!\!\!\!{\scriptstyle\kappa_2\quad\mu_1\atop\scriptstyle 
\kappa_2\mu_1 =\beta_1}}
w_2\left( \left|\kappa_2\right|^2\right) a\left(\mu_1\right)\Biggr)   
\overline{\Biggl(\sum\sum_{\!\!\!\!\!\!\!\!\!\!\!{\scriptstyle 
\kappa_4\quad\mu_2\atop\scriptstyle \kappa_4\mu_2 =\beta_2}}
w_2\left(\left|\kappa_4\right|^2\right) a\left(\mu_2\right)\Biggr) } \ \times \cr 
 &\quad\ \times \int_{\Bbb C}\int_{\Bbb C} 
\left(\Delta_{\Bbb C}^2\left( w_0\circ{\frak N}\right)\right)\!(s) 
\,p\!\left( {\sqrt{T} s\over\left[\beta_1 , \beta_2\right] z}\right) {\rm d}_{+}s 
\ w_1\!\left( \left|\beta_2 z\right|^2\right) 
\overline{w_1\!\left(\left|\beta_1 z\right|^2\right)} \,|z|^4\,{\rm d}_{+}z =\cr 
 &={\cal D}_2^{\star}\;, &(6.42)}$$
where ${\cal D}_2^{\star}$ is as defined in Theorem~2, Equations~(1.23), (1.24) 
and (1.26)-(1.28). 
By (6.38), (6.41) and (6.42), the result (6.29) follows 
\quad$\square$ 

\bigskip 

\goodbreak\proclaim{\Smallcaps Lemma~22}. Suppose that  
$$Z=T^{\varepsilon} M_1\;,\eqno(6.43)$$
and that $\Xi$ and ${\cal E}^{\star}$ are given by the equations 
(6.2), (6.3) and (6.6)-(6.8) of Lemma~20.  
Then either 
$${\cal E}^{\star}\ll_{\eta,\varepsilon} T^{6\varepsilon} K_1 K_2^2 M_1 \|a\|_2^2\;,$$ 
or else, for some $\varphi,\Phi_1,\varphi_1,\Phi_2,\varphi_2\in{\frak O}$ and 
some ${\bf z}\in{\Bbb C}^3$, one has: 
$$0<|\varphi|^2\leq e^{\eta}N\;,\eqno(6.44)$$ 
$$\Phi_1\varphi_1 =\Phi_2\varphi_2 =\varphi\;,\eqno(6.45)$$ 
$$\left| z_1\right|^2\ll {|\varphi|^2 T^{\varepsilon}\over N}\;,\qquad 
\left| z_2\right|^2\ll {|\varphi|^2 T^{\varepsilon}\over\Xi}\;,\qquad 
\left| z_3\right|^2\ll {\left|\Phi_1\right|^2 T^{\varepsilon}\over Z}
\eqno(6.46)$$ 
and 
$${\cal E}^{\star}\ll_{\varepsilon} 
T^{1+\varepsilon} \left|\varphi z_2 z_3\right|^2 |Y|\;,$$ 
where 
$$\eqalign{Y &=\sum_{\scriptstyle\rho_1\atop\scriptstyle\left(\rho_1,\varphi_1\right)\sim 1} 
\!\!\!\!\!\sum_{\,{\scriptstyle\rho_2\atop\scriptstyle 
\ \sim\left(\rho_2,\varphi_2\right)}} 
\!\!\!\!\!\!a\!\left(\Phi_1\rho_1\right) \overline{a\!\left(\Phi_2\rho_2\right)} 
\,\sum_{\nu\neq 0} \sum_{\xi\neq 0} 
W_{\eta}\!\!\left({|\varphi\nu|^2\over N}\right) 
W_{\eta}\!\!\left({|\varphi\xi|^2\over\Xi}\right) 
{\rm e}\left(\Re \left( z_1\nu + z_2\xi\right)\right) \,\times \cr
 &\quad\times 
\sum_{\zeta\neq 0} 
W_{\eta}\!\!\left({|\Phi_1\zeta|^2\over Z}\right){\rm e}\left(\Re \left( z_3\zeta\right)\right) 
\sum_{\varpi\neq 0} g_{\rho_1,\rho_2}\!\!\left(\left|\varpi\right|^2\right) 
S\!\left(\nu\xi\rho_1^{*} , \zeta ; \rho_2\varpi\right) 
}\eqno(6.47)$$
with 
$$\eqalignno{g_{\rho_1,\rho_2}(x) 
=g_{\rho_1,\rho_2}\!\left(\varphi_1,\varphi_2;z_1,z_2,z_3;x\right)   
 &={\left| z_1 z_2 z_3\rho_1\right|^2\left|\rho_2\right|^4 x^2\over T} 
\,{w}_0\!\left({\left| z_1 z_2 z_3\rho_1\right|^2\left|\rho_2\right|^4 x^2\over T}\right) 
w_1\!\left(\left| z_2\rho_2\right|^2 x\right) \times\qquad\ \cr  
 &\quad\,\times\,\overline{w_1\!\left(\left| z_2 z_3\rho_1\rho_2\right|^2 x\right)} 
\,w_2\!\left(\left|\varphi_1 z_3\rho_2\right|^2 x\right) 
\overline{w_2\!\left(\left|\varphi_2\right|^2 x\right)}\;, &(6.48)}$$ 
and with $S(\alpha,\beta;\gamma)$ denoting the Kloosterman sum defined by  
the equation (5.5) of Lemma~14. 

\medskip 

\goodbreak
\noindent{\Smallcaps Proof.}\  
We first clarify the effect of the condition
$\bigl(\kappa_2\mu_1,\kappa_4\mu_2\bigr)\sim\varphi$ in (6.6), by defining new  
${\frak O}$-valued  variables $\Phi_1,\Phi_2,\varphi_1,\varphi_2$,  
which are rendered dependent by being required to satisfy  
$$\Phi_1\sim\left(\mu_1 , \varphi\right) ,\qquad 
\Phi_2\sim\left(\mu_2 , \varphi\right)$$
and the condition (6.45). 
These relations imply that 
$$\mu_1 =\Phi_1\rho_1\qquad\hbox{and}\qquad 
\mu_2 =\Phi_2\rho_2\;,$$
where $\rho_1,\rho_2\in{\frak O}$ are such that 
$$\left(\rho_1 , \varphi_1\right)\sim \left(\rho_2 , \varphi_2\right)\sim 1\;.\eqno(6.49)$$
It follows that if $\bigl(\kappa_2\mu_1,\kappa_4\mu_2\bigr)\sim\varphi$, 
then $\varphi_1\mid\kappa_2$ and $\varphi_2\mid\kappa_4$, and so, 
for a unique pair $\varpi_1,\varpi_2\in{\frak O}$, one has 
$$\kappa_{2} =\varphi_1 \varpi_1\qquad\hbox{and}\qquad 
\kappa_{4} =\varphi_2 \varpi_2\;.$$
Given (6.45) and 
the last two equations above, one has (simultaneously) 
$$\bigl(\kappa_2\mu_1,\kappa_4\mu_2\bigr)\sim\varphi\;,\quad 
\bigl(\mu_1 , \varphi\bigr)\sim\Phi_1\quad\hbox{and}\quad \bigl(\mu_2 , \varphi\bigr)\sim\Phi_2$$
if and only if the dependent variables $\varphi_1,\varphi_2,\rho_1,\rho_2,\varpi_1,\varpi_2$ 
satisfy both (6.49) and 
$$\left( \varpi_1\rho_1 , \varpi_2\rho_2\right)\sim 1\;.$$
\par 
Using the above to rewrite the right-hand side of equation (6.6), one obtains:  
$$\eqalign{ 
{\cal E}^{\star} &={1\over 16}\sum_{\varphi\neq 0}
\sum\!\!\!\sum_{\!\!\!\!\!\!\!\!\!\!{\scriptstyle\Phi_1\ \ \varphi_1\atop\scriptstyle 
\ \ \Phi_1\varphi_1=\varphi}}
\!\!\!\sum\sum_{\!\!\!\!\!\!\!\!\!\!{\scriptstyle\Phi_2\ \ \,\varphi_2\atop\scriptstyle 
\!\!=\Phi_2\varphi_2}}
\sum_{\scriptstyle\rho_1\atop\scriptstyle\left(\rho_1,\varphi_1\right)\sim 1} 
\!\!\!\!\!\sum_{\,{\scriptstyle\rho_2\atop\scriptstyle 
\ \sim\left(\rho_2,\varphi_2\right)}} 
\!\!\!\!\!\!a\!\left(\Phi_1\rho_1\right) \overline{a\!\left(\Phi_2\rho_2\right)}
\,\sum_{\varpi_2}\overline{w_2\!\!\left(\left|\varphi_2\varpi_2\right|^2\right)} \ \times \cr 
 &\quad\ \times\sum_{\varpi_1} w_2\!\!\left(\left|\varphi_1\varpi_1\right|^2\right)
H_{\varphi}\!\left(\rho_1\varpi_1 , \rho_2\varpi_2\right) ,}\eqno(6.50)$$
where $H_{\varphi}(\psi_1,\psi_2)$ is as defined in 
(6.7) and (6.8) and is (therefore) equal to zero   
whenever $\psi_1,\psi_2\in{\frak O}$ fail to satisfy the condition 
$(\psi_1,\psi_2)\sim 1$. 
\par 
Since $W_{\eta}(u)$ is equal to zero for $u\geq e^{\eta}$, the definitions 
(6.7) and (6.8) imply that $H_{\varphi}(\psi_1,\psi_2)$ is (trivially) 
equal to zero whenever one has $|\varphi|^2\geq e^{\eta}N$. This, combined 
with the relations 
$$\sum_{\varphi\neq 0} 
\,{\tau_2^2(\varphi)\over |\varphi|^{2+2\varepsilon}} 
\ll_{\varepsilon}\sum_{\varphi\neq 0}\,{1\over |\varphi|^{2+\varepsilon}}  
\ll_{\varepsilon} 1$$ 
enables us to deduce from (6.50) that, for some 
$\varphi,\Phi_1,\varphi_1,\Phi_2,\varphi_2\in{\frak O}$ satisfying 
(6.44) and (6.45), one has 
$${\cal E}^{\star}\ll_{\varepsilon}
|\varphi|^{2+2\varepsilon}\Biggl| 
\sum_{\scriptstyle\rho_1\atop\scriptstyle\left(\rho_1,\varphi_1\right)\sim 1} 
\!\!\!\!\!\sum_{\,{\scriptstyle\rho_2\atop\scriptstyle 
\ \sim\left(\rho_2,\varphi_2\right)}} 
\!\!\!\!\!\!a\!\left(\Phi_1\rho_1\right) \overline{a\!\left(\Phi_2\rho_2\right)}
\,\sum_{\varpi_2}\overline{w_2\!\!\left(\left|\varphi_2\varpi_2\right|^2\right)} 
\,Q\!\left(\rho_1 , \rho_2\varpi_2\right)\Biggr|\;,\eqno(6.51)$$
where 
$$Q\!\left(\rho , \psi\right) 
=\sum_{\varpi_1}w_2\!\!\left(\left|\varphi_1\varpi_1\right|^2\right)
H_{\varphi}\!\left(\rho\varpi_1 , \psi\right) .\eqno(6.52)$$
Accordingly, we assume henceforth (as we may) that 
$\varphi$, $\Phi_1$, $\varphi_1$, $\Phi_2$ and $\varphi_2$ are given 
Gaussian integers satisfying   
(6.44) and (6.45), and that what is stated in (6.51) and (6.52) is valid.  
\par 
By (6.52), (6.8), (6.7) and (6.3) it follows that, for 
$\rho,\psi\in{\frak O}-\{ 0\}$, one has 
$$Q\!\left(\rho , \psi\right) 
=\sum_{\nu\neq 0}\sum_{\xi\neq 0}W_{\eta}\left({|\varphi\nu|^2\over N}\right)
W_{\eta}\left({|\varphi\xi|^2\over\Xi}\right) 
\int_{\Bbb C}\int_{\Bbb C} R\left(\nu\xi , \rho ; \psi ; z_1 , z_2\right) 
{\rm e}\left(\Re \left(\nu z_1+\xi z_2\right)\right){\rm d}_{+}z_1\,{\rm d}_{+}z_2\;,
\eqno(6.53)$$
where 
$$R\left(\beta , \rho ; \psi ; z_1 , z_2\right) 
=\sum_{\varpi_1\in{\frak O}} r_{\rho,\psi,z_1,z_2}(\varpi_1)\,
{\rm e}\!\left(\Re \left({\beta\rho^{*}\varpi_1^{*}\over\psi}\right)\right) , 
\eqno(6.54)$$
with
$$r_{\rho,\psi,z_1,z_2}(s)
=w_1\!\!\left(\left| z_2\psi\right|^2\right)\left| z_2 \psi\rho s\right|^2 
w_0\!\!\left({\left| z_1 z_2\psi\rho s\right|^2\over T}\right) 
\overline{w_1\!\!\left(\left| z_2\rho s\right|^2\right)} 
\,w_2\!\left( |\varphi_1 s|^2\right)\qquad\quad\hbox{($s\in{\Bbb C}$).}\eqno(6.55)$$ 
Therefore, given (1.19), one may apply the result (5.4) of 
Lemma~14 to the sum on the right-hand side of Equation~(6.54), so as to obtain:
$$R\left(\beta , \rho ; \psi ; z_1 , z_2\right) 
={1\over |\psi|^2}\sum_{\zeta\in{\frak O}}
\hat{r}_{\rho,\psi,z_1,z_2}\!\left({\zeta\over\psi}\right)
S\left( \beta \rho^{*} , \zeta ; \psi\right) .\eqno(6.56)$$
\par 
The condition (1.19) applies to each of the functions $w_0$, $w_1$ and $w_2$,  
and so, as is evident from a comparison of (6.55) with the definition (6.19) of 
the function $f_{\psi_1,\psi_2,z}(\kappa)$ considered in our proof of Lemma~20, 
one may bound $\hat{r}_{\rho,\psi,z_1,z_2}(\zeta/\psi)$  
through steps similar to the steps described in the paragraph below (6.22) 
(where we have estimated the Fourier transform $\hat{f}_{\psi_1,\psi_2,z}$).  
Indeed, in place of what was noted (concerning the function 
$f_{\psi_1,\psi_2,z}$) in (6.21) and (6.22), 
we are now able to note that, for each fixed choice of $\rho$, $\psi$, $z_1$ and $z_2$, one has 
$r_{\rho,\psi,z_1,z_2}(z)=h\bigl(|z|^2\bigr)\,$ (say), where the function  
$h : [0,\infty)\rightarrow{\Bbb C}$ is such that, 
for $x>0$ and $k\in{\Bbb N}\cup\{ 0\}$, one has  
$$h^{(k)}(x)=\cases{O_k\!\left( K_1 |\psi|^2 (\eta x)^{-k}\right) 
 &if $e^{-\eta}|\varphi_1|^{-2} K_2\leq x\leq e^{\eta}|\varphi_1|^{-2} K_2$ , \cr
0 &otherwise,}$$
and where $h$ is identically equal to zero unless one has both
$$e^{-2\eta}T/K_1\leq\left|\psi z_1\right|^2\leq e^{2\eta}T/K_1\qquad 
{\rm and}\qquad 
e^{-\eta} K_1\leq\left|\psi z_2\right|^2\leq e^{\eta} K_1\;.\eqno(6.57)$$
Hence, by means of the result (5.7) of Lemma~15, one finds that if $\zeta\neq 0$ then, 
for all $j\in{\Bbb N}\cup\{ 0\}$, 
$$\eqalign{ 
\left|\hat{r}_{\rho,\psi,z_1,z_2}\left({\zeta\over\psi}\right)\right|
 &\ll |\psi|^{2j} |\zeta|^{-2j}\,O_{j}\!\left(\left|\varphi_1\right|^{-2} K_2 K_1 |\psi|^2 
\left( \eta\left|\varphi_1\right|^{-1} K_2^{1/2}\right)^{-2j}\right) \ll_{\eta,j}\cr 
 &\ll_{\eta,j} K_1 \left( \left|\varphi_1\right|^{-2} K_2\right)^{\!\!1-j} 
|\psi|^{2j+2} |\zeta|^{-2j}\;.}$$
By using the last bound above, and estimating 
the Kloosterman sums  in (6.56) trivially (i.e. by the result (5.28) of Lemma~19), 
one finds that  
$$\sum_{|\zeta|^2\geq X}
{\left| \hat{r}_{\rho,\psi,z_1,z_2}(\zeta /\psi) 
\,S\left( \beta \rho^{*} , \zeta ; \psi\right)\right|\over |\psi|^2} 
\ll_{\eta,j}\ {K_1 |\psi|^4 \over T^{(j-1)\varepsilon}}\qquad\  
\hbox{($X\asymp T^{\varepsilon}K_2^{-1}\left|\varphi_1\psi\right|^2\,$ 
and $\,j\geq 2$).}\eqno(6.58)$$
\par 
Since our principal concern (for the remainder of this proof) 
is with the sum over $\rho_1$, $\rho_2$ and $\varpi_2$ on the right-hand side of Equation~(6.51), 
it therefore follows by (1.19), (1.20) and (6.45) that the sum $Q(\rho,\psi)$ defined in  
(6.52) requires our attention only in cases where 
$|\psi|^2\asymp |\varphi|^{-2} K_2 M_1\leq K_2 M_1$ 
and $|\rho|^2\asymp |\Phi_1|^{-2} M_1\leq M_1$. 
In all such cases it follows by (6.56), (6.58), (6.45) and (1.18) that one has, in 
(6.53), 
$$R\left(\beta , \rho ; \psi ; z_1 , z_2\right) 
=\sum_{\zeta}
W_{\eta}\!\!\left({\left|\Phi_1\zeta\right|^2\over Z}\right) |\psi|^{-2} 
\,\hat{r}_{\rho,\psi,z_1,z_2}\left({\zeta\over\psi}\right)
S\left( \beta \rho^{*} , \zeta ; \psi\right) 
+O_{\eta,j}\left( {T^{2-(j-1)\varepsilon} |\psi|^2\over |\rho|^2}\right) ,\eqno(6.59)$$
with $W_{\eta} : [0,\infty)\rightarrow [0,1]$ as described below (5.27), 
and with $Z=T^{\varepsilon} M_1$ (as in (6.43)).  
Moreover, by a change of variable of integration , 
$$|\psi|^{-2}\,\hat{r}_{\rho,\psi,z_1,z_2}\left({\zeta\over\psi}\right)
=\int_{\Bbb C}r_{\rho,\psi,z_1,z_2}\left( -\psi z_3\right)
{\rm e}\left( \Re \left(\zeta z_3\right)\right) {\rm d}_{+}z_3\;.$$
Upon substituting the last expression for 
$|\psi|^{-2}\hat{r}_{\rho,\psi,z_1,z_2}(\zeta/\psi)$
into the case $j=[2/\varepsilon]+3$ of (6.59), and then substituting the 
resulting estimate for $R(\beta,\rho;\psi;z_1,z_2)$ 
into (6.53) (with due application of the
properties of the function $W_{\eta}(u)$ set out below (5.27), 
while also bearing in mind what was noted in, and just above, (6.57)), one obtains:
$$\eqalign{Q(\rho,\psi)
 &=\sum_{\nu\neq 0}\sum_{\xi\neq 0}\sum_{\zeta}
W_{\eta}\!\!\left({|\varphi\nu|^2\over N}\right) 
W_{\eta}\!\!\left({|\varphi\xi|^2\over\Xi}\right) 
W_{\eta}\!\!\left({|\Phi_1\zeta|^2\over Z}\right) 
{\cal I}(\rho,\psi;\nu,\xi,\zeta)\,S\!\left(\nu\xi\rho^{*} , \zeta ; \psi\right) \,+ \cr 
 &\quad +O_{\eta,\varepsilon}\left( 
|\varphi|^{-4} N \Xi \,|\rho\psi|^{-2} T^{1-\varepsilon}\right)\;,}\eqno(6.60)$$
with, for ${\bf s}\in{\Bbb C}^3$, 
$${\cal I}\!\left(\rho,\psi;s_1,s_2,s_3\right) 
={\cal I}(\rho,\psi;{\bf s})
=\int_{\Bbb C}\int_{\Bbb C}\int_{\Bbb C} 
r_{\rho,\psi,z_1,z_2}\left( -\psi z_3\right) 
{\rm e}\left(\Re \left({\bf s}\cdot{\bf z}\right)\right) 
{\rm d}_{+}z_1 {\rm d}_{+}z_2 {\rm d}_{+}z_3\eqno(6.61)$$
where ${\b z}=\bigl( z_1,z_2,z_3\bigr)$, and where, by (6.55),
$$r_{\rho,\psi,z_1,z_2}\left( -\psi z_3\right) 
=\left|\rho\psi^2 z_2 z_3\right|^2 w_0\!\!\left({\left|\rho\psi^2 z_1 z_2 z_3\right|^2\over T}\right) 
w_1\!\left(\left| z_2\psi\right|^2\right) 
\,\overline{w_1\!\!\left(\left|\rho\psi z_2 z_3\right|^2\right)}
w_2\!\!\left(\left|\varphi_1\psi z_3\right|^2\right)\;.\eqno(6.62)$$
\par 
By (6.57), (6.62) and (1.19), it follows that 
$$\left| r_{\rho,\psi,z_1,z_2}\bigl( -\psi z_3\bigr)\right| > 0\quad\ {\rm only\ if}\ \quad 
\left| z_1\right|^2\asymp{T\over \left|\psi\right|^2 K_1}\;,\quad 
\left| z_2\right|^2\asymp{K_1\over |\psi|^2}\quad\hbox{and}\quad  
\left| z_3\right|^2\asymp{1\over |\rho|^2}\;.\eqno(6.63)$$
By (6.61)-(6.63) and (1.19) (again), one has 
$${\cal I}(\rho,\psi;{\bf s})
\ll\left|\rho\psi\right|^{-2} T K_1\qquad\qquad\hbox{($\,\rho,\psi\in{\frak O}-\{ 0\}$, 
$\,{\bf s}\in{\Bbb C}^3$).}
\eqno(6.64)$$
\par 
The result (5.31) of Lemma~19 implies that, 
for all $\nu,\xi,\rho,\psi\in{\frak O}-\{ 0\}$ such that $(\rho,\psi)\sim 1$, one has 
$$\left|S\!\left(\nu\xi\rho^{*} , 0 ; \psi\right)\right|
\leq \left|\left(\nu\xi\rho^{*} , \psi\right)\right|^2
=\left|\left(\nu\xi , \psi\right)\right|^2\;.$$
Therefore,  given (6.64) and the properties of $W_{\eta}(\nu)$ 
summarised below (5.27), 
it follows that the total contribution to the sum in (6.60) arising from
terms with $\zeta =0$ is 
$$\eqalignno{O\left(\sum_{\scriptstyle\nu\neq 0\atop\scriptstyle |\varphi\nu|^2\ll N}
\sum_{\scriptstyle\xi\neq 0\atop\scriptstyle |\varphi\xi|^2\ll \Xi}
|\rho\psi|^{-2} T K_1\left| (\nu\xi , \psi)\right|^2\right) 
 &\ll |\rho\psi|^{-2} T K_1\sum_{0<\left|\lambda\right|^2\ll N\Xi /|\varphi|^4} \tau_2(\lambda)  
\,|(\lambda , \psi)|^2 \,\ll_{\varepsilon} \cr
 &\ll_{\varepsilon} |\rho\psi|^{-2} T K_1 |\varphi|^{-4} (N\Xi)^{1+\varepsilon}\tau_2(\psi)\;.
 &(6.65)}$$
Moreover, given 
(1.24), (6.1), (6.2) and the hypotheses that $K_1\geq 1$, $K_2\geq 1$ and $M_1\geq 1$, 
it follows that the final upper bound in (6.65) is greater than 
what appears inside the brackets of the $O$-term in (6.60). 
The result (6.60) therefore yields the estimate  
$$\eqalign{Q(\rho , \psi)
 &=\sum_{\nu\neq 0}\sum_{\xi\neq 0}\sum_{\zeta\neq 0}
W_{\eta}\!\!\left({|\varphi\nu|^2\over N}\right) 
W_{\eta}\!\!\left({|\varphi\xi|^2\over\Xi}\right) 
W_{\eta}\!\!\left({|\Phi_1\zeta|^2\over Z}\right) 
{\cal I}(\rho,\psi;\nu,\xi,\zeta)\,S\!\left(\nu\xi\rho^{*} , \zeta ; \psi\right) \,+ \cr 
 &\quad +O_{\eta,\varepsilon}\left( |\rho\psi|^{-2} T K_1 |\varphi|^{-4} 
(N\Xi)^{1+\varepsilon}\tau_2(\psi)
\right)\;,}$$
which, when substituted into (6.51), gives
$${\cal E}^{\star} 
\ll_{\varepsilon} |\varphi|^{2+2\varepsilon}\bigl| 
{\cal E}^{\flat}+{\cal E}^{\sharp}\bigr|\;, 
\eqno(6.66)$$
where  
$$\eqalign{{\cal E}^{\flat} 
 &=\sum_{\scriptstyle\rho_1\atop\scriptstyle\left(\rho_1,\varphi_1\right)\sim 1} 
\!\!\!\!\!\sum_{\,{\scriptstyle\rho_2\atop\scriptstyle 
\ \sim\left(\rho_2,\varphi_2\right)}} 
\!\!\!\!\!\!a\!\left(\Phi_1\rho_1\right) \overline{a\!\left(\Phi_2\rho_2\right)}
\sum_{\varpi_2}\overline{w_2\!\!\left(\left|\varphi_2\varpi_2\right|^2\right)} 
\sum_{\nu\neq 0}\sum_{\xi\neq 0}W_{\eta}\!\!\left({|\varphi\nu|^2\over N}\right) 
W_{\eta}\!\!\left({|\varphi\xi|^2\over\Xi}\right)\times \cr
 &\quad\times\sum_{\zeta\neq 0}
W_{\eta}\!\!\left({|\Phi_1\zeta|^2\over Z}\right) 
{\cal I}\!\left(\rho_1,\rho_2\varpi_2;\nu,\xi,\zeta\right) 
\,S\!\left(\nu\xi\rho_1^{*} , \zeta ; \rho_2\varpi_2\right) ,}\eqno(6.67)$$
and 
$$\eqalignno{{\cal E}^{\sharp} 
 &=\sum_{\scriptstyle\rho_1\atop\scriptstyle\left(\rho_1,\varphi_1\right)\sim 1} 
\!\!\!\!\!\sum_{\,{\scriptstyle\rho_2\atop\scriptstyle 
\ \sim\left(\rho_2,\varphi_2\right)}} 
\!\!\!\!\!\!a\!\left(\Phi_1\rho_1\right) \overline{a\!\left(\Phi_2\rho_2\right)}
\sum_{\varpi_2}\overline{w_2\!\!\left(\left|\varphi_2\varpi_2\right|^2\right)} 
\,O_{\eta,\varepsilon}\!\!\left( |\varphi|^{-4} (N\Xi)^{1+\varepsilon} T K_1 
\left|\rho_1\rho_2\varpi_2\right|^{-2}\tau_2\left(\rho_2\varpi_2\right)\right) 
\ll_{\eta,\varepsilon} \cr
 &\ll_{\eta,\varepsilon}
|\varphi|^{-4} (N\Xi)^{1+\varepsilon} T K_1\sum_{\rho_1}\sum_{\rho_2}
{\left| a\!\left(\Phi_1\rho_1\right) a\!\left(\Phi_2\rho_2\right)\right|\over 
\left|\rho_1\rho_2\right|^2} 
\sum_{\varpi_2}\left| w_2\!\!\left(\left|\varphi_2\varpi_2\right|^2\right)\right| 
\left|\varpi_2\right|^{-2}\tau_2\left(\rho_2\varpi_2\right) \ll_{\varepsilon} \cr
 &\ll_{\varepsilon} |\varphi|^{-4} (N\Xi)^{1+\varepsilon} T K_1 
\| a\|_2^2 
\,\Biggl( \sum_{\scriptstyle\rho\atop\scriptstyle 
\left|\Phi_1\rho_1\right|^2\asymp M_1}{1\over\left|\rho_1\right|^4}\Biggr)^{\!\!1/2} 
\Biggl( \sum_{\scriptstyle\rho_2\atop\scriptstyle 
\left|\Phi_2\rho_2\right|^2\asymp M_1}{1\over |\rho_2|^4}\Biggr)^{\!\!1/2} 
\left( M_1 K_2\right)^{\varepsilon} \ll \cr 
 &\ll |\varphi|^{-4} \left|\Phi_1\Phi_2\right| M_1^{-1} 
\left( M_1 K_2\right)^{\varepsilon} (N\Xi)^{1+\varepsilon} T K_1 \| a\|_2^2  
\ll |\varphi|^{-2} 
T^{5\varepsilon} K_1 K_2^2 M_1 \| a\|_2^2 &(6.68) 
}$$
(with the last four estimates following by virtue of the estimate (2.13) and the Cauchy-Schwarz inequality, 
given that we have $\varepsilon\in (0,1/6]$, $K_1,K_2,M_1\in [1,\infty)$,  
$\varphi , \Phi_1 , \varphi_1 , \Phi_2 , \varphi_2\in{\frak O}-\{ 0\}$ 
and (1.18)-(1.20), (1.24), (6.2) and (6.45)). 
\par  
By using (6.61) and (6.62) we are able 
to reformulate (6.67) as the equation 
$${\cal E}^{\flat} 
= \int_{{\Bbb C}^{*}}\int_{{\Bbb C}^{*}}\int_{{\Bbb C}^{*}}
{\cal F}\!\left( z_1,z_2,z_3\right) 
{\rm d}_{\times}z_1
{\rm d}_{\times}z_2
{\rm d}_{\times}z_3\;,\eqno(6.69)$$
in which we have ${\rm d}_{\times}z=|z|^{-2}{\rm d}_{+}z$ and,  
for ${\bf z}=\left( z_1,z_2,z_3\right)\in{\Bbb C}^3$,
$$\eqalignno{{\cal F}\!\left( z_1,z_2,z_3\right)
 &=\sum_{\scriptstyle\rho_1\atop\scriptstyle\left(\rho_1,\varphi_1\right)\sim 1} 
\!\!\!\!\!\sum_{\,{\scriptstyle\rho_2\atop\scriptstyle 
\ \sim\left(\rho_2,\varphi_2\right)}} 
\!\!\!\!\!\!a\!\left(\Phi_1\rho_1\right) \overline{a\!\left(\Phi_2\rho_2\right)} 
\sum_{\nu\neq 0}\sum_{\xi\neq 0}\sum_{\zeta\neq 0}
W_{\eta}\!\!\left({|\varphi\nu|^2\over N}\right)  
W_{\eta}\!\!\left({|\varphi\xi|^2\over\Xi}\right) 
W_{\eta}\!\!\left({|\Phi_1\zeta|^2\over Z}\right) \times \cr 
 &\quad\ \times {\rm e}\left(\Re \left( {\bf z}\cdot (\nu,\xi,\zeta)\right)\right) 
\sum_{\varpi_2} f\!\left( {\bf z} ; \rho_1 , \rho_2 , \varpi_2\right)
\,S\!\left(\nu\xi\rho_1^{*} , \zeta ; \rho_2\varpi_2\right) , &(6.70)
}$$
with  
$$f\!\left( {\bf z} ; \rho_1 , \rho_2 , \varpi\right)
=\left| z_1 z_2 z_3\right|^2 r_{\rho_1,\rho_2\varpi,z_1,z_2}\!\left( -\rho_2\varpi z_3\right) 
\,\overline{w_2\!\!\left(\left|\varphi_2\varpi\right|^2\right)}\;.\eqno(6.71)$$ 
\par 
By (6.71), (6.63) and (1.19), one has 
$f( {\bf z} ; \rho_1 , \rho_2 , \varpi_2) \neq 0$ 
on the right-hand side of Equation~(6.70) only when it is the case that 
$\bigl|\varpi_2\bigr|^2\asymp\bigl|\varphi_2\bigr|^{-2} K_2$ and 
$$\left| z_1\right|^2 \asymp{T\over \left|\rho_2\varpi_2\right|^2 K_1}
\asymp{\left|\varphi_2\right|^2 T\over \left|\rho_2\right|^2 K_1 K_2}\;,\quad 
\left| z_2\right|^2 \asymp{K_1\over \left|\rho_2 \varpi_2\right|^2} 
\asymp{\left|\varphi_2\right|^2 K_1\over\left|\rho_2\right|^2 K_2}\quad 
{\rm and}\quad  
\left| z_3\right|^2 \asymp{1\over\left|\rho_1\right|^2}\;.$$
Furthermore, the hypothesis (1.20) implies that 
$a(\Phi_1\rho_1) \overline{a(\Phi_2\rho_2)}=0$ unless 
$|\rho_i|^2\asymp |\Phi_i|^{-2} M_1$ for $i=1,2$.
Given what has just been noted, and given the condition (6.45) 
and the definitions of $N$, $\Xi$, $Z$ and ${\cal F}(z_1,z_2,z_3)$  
in (1.24), (6.2), (6.43) and (6.70),  
it follows that, 
for $z_1,z_2,z_3\in{\Bbb C}$, one has ${\cal F}(z_1,z_2,z_3)\neq 0$
only if 
$|z_1|^2\asymp T^{\varepsilon} |\varphi|^2 N^{-1}$, 
$|z_2|^2\asymp T^{\varepsilon} |\varphi|^2 \Xi^{-1}$ and 
$|z_3|^2\asymp T^{\varepsilon} |\Phi_1|^2 Z^{-1}$; 
since one can (moreover) verify that the function ${\cal F}(z_1,z_2,z_3)$  
defined by (6.70) is continuous on ${\Bbb C}^3$, it therefore 
follows from (6.69) that one has:  
$${\cal E}^{\flat}\ll \left| {\cal F}\!\left( z_1,z_2,z_3\right)\right|  
\ \,\hbox{for some}\ \,\left( z_1,z_2,z_3\right)\in ({\Bbb C}^{*})^3\ \,
\hbox{such that the conditions in (6.46) are satisfied.}\eqno(6.72)$$ 
\par 
We observe that, by (1.18), (1.24), and the hypotheses that 
$K_1,K_2,M_1\geq 1$ and $0<\varepsilon\leq 1/6$, it is certainly the case 
that $N\ll T$. Therefore, given that $\varphi$ satisfies (6.44), 
it follows from (6.66), (6.68) and (6.72) 
that, for some $(z_1,z_2,z_3)\in{\Bbb C}^3$ satisfying (6.46), one has 
$${\cal E}^{\star} 
\ll_{\varepsilon} \max\left\{ O_{\eta,\varepsilon}\left( 
T^{6\varepsilon} K_1 K_2^2 M_1 \| a\|_2^2\right) \,,\, 
|\varphi|^2 T^{\varepsilon} \left| {\cal F}\left( z_1,z_2,z_3\right)\right|\right\} . 
\eqno(6.73)$$ 
\par 
By (6.70), (6.71), (6.62) and the case $i=2$ of (1.19), we have  
$\,{\cal F}(z_1,z_2,z_3) = T | z_2 z_3|^2 Y$, where $Y$ is as defined 
in the equations (6.47) and (6.48). Recall moreover that, 
in all of our discussion since obtaining the estimate (6.51)-(6.52), 
we are able to assume a fixed choice of      
$\varphi,\varphi_1,\Phi_1,\varphi_2,\Phi_2\in{\frak O}$ such that   
both of the conditions (6.44) and (6.45) are satisfied. 
The results of the lemma are therefore an 
immediate corollary of our finding that (6.73) holds  
for some $(z_1,z_2,z_3)\in{\Bbb C}^3$ satisfying (6.46)\ $\square$  

\bigskip

\goodbreak 
\noindent{\SectionHeadingFont 7. Completing the proof of Theorem~2}

\medskip

\noindent In this section we first work to establish certain upper bounds 
for the absolute value of $Y$ (the sum of Kloosterman sums 
defined by (6.47)-(6.48), in Lemma~22).  
Our proof of these bounds is principally an exercise in 
the application of a theorem proved in [44]: for the 
convenience of the reader, this theorem [44, Theorem~11] 
is reproduced below (as Lemma~23). 
There are, however, certain extreme cases in which 
this method of obtaining bounds for $Y$ fails, 
because of a limitation on the scope of application of Lemma~23 that is 
an unfortunate artefact of our own construction, rather than being an essential 
feature of the method of proof used in [44] (for more details of this matter, 
see Remarks~24, below). In these extreme cases we obtain an adequate bound 
for $Y$ through the application of Lemma~25, below: although the bound in (7.12)   
is weaker than that in (7.8), it has the advantage 
of being valid under more general conditions.   
\par 
Lemma~26 and Lemma~27 are technical in nature (they enable 
us to ascertain that Lemma~23 and Lemma~25 are applicable, 
where needed). Bounds for the sum $Y$ are obtained in Lemma~28 
(it is in the proof of this lemma that Lemma~23 and Lemma~25 are applied). 
At the end of the section we complete the proof of Theorem~2 
by showing that the results (1.30)-(1.32) may be 
deduced from the bounds (7.33) and (7.34) of Lemma~28 and the 
three lemmas of Section~6. 
\par 
Our assumptions throughout this final section  
are those detailed in the second and third paragraphs of Section~6. 

\bigskip 

\goodbreak\proclaim{\Smallcaps Lemma~23}. Let $\vartheta$ be the real absolute constant 
defined in (1.10)-(1.11), let $\Psi_0 ,\Psi_1,\Psi_2,\Psi_3\geq 1\geq\delta >0$, 
and let $\varepsilon_1  >0$. Let $D$ be a complex valued function
with domain ${\frak O}-\{ 0\}$; for $h=1,2,3$, 
let $A_h : {\Bbb C}\rightarrow{\Bbb C}$ be a smooth function such that,  
for $j,k\in{\Bbb N}\cup\{ 0\}\ $ and $\ x,y\in{\Bbb R}$, one has 
$$(\delta |x+iy|)^{j+k}{\partial^{j+k}\over\partial x^j\partial y^k}\,A_h(x+iy) 
=\cases{ O_{j,k}(1) &if $\,\Psi_h /2<|x+iy|^2<\Psi_h$, \cr 0 &otherwise.}\eqno(7.1)$$
Let $P,Q,R,S\geq 1$ and $X>0$ satisfy 
$$Q=RS\geq\max\left\{ \sqrt{\Psi_0 }\,,\,\sqrt{\Psi_1}\right\}\eqno(7.2)$$
and 
$$X={PS\sqrt{R}\over 4\pi^2\sqrt{\Psi_0  \Psi_1}}\geq 2\;,\eqno(7.3)$$ 
and let $b$ be a complex-valued function with domain 
$${\cal B}(R,S)
=\left\{ (\rho,\sigma)\in{\frak O}\times{\frak O} : 
\,R/2<|\rho|^2\leq R\,,\ S/2<|\sigma|^2\leq S\ {\rm and}\ (\rho,\sigma)\sim 1\right\}\;.\eqno(7.4)$$
For each pair 
$(\rho,\sigma)\in{\cal B}(R,S)$, 
let $g_{\rho,\sigma} : (0,\infty)\rightarrow{\Bbb C}$ 
be an infinitely differentiable function such that, 
for $j\in{\Bbb N}\cup\{ 0\}$ and $x>0$,  one has 
$$g_{\rho,\sigma}^{(j)}(x)=\cases{O_j\left( (\eta x)^{-j}\right) &if $\,P/2<x<P$, \cr 
0 &otherwise.}\eqno(7.5)$$
Put 
$$Y_{\dagger}=\sum_{(\rho,\sigma)\in{\cal B}(R,S)} b(\rho,\sigma) 
\sum_{\Psi_0 /4<|\psi_0 |^2\leq \Psi_0 } D(\psi_0 ) \sum_{\Psi_1/2<|\psi_1 |^2 < \Psi_1} A_1(\psi_1 )  
\,G_{\rho,\sigma}(\psi_0 ,\psi_1 )\;,\eqno(7.6)$$
where 
$$G_{\rho,\sigma}(\psi_0 ,\psi_1 ) 
=\sum_{0\neq\varpi\in{\frak O}} g_{\rho,\sigma}\!\left( |\varpi|^2\right) 
S\!\left( \rho^{*}\psi_0  , \psi_1  ; \varpi\sigma\right)\;,\eqno(7.7)$$
with $S(\alpha,\beta;\gamma)$ being 
the Kloosterman sum that is defined in Equation~(5.5), 
and with $\rho^{*}\in{\frak O}$ satisfying 
$\rho^{*} \rho\equiv 1\bmod \varpi\sigma{\frak O}\,$    
(so that $\varpi$, the variable of summation in (7.7), is 
implicitly constrained to satisfy the condition $(\varpi,\rho)\sim 1$). 
Then 
$$Y_{\dagger}^2\ll_{\varepsilon_1,\eta} 
\,\left( 1+{X^2\over \left( 1+Q \Psi_0 ^{-1}\right) 
\!\left( 1+Q \Psi_1^{-1}\right)^{2}\!\Psi_1}\right)^{\!\!\vartheta} 
\!(\Psi_0 +Q)\!\left( \Psi_1 +Q\right)\delta^{-11} 
Q^{\varepsilon_1}\!\left\|b\right\|_2^2 \left\| D\right\|_2^2 \Psi_1 P^2 S\log^2(X)\; .\eqno(7.8)$$
\indent If it is moreover the case that 
$$\Psi_0 =\Psi_2 \Psi_3\;,\eqno(7.9)$$ 
and that 
$$D\left( \psi\right) 
=\sum_{\psi'\mid\psi} A_2\left(\psi'\right) 
A_3\!\left( {\psi\over\psi'}\right)\qquad\qquad 
\hbox{($0\neq\psi\in{\frak O}$),}\eqno(7.10)$$ 
then one has also:   
$$\eqalign{Y_{\dagger}^2
 &\ll_{\varepsilon_1,\eta} 
\,\left(\!\Biggl( 1+{X^2\over\left( \Psi_2 +\Psi_3\right) 
\left( 1+Q \Psi_1^{-1}\right)^{2} \Psi_1}\Biggr)^{\!\!\vartheta} \Psi_0   
+\Biggl( 1+{X^2\over Q^2 \Psi_0 ^{-1}\left( 1+Q \Psi_1^{-1}\right)^{2} \Psi_1}\Biggr)^{\!\!\vartheta} 
Q\!\right) \times \cr 
&\qquad\quad\times \left( \Psi_1 +Q\right)\delta^{-22} Q^{1+\varepsilon_1 }\left\| b\right\|_{\infty}^2 
\Psi_0  \Psi_1 P^2 S \log^2(X)\;.}\eqno(7.11)$$

\medskip

\goodbreak
\noindent{\Smallcaps Proof.}\  
This lemma is essentially just a reformulation of 
[44, Theorem~11]. It is however worth mentioning that, 
since $\eta$ may be arbitrarily small, the condition (7.5) 
is weaker than the corresponding condition in the statement of 
[44, Theorem~11] (i.e. in [44, Condition~(1.4.9)] 
one has $O_j(x^{-j})$, in place of the term $O_j((\eta x)^{-j})$ 
which occurs in (7.5)). To justify this weaker condition we observe that 
if $\eta$ is assigned any fixed numerical value from the interval 
$(0, (\log 2)/3]$ (that being the range of values permitted by our 
current assumptions) then $\eta$ effectively becomes an absolute positive 
constant, and so, subject to that assignment having been made, the 
condition (7.5) will in fact imply that one has  
$g_{\rho,\sigma}^{(j)}(x)\ll_j x^{-j}$ for all $x\in(P/2,P)$ and all $j\in{\Bbb N}\cup\{ 0\}$. 
We therefore obtain bounds which, except in 
respect of the dependence of the relevant implicit constants upon $\eta$,  
are equivalent to the bounds in [44, (1.4.12) and (1.4.14)]\quad$\square$ 

\bigskip 

\goodbreak\proclaim{\Smallcaps Remarks~24}. The condition $RS\geq\max\{\sqrt{\Psi_0 },\sqrt{\Psi_1}\}$ 
occurring in (7.2) is somewhat artificial:  it is not of any significance as regards 
the method of proof of [44, Theorem~11], but without it we would 
have had to include certain additional terms in the  
upper bounds (7.8) and (7.11) for $Y_{\dagger}^2$. 
At the time of writing [44] we had thought that the 
condition (7.2) (which corresponds to the condition [44, (1.4.6)]) 
would not hinder in any way our intended future use of 
Lemma~23 in proving Theorem~2 of the present paper. 
However we subsequently found that the constraint 
$RS\geq\max\{\sqrt{\Psi_0 },\sqrt{\Psi_1}\}$ does prevent the 
application of Lemma~23 in respect of certain `extreme' cases 
that arise in our proof of Theorem~2; we therefore regret our 
earlier decision to simplify [44, Theorem~11] through the 
inclusion of the condition `$RS\geq\max\{\sqrt{N},\sqrt{L}\}$' 
which appears in [44, (1.4.6)]. Rather than to show now how (7.8) or 
(7.11) should be modified when one has $RS<\max\{\sqrt{\Psi_0 },\sqrt{\Psi_1}\}$,  
we instead prefer an ad hoc solution to problem of the above mentioned 
deficiency of Lemma~23; our chosen solution makes use of the following lemma. 

\bigskip 

\goodbreak 
\proclaim{\Smallcaps Lemma~25}. Let $\vartheta$, $\varepsilon_1 $, $\delta$, $\Psi_0 $, 
$\Psi_1$, $A_1$, $D$, $P$, $Q$, $R$, $S$, $X$, 
$b$, ${\cal B}(R,S)$ and the family 
$(g_{\rho,\sigma})_{(\rho,\sigma)\in {\cal B}(R,S)}$ 
be such that, if one excludes the both the condition 
$RS\geq\max\{\sqrt{\Psi_0 },\sqrt{\Psi_1}\}$,  
occurring in (7.2), and those of the hypotheses of Lemma~23 that are concerned with  
`$\Psi_2$', or `$\Psi_3$', or `$A_2$', or `$A_3$', then all of the remaining hypotheses of 
that lemma are satisfied. 
Suppose moreover that $Y_{\dagger}$ is as defined in (7.6)-(7.7). 
Then 
$$Y_{\dagger}^2\ll_{\varepsilon_1 ,\eta} 
\,\left( \Psi_0  \Psi_1\right)^{\varepsilon_1 } \| b\|_2^2 \,\| D\|_2^2 \,\Psi_1 P^2 S 
\,(\Psi_0 +Q) \left( \Psi_1 +Q\right) X^{2\vartheta}\log^2 (X)\;.\eqno(7.12)$$ 

\medskip 

\goodbreak
\noindent{\Smallcaps Proof.}\  
By essentially the same steps as those through which the result 
[44, (9.75)-(9.76)] was arrived at (within the proof of 
[44, Theorem~11]), we find here that either 
$$Y_{\dagger}\ll_{\varepsilon_1 } \,\eta^{-5} \log(X) P S R^{1/2} 
\left\| D\right\|_2 
\left\| A_1\right\|_2 
\sum_{(\rho,\sigma)\in{\cal B}(R,S)}|b(\rho,\sigma)| 
\left( 1+{\Psi_0 ^{(1+\varepsilon_1 )/2}\over |\rho\sigma|}\right) 
\!\left( 1+{\Psi_1^{(1+\varepsilon_1 )/2}\over |\rho\sigma|}\right) ,\eqno(7.13)$$ 
or else  
$$Y_{\dagger}\ll \eta^{-2} \log(X) P S R^{1/2} \int\limits_{-\infty}^{\infty} 
{U(t)\,{\rm d}t\over (1+|t|)^2}\eqno(7.14)$$
where  
$$U(t)=\sum_{(\rho,\sigma)\in{\cal B}(R,S)}\!\!\!\!|b(\rho,\sigma)| 
\!\sum_{\scriptstyle V\atop\scriptstyle\nu_V>0}^{\left(\Gamma_0(\rho\sigma)\right)}
\!\!X^{\nu_V}   
\Biggl|\sum_{{\textstyle{\Psi_1\over 2}}<|\psi_1 |^2 < \Psi_1} 
\!\!\!\!{A_1(\psi_1 )\,c_V^{\infty}\!\left(\psi_1 ;\nu_V,0\right)\over 
\left| \psi_1\right|^{it}}  
\sum_{{\textstyle{\Psi_0 \over 4}}<|\psi_0 |^2\leq \Psi_0 } 
\!\!\!\!{\overline{D(\psi_0 )}\,c_V^{1/\sigma}\!\!\left( \psi_0 ;\nu_V,0\right)\over 
\left|\psi_0\right|^{-it}}\Biggr|\;,\eqno(7.15)$$ 
with the summation to which the superscript `$(\Gamma_0(\rho\sigma))$' attaches 
being over cuspidal irreducible subspaces 
$V$ of $L^2(\Gamma_0(\rho\sigma)\backslash SL(2,{\Bbb C}))$, while 
the `spectral parameter' $\nu_V$ and modified Fourier coefficients   
$c_V^{\infty}(\psi ;\nu_V,0)$, $c_V^{1/\sigma}(\psi ;\nu_V,0)\,$ 
($0\neq\psi\in{\frak O}$) are as defined in [44, Section~1.1]. We have, therefore, just 
two cases to consider: the case in which (7.13) holds and the case in which 
what is stated in (7.14)-(7.15) holds. 
\par 
If (7.13) holds then, by the definition (7.4), the Cauchy-Schwarz inequality and 
the case $h=1$, $j=k=0$ of the hypothesis (7.1), one has 
$$Y_{\dagger}^2\ll_{\varepsilon_1 , \eta} 
\log^2(X) P^2 S^2 R \,\| D\|_2^2 \,\Psi_1 \left( 1+{\Psi_0 ^{1+\varepsilon_1 }\over RS}\right) 
\left( 1+{\Psi_1^{1+\varepsilon_1 }\over RS}\right) 
RS \,\| b\|_2^2\;,$$ 
and so (given that we have $X\in[2,\infty)$, $RS=Q$ and, by (1.13), $\vartheta\in[0,2/9]$)  
the result (7.12) follows. 
\par 
We now have only to consider the 
case in which what is stated in (7.14)-(7.15) holds. 
By the discussion in [44, Section~1.1] leading up to the point noted 
in [44, (1.1.11)], it follows that in (7.15) one has $1-\nu_V^2 =\lambda_V$ (say), 
where $\lambda_V$ is some positive eigenvalue of the operator 
$-\Delta_3 : {\frak D}_{\Gamma_0(\rho\sigma)}\rightarrow L^2(\Gamma_0(\rho\sigma)\backslash{\Bbb H}_3)\,$ 
(the notation just used being that introduced between (1.8) and (1.9)). 
In particular, by (1.12) one has $1-\nu_V^2\geq 1-\vartheta^2$ in the sum occurring in (7.15), 
and so (given the explicit condition $\nu_V >0$ attached to that sum) 
each of the relevant spectral parameters $\nu_V$ must satisfy $\nu_V\leq\vartheta$. 
Therefore, since $X\geq 2$, and since  (by (1.13)) one has $\vartheta\leq 2/9<1$, 
it follows by way of the triangle inequality from (7.15) that, for all $t\in{\Bbb R}$, 
we have: 
$$\eqalign{|U(t)| &\leq X^{\vartheta}
\sum_{j=0}^1\ \sum_{(\rho,\sigma)\in{\cal B}(R,S)}\!\!\!\!|b(\rho,\sigma)|\ \times \cr 
 &\quad\times\sum_{\scriptstyle V\atop\scriptstyle 
\left| p_V\right| ,\left|\nu_V\right|\leq 1}^{\left(\Gamma_0(\rho\sigma)\right)}  
\Biggl|\sum_{{\textstyle{\Psi_1\over 2}}<|\psi_1 |^2 < \Psi_1}\!\!\!\!\!\!\!A_1(\psi_1 )  
|\psi_1 |^{-it}\,c_V^{\infty}\!\left(\psi_1 ;\nu_V,p_V\right) 
\!\!\sum_{{\textstyle{\Psi_0 \over 2^{j+1}}}<|\psi_0 |^2\leq\textstyle{\Psi_0 \over 2^j}} 
\!\!\!\!\!\!\!\overline{D(\psi_0 )} 
\,|\psi_0 |^{it} c_V^{1/\sigma}\!\left( \psi_0 ;\nu_V,p_V\right)\Biggr|\;.}$$
Hence, by an application of the Cauchy-Schwarz inequality, we are able 
to deduce that 
$$|U(t)|\leq X^{\vartheta}\sum_{j=0}^1\ \sum_{(\rho,\sigma)\in{\cal B}(R,S)} |b(\rho,\sigma)| 
\,U_{\rho,\sigma}^{1/2} V_{\rho,\sigma,j}^{1/2}\qquad\qquad\hbox{($t\in{\Bbb R}$),}\eqno(7.16)$$ 
where 
$$U_{\rho,\sigma}=U_{\rho,\sigma}(t)=\sum_{\scriptstyle V\atop\scriptstyle 
\left| p_V\right| ,\left|\nu_V\right|\leq 1}^{\left(\Gamma_0(\rho\sigma)\right)}  
\Biggl|\sum_{{\textstyle{\Psi_1\over 2}}<|\psi_1 |^2 < \Psi_1}\!\!\!A_1(\psi_1 )  
|\psi_1 |^{-it}\,c_V^{\infty}\!\left(\psi_1 ;\nu_V,p_V\right) \Biggr|^2\;,$$
$$V_{\rho,\sigma,j}=V_{\rho,\sigma,j}(t)=\sum_{\scriptstyle V\atop\scriptstyle 
\left| p_V\right| ,\left|\nu_V\right|\leq 1}^{\left(\Gamma_0(\rho\sigma)\right)}  
\Biggl|\sum_{{\textstyle{\Psi_0 \over 2^{j+1}}}<|\psi_0 |^2\leq\textstyle{\Psi_0 \over 2^j}}\!\!\overline{D(\psi_0 )} 
\,|\psi_0 |^{it} c_V^{1/\sigma}\!\left( \psi_0 ;\nu_V,p_V\right) \Biggr|^2\;.$$
By (7.4) and the case `$P=K=1$' of [45, Theorem~1] one has, in (7.16), 
$$U_{\rho,\sigma}\ll\left( 1+O_{\varepsilon_1 }\left( {\Psi_1^{1+\varepsilon_1 }\over |\rho\sigma|^2}\right)\right) 
\left\| A_1\right\|_2^2\qquad{\rm and}\qquad 
V_{\rho,\sigma,j}\ll\left( 1+O_{\varepsilon_1 }\left( {\Psi_0 ^{1+\varepsilon_1 }\over |\rho\sigma|^2}\right)\right) 
\left\| D\right\|_2^2\;.\eqno(7.17)$$ 
\par 
Since one has $\int_{-\infty}^{\infty}(1+|t|)^{-2}{\rm d}t=2<\infty$, 
it follows that from (7.16), (7.17) and the assumed upper bound (7.14) one may 
deduce the upper bound      
$$Y_{\dagger}\ll_{\varepsilon_1 , \eta} 
\,X^{\vartheta}\log(X) P S R^{1/2} 
\left\| D\right\|_2 
\left\| A_1\right\|_2 
\sum_{(\rho,\sigma)\in{\cal B}(R,S)}|b(\rho,\sigma)| 
\left( 1+{\Psi_0 ^{(1+\varepsilon_1 )/2}\over |\rho\sigma|}\right) 
\!\left( 1+{\Psi_1^{(1+\varepsilon_1 )/2}\over |\rho\sigma|}\right) .\eqno(7.18)$$  
By a calculation similar to that employed in dealing with the 
case in which (7.13) holds, it follows from (7.18) that 
we obtain the result (7.12)\quad$\square$ 

\bigskip 

\goodbreak\proclaim{\Smallcaps Lemma~26}. Let $\Upsilon_{\eta}$ be the function 
defined on the interval $[0,\infty)$ by:  
$$\Upsilon_{\eta}(u) = W_{\eta}(u)-W_{\eta}\!\left( e^{\eta} u\right)\qquad\qquad 
\hbox{($u\geq 0$).}\eqno(7.19)$$ 
Then $\Upsilon_{\eta}$ is real valued and infinitely differentiable on $[0,\infty)$; 
the support of $\Upsilon_{\eta}$ is contained in the interval $[e^{-\eta},e^{\eta}]$,  
and one has  both 
$$\Upsilon_{\eta}^{(j)}(u)\ll_j\,\eta^{-j}\qquad\qquad 
\hbox{($u\geq 0\,$ and $\,j=0,1,2,\ldots\ $)}\eqno(7.20)$$ 
and  
$$\sum_{0\leq h<1+\eta^{-1}\log B} 
\Upsilon_{\eta}\!\!\left( {|\beta|^2\over B e^{-h\eta}}\right) 
= W_{\eta}\!\!\left( {|\beta|^2\over B}\right) 
\qquad\qquad\hbox{($0\neq\beta\in{\frak O}\,$ and $\,0<B<\infty$).}\eqno(7.21)$$ 
Moroever, if $B_1\in(0,\infty)$ and $\gamma\in{\Bbb C}$, and if 
$A$ is the complex function given by 
$$A(z)=\Upsilon_{\eta}\!\!\left( {|z|^2\over B_1}\right) 
{\rm e}\left( \Re (\gamma z)\right)\qquad\qquad 
\hbox{($z\in{\Bbb C}$),}\eqno(7.22)$$ 
then $A$ is smooth on ${\Bbb C}$ and, for 
$$\delta ={1\over \eta^{-1} + B_1^{1/2} |\gamma|}\eqno(7.23)$$ 
and all $j,k\in{\Bbb N}\cup\{ 0\}$, one has: 
$$\left(\delta |x+iy|\right)^{j+k}\,{\partial^{j+k}\over\partial x^j \partial y^k} 
\,A(x+iy) 
=\cases{ O_{j,k}(1) &if $\,B_1 e^{-\eta}<|x+iy|^2<B_1 e^{\eta}$, \cr 
0 &otherwise,}\eqno(7.24)$$ 
at all points $(x,y)\in{\Bbb R}^2$.

\medskip 

\goodbreak
\noindent{\Smallcaps Proof.}\  
We omit the proofs of the results stated in, or above, (7.20): those results 
are straightforward consequences of our hypothesis that $W_{\eta}$ is as 
described in Remarks~18. In order to prove (7.21) we 
first note that, given (7.19), one can show by induction that 
$$\sum_{h=0}^{H-1} \Upsilon_{\eta}\!\left( e^{h\eta} u\right) 
=W_{\eta}(u)-W_{\eta}\!\left( e^{H\eta} u\right) \qquad\qquad 
\hbox{($u\geq 0\,$ and $\,H\in{\Bbb N}$).}\eqno(7.25)$$ 
Then we observe that if $B>0$ and $\beta\in{\frak O}-\{ 0\}$ then 
$|\beta|^2 /B\geq 1/B$, and so (given that 
$W_{\eta}(u)=\Upsilon_{\eta}(u)=0$ for $u\geq e^{\eta}$) 
one has 
$\Upsilon_{\eta}\!\left( e^{h\eta} |\beta|^2 /B\right) 
=W_{\eta}\!\left( e^{h\eta} |\beta|^2 /B\right) =0$ 
for all $h\in [1+\eta^{-1}\log B , \infty)$. The result (7.21) 
therefore follows from the case $H=1+[\eta^{-1}|\log B|]$, $u=|\beta|^2/B$ of 
the equality in (7.25). 
\par 
In our proofs of the remaining results of the lemma we may assume that 
$B_1 >0$, that $\gamma\in{\Bbb C}$, and that $A : {\Bbb C}\rightarrow{\Bbb C}$ 
is the function defined in (7.22). Since $\Upsilon_{\eta}$ is infinitely differentiable 
on $[0,\infty)$, and since the support of $\Upsilon_{\eta}$ is a compact subset of $(0,\infty)$, 
it therefore follows by (for example) the case $t=0$ of [44, Lemma~9.4] that 
the function $z\rightarrow\Upsilon_{\eta}(B_1^{-1} |z|^2)$ is smooth and compactly 
supported in ${\Bbb C}-\{ 0\}$. Moreover, for $x,y\in{\Bbb R}$, we have 
${\rm e}\left( \Re ((x+iy)\gamma))\right) 
=\exp\left( 2\pi i (ux-vy)\right) 
=\exp(2\pi i ux)\exp(-2\pi i vy)$, 
where $u=\Re (\gamma)$, $v=\Im (\gamma)$;  
note that, since $\exp(z)$ is a holomorphic function on ${\Bbb C}$, the 
mappings $z\mapsto\exp(2\pi iu\Re (z))$ and 
$z\mapsto\exp(-2\pi iv\Im (z))$ are smooth functions on ${\Bbb C}$.
Therefore, given that any product of two smooth functions is smooth, we may 
conclude that $A(z)$ (as defined in (7.22)) is a product of two 
functions that are smooth on ${\Bbb C}$, and so is itself smooth on 
${\Bbb C}$. 
\par 
Since the support of $\Upsilon_{\eta}$ is contained in the interval 
$[e^{-\eta},e^{\eta}]$, and since $\Upsilon_{\eta}$ is infinitely differentiable 
(and so continuous) on $[0,\infty)$, it follows that one has $\Upsilon_{\eta}(u)=0$ for 
all $u\in[0,e^{-\eta}]\cup[e^{\eta},\infty)$. Therefore, given the definition (7.22), 
we find that for all $z\in{\Bbb C}$ one has $A(z)\neq 0$ only if 
$B_1 e^{-\eta} < |z|^2 < B_1 e^{\eta}$. This last observation implies part of what 
is asserted in (7.23)-(7.24); in order to prove the remaining part of (7.23)-(7.24), we 
need only show that if $\delta$ is given by the equation (7.23) then one has the upper bounds  
$${\partial^{j+k}\over\partial x^j \partial y^k}\,A(x+iy) 
\ll_{j,k} \,(\delta |x+iy|)^{-(j+k)}\qquad\qquad\hbox{($j,k\in{\Bbb N}\cup\{ 0\}$)}$$ 
at all points $(x,y)\in{\Bbb R}^2$ satisfying 
$$B_1 e^{-\eta} < |x+iy|^2 <B_1 e^{\eta}\;.$$
Since $A$ is a smooth complex function, and since, when  
$x$ and $y$ are real variables and $z$ is the dependent complex variable $x+iy$, one has 
$${\partial\over\partial x}={\partial\over\partial z}+{\partial\over\partial\overline{z}}\qquad\quad 
{\rm and}\qquad\qquad 
{\partial\over\partial y}= i\left( 
{\partial\over\partial z}-{\partial\over\partial\overline{z}}\right) ,$$ 
where 
$${\partial\over\partial z}
={1\over 2}\left( {\partial\over\partial\Re(z)}-i{\partial\over\partial\Im(z)}\right)\qquad\quad 
{\rm and}\qquad\quad  
{\partial\over\partial\overline{z}} 
={1\over 2}\left( {\partial\over\partial\Re(z)}+i{\partial\over\partial\Im(z)}\right) ,$$    
it therefore follows from the preceding discussion that we may complete the proof of the 
lemma by showing that, for all $m,n\in{\Bbb N}\cup\{ 0\}$, one has 
$${\partial^{m+n}\over\partial z^m \partial\overline{z}^n}\,A(z) 
\ll_{m,n} \,(\delta |z|)^{-(m+n)}\qquad\qquad 
\hbox{($B_1 e^{-\eta} < |z|^2 < B_1 e^{\eta}$),}\eqno(7.26)$$
with $\delta$ as in (7.23) (and with $z$ understood to be the 
complex variable having real part $x$ and imaginary part $y$).  
Accordingly, we assume henceforth that 
$m$ and $n$ are non-negative integers, and that $\delta$ is the positive 
real number given by the equation (7.23). 
\par 
Given how the operators 
$\partial /\partial z$ and $\partial /\partial\overline{z}$ 
act upon holomorphic or anti-holomorphic functions 
(see, for example, [44, Equations~(5.21)]), 
it follows from the definition (7.22) (by way of a short calculation in which the 
identity ${\rm e}(\Re (\gamma z)) 
=\exp(\pi i \gamma z) \exp(\pi i\,\overline{\gamma}\,\overline{z})$ is utilized) 
that one has 
$${\partial^{m+n}\over\partial z^m \partial\overline{z}^n}\,A(z) 
=\sum_{r=0}^m \sum_{s=0}^n \pmatrix{m \cr r}\pmatrix{n \cr s} 
(\pi i\gamma)^r \left(\pi i\,\overline{\gamma}\right)^s 
{\rm e}\left( \Re (\gamma z)\right) 
\Biggl(\left( {\partial\over\partial z}\right)^{\!\!m-r} 
\!\!\left( {\partial\over\partial\overline{z}}\right)^{\!\!n-s}   
\!\Upsilon_{\eta}\!\!\left( {|z|^2\over B_1}\right)\Biggr) .\eqno(7.27)$$ 
Moreover,  by a direct application of the definitions of 
the operators $\partial /\partial z$, $\partial /\partial\overline{z}$,  
one can show that 
$$\left( {\partial\over\partial z}\right)^{\!\!m'} 
\!\!\left( {\partial\over\partial\overline{z}}\right)^{\!\!n'}   
\!\Upsilon_{\eta}\!\!\left( {|z|^2\over B_1}\right) 
=\sum_{0\leq t\leq\min\left\{ m' , n'\right\} } 
{1\over t!}\pmatrix{m' \cr t}\pmatrix{ n' \cr t} 
B_1^{t-\left( m'+n'\right) }\left(\overline{z}\right)^{m'-t} z^{n'-t}   
\Upsilon_{\eta}^{\left( m'+n'-t\right)}\!\!\left( B_1^{-1} |z|^2\right) , 
$$
whenever $m'$ and $n'$ are non-negative integers. By this last result, 
combined with (7.27) and the hypothesis (7.20), it follows that at all 
points $z\in{\Bbb C}$ such that $0<|z|^2<B_1 e^{\eta}$ one has: 
$$\eqalign{ 
{\partial^{m+n}\over\partial z^m \partial\overline{z}^n}\,A(z) 
 &\ll_{m,n}\,\sum_{r=0}^m \sum_{s=0}^n \pmatrix{m \cr r}\pmatrix{n \cr s} 
|\gamma|^{r+s} \sum_{0\leq t\leq\min\left\{ m-r , n-s\right\} } 
{|z|^{(m-r)+(n-s)-2t}\over\left(\eta B_1\right)^{(m-r)+(n-s)-t} } < \cr 
 &<  e^{(m+n)\eta} \sum_{r=0}^m \sum_{s=0}^n \pmatrix{m \cr r}\pmatrix{n \cr s} 
|\gamma|^{r+s} |z|^{(r+s)-(m+n)}\sum_{t=0}^{\infty}  
\eta^{r+s+t-(m+n)} = \cr 
 &= e^{(m+n)\eta} \left(\eta |z|\right)^{-(m+n)} 
\left( 1+ \eta |\gamma z|\right)^{m+n} (1-\eta)^{-1} \leq \cr 
 &\leq (1-\eta)^{-1} e^{(m+n)\eta} \left(\eta |z|\right)^{-(m+n)} 
\bigl( 1+ \eta |\gamma| B_1^{1/2} e^{\eta /2}\bigr)^{m+n} \leq \cr 
 &\leq  (1-\eta)^{-1} e^{(3/2)(m+n)\eta} \left(\delta |z|\right)^{-(m+n)}\;. 
}$$
Since $0<\eta\leq (\log 2)/3<1$, the above bounds 
imply that (7.26) holds, and so are sufficient to complete the proof of the lemma  
\quad$\square$ 

\bigskip 
 
\goodbreak\proclaim{\Smallcaps Lemma~27}. Let $\rho_1,\rho_2,\varphi_1,\varphi_2\in{\frak O}-\{ 0\}$, 
let ${\bf z}\in{\Bbb C}^3$, and let the function $g : (0,\infty)\rightarrow{\Bbb C}$ 
satisfy 
$$g(x)=g_{\rho_1,\rho_2}\!\left(\varphi_1,\varphi_2;z_1,z_2,z_3;x\right)\qquad\qquad 
\hbox{($0<x<\infty$),}$$ 
where $g_{\rho_1,\rho_2}\!\left(\varphi_1,\varphi_2;z_1,z_2,z_3;x\right)$ 
is defined by the equation (6.48) of Lemma~22. 
Then $g$ is infinitely differentiable on $(0,\infty)$, and, for 
$j\in{\Bbb N}\cup\{ 0\}$ and $x>0$, one has: 
$$g^{(j)}(x)=\cases{ O_j\!\left( (\eta x)^{-j}\right) &if $\,P/2<x<P$, \cr 
0 &otherwise,}\eqno(7.28)$$ 
where 
$$P=2^{1/2} \left|\varphi_2\right|^{-2} K_2\;.\eqno(7.29)$$ 

\medskip 

\goodbreak
\noindent{\Smallcaps Proof.}\  
We observe that one of the factors of the product on the right-hand side of equation (6.48) is 
the complex conjugate of $w_2(|\varphi_2|^2x)$. Therefore, 
by our hypothesis that (1.19) holds for $i=2$, it follows that,  
for $x\in(0,\infty)$, one has 
$g(x)\neq 0$ only if $e^{-\eta} |\varphi_2|^{-2} K_2\leq x\leq e^{\eta} |\varphi_2|^{-2} K_2$; 
this shows, since $1<e^{\eta}\leq 2^{1/3}$, that if 
$P$ is given by (7.29) then one has $g(x)=0$ for $x\in(0,2^{-6/7}P)\cup(2^{-1/7}P,\infty)$, and so 
one obtains the result (7.28) for all $x\in(0,P/2]\cup[P,\infty)$ and all $j\in{\Bbb N}\cup\{ 0\}$. 
\par 
Similarly, in light of the occurrence of the factor 
$w_0(T^{-1} |z_1 z_2 z_3 \rho_1|^2 |\rho_2|^4 x^2)$ on the right-hand side of 
equation (6.48), and given that we have $K_0=1$, we may infer from (1.19) 
that 
$$g(x)=0\qquad{\rm for\ all}\ \,x\in(0,\infty)\ \,{\rm satisfying}\quad\  
{\left| z_1 z_2 z_3 \rho_1\right|^2 \left|\rho_2\right|^4 x^2 
\over T}> 2^{1/2}\;.\eqno(7.30)$$ 
We postpone making use of this observation until a later stage of this proof. 
\par 
Since $w_0$, $w_1$, $w_2$ and the mapping $x\mapsto x^2$ are all 
infinitely differentiable complex valued functions on the interval $(0,\infty)$, 
it follows from the stated hypotheses of the lemma and the definition (6.48) 
that $g$ is an infinitely differentiable complex valued function
on the same interval. Indeed, by Leibniz's rule for higher order derivatives 
of products, one has, for $j=0,1,2,\ldots\ $, the equality  
$$\eqalign{ 
g^{(j)}(x) 
 &=\sum_{j_1\geq 0}\ \cdots\ \sum_{\scriptstyle j_6\geq 0\atop 
\!\!\!\!\!\!\!\!\!\!\!\!\!\!\!\!\!\!\!\!\!\!\!\!\!\!\!\!\!\!\! 
{\scriptstyle j_1+\ \cdots\ +j_6=j}} {j!\over\left( j_1 !\right)\cdots\left( j_6 !\right)} 
\left( {\left| z_1 z_2 z_3 \rho_1\right|^2 \left|\rho_2\right|^4  
\over T}{{\rm d}^{j_1}\over{\rm d} x^{j_1}}\,x^2\right) 
\left( {{\rm d}^{j_2}\over{\rm d} x^{j_2}}\,{w}_0\!\left({\left| z_1 z_2 z_3 
\rho_1\right|^2\left|\rho_2\right|^4 x^2\over T}\right)\right) \times \cr 
 &\quad\,\times\left( {{\rm d}^{j_3}\over{\rm d} x^{j_3}} 
\,w_1\!\left(\left| z_2\rho_2\right|^2 x\right)\right) 
\!\overline{\left( {{\rm d}^{j_4}\over{\rm d} x^{j_4}}
\,w_1\!\left(\left| z_2 z_3\rho_1\rho_2\right|^2 x\right)\right) }  
\!\left( {{\rm d}^{j_5}\over{\rm d} x^{j_5}} 
\,w_2\!\left(\left|\varphi_1 z_3\rho_2\right|^2 x\right)\right) 
\!\overline{\left( {{\rm d}^{j_6}\over{\rm d} x^{j_6}} 
\,w_2\!\left(\left|\varphi_2\right|^2 x\right)\right) } 
}\eqno(7.31)$$
at all points $x$ of the interval $(0,\infty)$. 
\par 
Given what has already been shown, the proof of the lemma will 
be complete if we can show next that $g^{(j)}(x)\ll_j (\eta x)^{-j}$  
whenever one has $j\in{\Bbb N}\cup\{ 0\}$ and $x\in(0,\infty)$. We begin this task by 
noting that it follows from the hypothesis (1.19)  
that, when $c$ is a non-negative real constant, when $i\in\{1,2\}$, 
and when $\,j\in{\Bbb N}\cup\{ 0\}$,  one has    
$${{\rm d}^j\over{\rm d}x^j}\,w_i(cx)\ll_j (\eta x)^{-j}\qquad\qquad 
\hbox{($x>0$).}\eqno(7.32)$$ 
Moreover, it can be proved by induction that  
for each $j\in{\Bbb N}\cup\{ 0\}$ one has the identity 
$${{\rm d}^j\over{\rm d}y^j}\,w_0\!\left( y^2\right) 
=\sum_{j/2\leq k\leq j}\alpha(j,k) w_0^{(k)}\!\left( y^2\right) y^{2k-j}\qquad\qquad 
\hbox{($y>0$),}$$ 
where $\alpha$ is the function on 
$({\Bbb N}\cup\{ 0\})\times({\Bbb N}\cup\{ 0\})$ 
determined by the recurrence relation  
$$\alpha(j,k)=\cases{ (2k+1-j)\alpha(j-1,k)+2\alpha(j-1,k-1) &if $\,0<j/2\leq k\leq j$, \cr 
1 &if $\,j=k=0$, \cr 
0 &otherwise.}$$
It therefore follows from the hypothesis (1.19) that, when $c$ is 
a non-negative real constant, and when $j\in{\Bbb N}\cup\{ 0\}$, one has 
$${{\rm d}^j\over{\rm d}x^j}\,w_0\!\left( cx^2\right) 
=\sum_{j/2\leq k\leq j}\alpha(j,k) w_0^{(k)}\!\left( c x^2\right) c^k x^{2k-j} 
=\sum_{j/2\leq k\leq j} O_j\!\left( \eta^{-k} x^{-j}\right) 
\ll_j (\eta x)^{-j}\qquad\qquad\hbox{($x>0$).}$$ 
By this last estimate, combined with (7.30), (7.31) and (7.32), 
we are able to deduce that 
$$g^{(j)}(x) 
=\sum_{j_1\geq 0}\ \cdots\ \sum_{\scriptstyle j_6\geq 0\atop 
\!\!\!\!\!\!\!\!\!\!\!\!\!\!\!\!\!\!\!\!\!\!\!\!\!\!\!\!\!\!\! 
{\scriptstyle j_1+\ \cdots\ +j_6=j}} 
O_{j}\!\left( x^{-j_1}\right)\prod_{\ell =2}^6 
O_{j_{\ell}}\left( (\eta x)^{-j_{\ell}}\right) 
\ll_j \,(\eta x)^{-j} \max\{ \eta^{j_1} : j_1\geq 0\} =(\eta x)^{-j}$$
whenever $j\in{\Bbb N}\cup\{ 0\}$ and $x>0$\quad$\square$ 

\bigskip 

\goodbreak\proclaim{\Smallcaps Lemma~28}. Let $\Xi$ and $Z$ be as defined in (6.2) and (6.43), respectively. 
Let $\varphi$, $\Phi_1,\varphi_1,\Phi_2,\varphi_2\in{\frak O}$ and 
${\bf z}\in{\Bbb C}^3$ be such that the conditions (6.44)-(6.46) of Lemma~22 
are satisfied, and let $Y$ be as defined in (6.47)-(6.48). 
Then one has both 
$$Y\ll_{\eta,\varepsilon} 
\left( {T^{8\varepsilon} K_2^2 M_1^2 \| a\|_2^2\over 
|\varphi|^{4} \left|\Phi_1\right|^{2}}\right)  
\!\left( {M_1^2\over T^{1/2}} + 
\left( {M_1^{2-(3/2)\vartheta}\over T^{(1-\vartheta)/2}}\right) 
\left( {K_2\over T^{1/2}}\right)^{\!\!\vartheta}\right)\eqno(7.33)$$ 
and 
$$Y\ll_{\eta,\varepsilon} 
\left( {T^{10\varepsilon} K_2^2 M_1^3 \| a\|_{\infty}^2\over 
|\varphi|^{4} \left|\Phi_1\right|^{2}}\right)  
\!\left( {M_1^2\over T^{1/2}} + 
\left( {M_1^2\over T^{1/2}}\right)^{\!\!1-\vartheta}
\!\left(\left( {K_1\over T^{1/2}}\right)^{\!\!\vartheta /2} 
\left( {K_2\over T^{1/2}}\right)^{\!\!1-(\vartheta /2)} + 
\left( {K_2\over T^{1/2} M_1^{1/2}}\right)^{\!\!\vartheta} 
\right)\right) ,\eqno(7.34)$$ 
where $\vartheta$ is the real absolute constant defined in (1.10)-(1.11). 

\medskip 

\goodbreak
\noindent{\Smallcaps Proof.}\  
Given the definition (6.47), it follows by virtue of the identity (7.21) of 
Lemma~26 that we have 
$$Y=\sum_{\scriptstyle {\bf h}\in{\Bbb Z}^3\atop\scriptstyle 
0\leq h_i<H_i\ \,(i=1,2,3)} Y_{\eta}\left( h_1,h_2,h_3\right)\;,\eqno(7.35)$$ 
where:  
$$H_1=1+\eta^{-1}\log\left(\left|\Phi_1\right|^{-2} Z\right) , \qquad 
H_2=1+\eta^{-1}\log\left( |\varphi|^{-2} \Xi\right) , \qquad 
H_3=1+\eta^{-1}\log\left( |\varphi|^{-2} N\right)\eqno(7.36)$$
and 
$$\eqalign{Y_{\eta}\left( h_1,h_2,h_3\right) 
 &=\sum_{\scriptstyle\rho_1\atop\scriptstyle\left(\rho_1,\varphi_1\right)\sim 1} 
\!\!\!\!\!\sum_{\,{\scriptstyle\rho_2\atop\scriptstyle 
\ \sim\left(\rho_2,\varphi_2\right)}} 
\!\!\!\!\!\!a\!\left(\Phi_1\rho_1\right) \overline{a\!\left(\Phi_2\rho_2\right)} 
\,\sum_{\nu} \sum_{\xi} \sum_{\zeta} 
\Upsilon_{\eta}\!\!\left({|\varphi\nu|^2\over e^{-h_3\eta} N}\right)   
\Upsilon_{\eta}\!\!\left({|\varphi\xi|^2\over e^{-h_2\eta} \Xi}\right) \,\times \cr 
 &\quad\ \times 
\Upsilon_{\eta}\!\!\left({|\Phi_1\zeta|^2\over e^{-h_1\eta} Z}\right) 
{\rm e}\left(\Re\bigl( {\bf z}\cdot(\nu,\xi,\zeta)\bigr)\right) 
\sum_{\varpi\neq 0} g_{\rho_1,\rho_2}\!\!\left(\left|\varpi\right|^2\right) 
S\!\left(\nu\xi\rho_1^{*} , \zeta ; \rho_2\varpi\right)\;, 
}\eqno(7.37)$$
with $S(\alpha,\beta;\gamma)$, $g_{\rho_1,\rho_2}(x)$ and $\Upsilon_{\eta}(u)$ 
as given by (5.5), (6.48) and (7.19), respectively. 
By (6.43)-(6.45), (6.1), (6.2), (1.18), (1.24) and the hypotheses of 
Theorem~2 concerning $\varepsilon$, $\eta$ and $K_1$, it follows from the 
definitions in (7.36) that the parameters $H_1$, $H_2$ and $H_3$ occurring 
in (7.35) satisfy 
$$1+H_i \leq O\left( \eta^{-1}\log(T)\right)\qquad\qquad 
\hbox{($i=1,2,3$).}\eqno(7.38)$$ 
Therefore either it is the case that the sum on the right-hand side of Equation~(7.35) is empty 
(i.e. void of terms), 
or else it is the case that, for some ${\bf h}\in{\Bbb Z}^3$ satisfying 
$$0\leq h_i<H_i\qquad\qquad\hbox{($i=1,2,3$),}\eqno(7.39)$$ 
one has 
$$Y\ll \left(\eta^{-1}\log(T)\right)^3 \left| Y_{\eta}\left( h_1,h_2,h_3\right)\right|\;. 
\eqno(7.40)$$ 
Only the latter of these two cases need concern us in what follows: 
for in the former case one has $Y=0$, so that 
the results (7.33) and (7.34) of the lemma are (in that case) certainly valid. 
Therefore we shall assume henceforth that 
${\bf h}=(h_1,h_2,h_3)$ is a given element of ${\Bbb Z}^3$ satisfying the 
conditions in (7.39), and that this ${\bf h}$ is  such 
that the bound (7.40) holds. 
\par 
We shall complete this proof with the aid of 
certain bounds for the sum $Y_{\eta}(h_1,h_2,h_3)$ that we have defined in (7.37);     
the bounds in question will be obtained 
through applications of Lemma~23 and Lemma~25.  
Having that in mind, we now put  
$$\Psi_1=e^{(1-h_1)\eta} \left|\Phi_1\right|^{-2} Z\;,\qquad 
\Psi_2=e^{(1-h_2)\eta} |\varphi|^{-2}\Xi\;,\qquad 
\Psi_3=e^{(1-h_3)\eta} |\varphi|^{-2} N\;,\eqno(7.41)$$ 
$$A_t(s)=\Upsilon_{\eta}\!\left( {|s|^2\over e^{-\eta} \Psi_t}\right) 
{\rm e}\left( \Re \left( s z_{4-t}\right)\right)\qquad\qquad 
\hbox{($s\in{\Bbb C}$, $\,t=1,2,3$),}\eqno(7.42)$$
$$\Psi_0 =\Psi_2 \Psi_3\;,\qquad\qquad D(\psi) =\sum_{\psi'\mid\psi} 
A_2(\psi') A_3\!\left( {\psi\over\psi'}\right)\qquad  
\hbox{($0\neq\psi\in{\frak O}$),}\eqno(7.43)$$ 
$$R=e^{\eta} \left|\Phi_1\right|^{-2} M_1\qquad{\rm and}\qquad 
S=e^{\eta} \left|\Phi_2\right|^{-2} M_1\;,\eqno(7.44)$$ 
and we take ${\cal B}(R,S)\subset{\frak O}\times{\frak O}$ to be defined by (7.4) and (7.44); 
for $(\rho_1,\rho_2)=(\rho,\sigma)\in{\cal B}(R,S)$, we put 
$$b(\rho,\sigma)=\cases{ a\left(\Phi_1\rho\right)\,\overline{a\left(\Phi_2\sigma\right)} 
 &if $\left(\rho_1 , \varphi_1\right)\sim\left(\sigma , \varphi_2\right)\sim 1$, \cr 
0 &otherwise,}\eqno(7.45)$$  
and we take 
$g_{\rho_1,\rho_2}$ to be (as in (7.37)) the function defined on the interval  
$(0,\infty)$ by the equation (6.48). 
\par 
Let $Y_{\dagger}$ be 
the sum of complex terms defined by the equations (7.6), (7.7) and (7.4) of Lemma~23,   
in conjunction with (7.41)-(7.45), (6.48) and (5.5).   
By (7.37), the hypothesis (1.20), the result (7.24) of Lemma~26 and 
the hypothesis that $0<\eta\leq(\log 2)/3$, it follows that we have 
$$Y_{\dagger}=Y_{\eta}\left( h_1,h_2,h_3\right)\;.\eqno(7.46)$$  
\par 
We seek next to verify as many as possible of the  
hypotheses of Lemma~23, assuming that 
$\Psi_0$, $\Psi_1$, $\Psi_2$, $\Psi_3$, $R$, $S$, the set ${\cal B}(R,S)$, 
the functions $A_1$, $A_2$, $A_3$, $D$, $b$ and $g_{\rho,\sigma}\,$ 
($(\rho,\sigma)\in{\cal B}(R,S)$) and the sum $Y_{\dagger}$ 
are as we have indicated in the last two paragraphs above 
(and also bearing in mind our hypotheses 
throughout the present section, and the additional hypotheses of the lemma).  
\par 
By (6.46), (7.39) and (7.41), we have 
$e^{-\eta} \Psi_t |z_{4-t}|^2\ll e^{-h_t \eta} T^{\varepsilon} \leq T^{\varepsilon}$, 
for $t=1,2,3$.  Therefore it follows by the results (7.22)-(7.24) of 
Lemma~26 that the complex functions $A_1$, $A_2$ and $A_3$ defined in (7.42) 
are smooth on ${\Bbb C}$, and are such that, 
for $t=1,2,3$ and all $j,k\in{\Bbb N}\cup\{ 0\}$, one has
$$\left( {|x+iy|\over\eta^{-1} + T^{\varepsilon /2}}\right)^{\!\!j+k} 
\!\!{\partial^{j+k}\over\partial x^j\partial y^k}\,A_t(x+iy) 
=\cases{ O_{j,k}(1) &if $\,e^{-2\eta} \Psi_t <|x+iy|^2<\Psi_t$, \cr 0 &otherwise,}\eqno(7.47)$$  
at all points $(x,y)\in{\Bbb R}^2$.
Independently of the point just noted, it is shown by Lemma~27 that, 
for $(\rho_1,\rho_2)=(\rho,\sigma)\in({\frak O}-\{ 0\})^2$, the 
function $g_{\rho_1,\rho_2} : (0,\infty)\rightarrow{\Bbb C}$ defined by 
Equation~(6.48) is infinitely differentiable on $(0,\infty)$, and is furthermore  
such that, for all $j\in{\Bbb N}\cup\{ 0\}$ and all $x>0$, one has 
$$g_{\rho_1,\rho_2}^{(j)}(x) 
=\cases{O_j\left( (\eta x)^{-j}\right) &if $\,2^{-1/2}\left| \varphi_2\right|^{-2} K_2 
<x<2^{1/2}\left| \varphi_2\right|^{-2} K_2$, \cr 
0 &otherwise.}\eqno(7.48)$$
By (7.36), (7.39), (7.41) and (7.43), we have $\Psi_1,\Psi_2,\Psi_3\geq 1$ and $\Psi_0 \geq 1$. 
Since the definition (7.4) implies that ${\cal B}(R,S)$ is the empty set if 
either $R<1$ or $S<1$, and since the definition (7.6) trivially implies 
that $Y_{\dagger}$ equals zero if ${\cal B}(R,S)=\emptyset$, 
we may therefore assume it is also the case that   
$$R,S\geq 1.\eqno(7.49)$$ 
\par 
We put now 
$$\delta =\bigl( \eta^{-1} + T^{\varepsilon /2}\bigr)^{-1}\;,\qquad\quad 
\varepsilon_1 =\varepsilon\;,\qquad\quad 
P=2^{1/2} \left|\varphi_2\right|^{-2} K_2 \eqno(7.50)$$ 
and (as hypothesised in Lemma~23) 
$$Q=RS\qquad{\rm and}\qquad X={PS\sqrt{R}\over 4\pi^2\sqrt{\Psi_0  \Psi_1}}\;.\eqno(7.51)$$ 
Given that $0<\varepsilon\leq 1/6$, it follows by (7.51), (7.50), (7.43), (7.44),  
(7.41), (6.45), (1.24), (6.2), (6.43), (7.39) and (6.1) that we have 
$$\eqalign{ 
4\pi^2 X &={2^{1/2} \left|\varphi_2\right|^{-2} K_2 S \sqrt{R}\over  
\sqrt{ \Psi_1 \Psi_2 \Psi_3}} = \cr 
 &={2^{1/2}  \left|\varphi_2 \Phi_2\right|^{-2} 
K_2 M_1^{3/2}\over\sqrt{e^{-(h_1+h_2+h_3)\eta} |\varphi|^{-4} Z \Xi N}}  
=2^{1/2} e^{(h_1+h_2+h_3)\eta/2} T^{(1-3\varepsilon)/2} 
\geq 2^{1/2} T^{1/4} \geq 8 \pi^2\;.}\eqno(7.52)$$
\par 
Since $0<\eta\leq (\log 2)/3<(\log 2)/2<1\,$, 
it follows from (7.47), (7.48), (7.49), (7.52) 
and the observations accompanying those points that, 
if one excludes the condition $RS\geq\max\{ \sqrt{\Psi_0 } , \sqrt{\Psi_1}\}$ 
occurring in (7.2), then the remaining hypotheses of Lemma~23  
(including the conditions (7.9) and (7.10) 
attached to the result (7.11))
are satisfied by  
the choice of $D$, $A_1$, $A_2$, $A_3$, ${\cal B}(R,S)$, $b$, 
$\delta$, $\varepsilon_1$, $g_{\rho,\sigma}\,$ ($(\rho,\sigma)\in{\cal B}(R,S)$),  
$\Psi_0 $, $\Psi_1$, $\Psi_2$, $\Psi_3$, $P$, $Q$, $R$, $S$ and 
$Y_{\dagger}$ declared in (7.41)-(7.45), between (7.45) and (7.46), 
and in (7.50) and (7.51).  
It therefore follows by Lemma~25 that we have the bound 
(7.12) for $Y_{\dagger}^2$. Similarly, it follows by Lemma~23 that 
either we have $\max\{\sqrt{\Psi_0 },\sqrt{\Psi_1}\} >RS$, or 
else it is the case that the bounds (7.8) and (7.11) for $Y_{\dagger}^2$ are valid. 
We conclude from this that there are only two cases requiring further consideration: 
one of these being the case in which both of the bounds (7.8) and (7.11) hold; 
the other being the case in which both 
the bound (7.12) and the inequality $\max\{\sqrt{\Psi_0 },\sqrt{\Psi_1}\} >RS$ hold. 
\par 
In order to facilitate our use of the results 
(7.8), (7.11) and (7.12) of Lemma~23 and Lemma~25, 
we note here that, by   
(7.49), (7.44), (7.51), (7.52), (1.18), (7.38), (7.39), (7.41), (6.43), (1.24), (7.43), 
(6.2) and (6.45), it follows that  
$$\left|\Phi_i\right|^2\leq e^{\eta} M_1 < 2 M_1\qquad\qquad\hbox{($i=1,2$),}\eqno(7.53)$$ 
$$Q\asymp\left|\Phi_1\Phi_2\right|^{-2} M_1^2\ll T\;,\qquad\qquad\quad  
2\leq X\ll e^{(h_1+h_2+h_3)\eta /2} T^{(1-3\varepsilon)/2} =T^{O(1)}\;,\eqno(7.54)$$ 
$$\Psi_1\ll e^{-h_1\eta} \left|\Phi_1\right|^{-2} T^{\varepsilon} M_1 
< T^{\varepsilon} R\leq T^{\varepsilon} Q\;, 
\qquad\qquad  
\Psi_3\ll e^{-h_3\eta} \left|\varphi\right|^{-2} T^{\varepsilon -1} K_1 K_2 M_1\eqno(7.55)$$   
and 
$$\Psi_0 =\Psi_2 \Psi_3\ll e^{-(h_2+h_3)\eta} \left|\varphi\right|^{-4} T^{2\varepsilon -1} K_2^2 M_1^2 
\ll\left|\varphi\right|^{-4} T^{2\varepsilon} M_1^2\ll T^{2\varepsilon} Q\ll T^2\;. 
\eqno(7.56)$$ 
By (6.44), (1.24) and (1.18), we also have   
$$\left|\varphi\right|^2\ll T^{\varepsilon -1} K_1 K_2 M_1 
\ll T^{\varepsilon -(1/2)} K_2 M_1\ll T^{\varepsilon} M_1\;. 
\eqno(7.57)$$ 
It moreover follows from (7.45), (7.43), (7.56) and the case $j=k=0$ of (7.47) that 
we have: 
$$\| b\|_2^2\leq\sum_{(\rho,\sigma)\in{\frak O}\times{\frak O}} 
\left| a\left(\Phi_1\rho\right) a\left(\Phi_2\sigma\right)\right|^2 
\leq \| a\|_2^4,\qquad\qquad 
\| b\|_{\infty}\leq \| a\|_{\infty}^2\eqno(7.58)$$ 
and 
$$\eqalign{ 
\| D\|_2^2 &=\sum_{\psi_0 \neq 0}\Biggl|\,\sum_{|\psi_2|^2\asymp \Psi_2} 
\sum_{\scriptstyle \,|\psi_3|^2\asymp \Psi_3\atop\scriptstyle 
\!\!\!\!\!\!\!\!\!\!\!\!\!\!\!\!\!\!\!\!\!\!\!\psi_2\psi_3 =\psi_0 } O(1)\Biggr|^2 \ll \cr 
 &\ll\ \,\sum\quad\ \,\sum_{\!\!\!\!\!\!\!\!\!\!\!\!\!\!\!\!\!\!\!\!\!\!\! 
|\psi_2|^2\asymp \Psi_2\asymp |\psi_2'|^2}
\ \,\sum\quad\ \,\sum_{\scriptstyle\!\!\!\!\!\!\!\!\!\!\!\!\!\!\!\!\!\!\!\!\!\!\! 
|\psi_3|^2\asymp \Psi_3\asymp |\psi_3'|^2\atop\scriptstyle 
\!\!\!\!\!\!\!\!\!\!\!\!\!\!\!\!\!\!\!\!\!\!\!\!\!\!\!\!\!\!
\!\!\!\!\!\!\!\!\!\!\!\!\!\!\!\!\!\!\!\!\!\!\!\!\!\!\!\!\!\!\!\!\!\!\!\!\!\!\!\! 
\psi_2'/\psi_2 =\psi_3 /\psi_3'} 1 
\ll \Psi_2 \Psi_3\left( 1 + \log\left(\min\left\{ \Psi_2 , \Psi_3\right\}\right)\right) 
\ll_{\varepsilon}  \Psi_0  T^{\varepsilon}\;.}\eqno(7.59)$$ 
By (7.44), (7.50), (7.54), (7.55), (7.56) and (6.45), we find furthermore that 
$$\Psi_0  \Psi_1 P^2 S Q^2 \ll e^{-(h_1+h_2+h_3)\eta} |\varphi|^{-8} 
\left|\Phi_1\right|^{-6} \left|\Phi_2\right|^{-2} T^{3\varepsilon -1} 
K_2^4 M_1^8\;.\eqno(7.60)$$ 
\par 
If both of the bounds (7.8) and (7.11) hold then,  
by (7.54), (7.55), (7.56) and the first two parts of (7.50), 
in combination with (7.58), (7.59), (7.60), (7.39), (6.45), (6.1)   
and (1.13), it follows that one has both 
$$\eqalignno{ 
Y_{\dagger}^2 &\ll_{\varepsilon ,\eta} 
Q^{\varepsilon} \| b\|_2^2 \| D\|_2^2 \Psi_1 P^2 S \log^2(X) 
\left( 1+ {X^2 \Psi_0  \Psi_1\over Q^3}\right)^{\!\!\vartheta} Q^2 T^{3\varepsilon} 
\left(\eta^{-1} + T^{\varepsilon /2}\right)^{\!\!11} 
\ll_{\varepsilon , \eta} \cr 
 &\ll_{\varepsilon , \eta} 
T^{10\varepsilon} \| b\|_2^2 \| D\|_2^2 \Psi_1 P^2 S Q^2 
\left( 1+ {T^{1-3\varepsilon} |\varphi|^{-4} T^{2\varepsilon -1} K_2^2 M_1^2 
\left|\Phi_1\right|^{-2} T^{\varepsilon} M_1\over \left|\Phi_1\Phi_2\right|^{-6} M_1^6} 
\right)^{\!\!\vartheta} \ll_{\varepsilon} \cr 
 &\ll_{\varepsilon} 
T^{11\varepsilon} \| a\|_2^4 \,\Psi_0  \Psi_1 P^2 S Q^2 
\left( 1+ {\left|\Phi_1\Phi_2\right|^3 K_2^2\over M_1^3}\right)^{\!\!\vartheta} \ll \cr 
 &\ll T^{14\varepsilon -1} \| a\|_2^4 \,|\varphi|^{-8} 
\left|\Phi_1\right|^{-4} K_2^4 M_1^8  
\left( 1+ {K_2^2\over M_1^3}\right)^{\!\!\vartheta} &(7.61)}$$ 
and 
$$\eqalignno{ 
Y_{\dagger}^2 &\ll_{\varepsilon ,\eta} 
Q^{1+\varepsilon} \| b\|_{\infty}^2 \Psi_0  \Psi_1 P^2 S \log^2(X) 
\left(\left( 1 + {X^2 \Psi_1\over Q^2 \Psi_2}\right)^{\!\!\vartheta} \!\Psi_0  
+\left( 1 + {X^2 \Psi_0  \Psi_1\over Q^4}\right)^{\!\!\vartheta} \!Q\right) 
Q T^{\varepsilon} 
\left(\eta^{-1} + T^{\varepsilon /2}\right)^{\!\!22} 
\ll_{\varepsilon , \eta} \cr 
 &\ll_{\varepsilon , \eta} 
T^{14\varepsilon} \| b\|_{\infty}^2 \Psi_0  \Psi_1 P^2 S Q^2 
\left(\left( {X^2 \Psi_0  \Psi_1 \Psi_3\over Q^2}\right)^{\!\!\vartheta} \!\Psi_0 ^{1-2\vartheta}  
+\left( 1 + {X^2 \Psi_0  \Psi_1\over Q^4}\right)^{\!\!\vartheta} \!Q T^{2\varepsilon}\right) \ll_{\varepsilon} \cr 
 &\ll_{\varepsilon} 
T^{14\varepsilon} \| a\|_{\infty}^4 |\varphi|^{-8} 
\left|\Phi_1\right|^{-6} \left|\Phi_2\right|^{-2} T^{3\varepsilon -1} 
K_2^4 M_1^8 \ \times \cr 
&\quad\ \ \,\times\left(\left( { \left|\Phi_1\right|^2 \left|\Phi_2\right|^4 
T^{\varepsilon -1} K_1 K_2^3
\over |\varphi|^6}\right)^{\!\!\vartheta} 
\!\left( T^{2\varepsilon -1} K_2^2 M_1^2\right)^{1-2\vartheta}   
+\left( 1 + {\left|\Phi_1\right|^6 \left|\Phi_2\right|^8 K_2^2\over 
|\varphi|^4 M_1^5}\right)^{\!\!\vartheta} \!\left|\Phi_1\Phi_2\right|^{-2} M_1^2 T^{2\varepsilon}\right) \leq \cr 
 &\leq T^{19\varepsilon -1} \| a\|_{\infty}^4 |\varphi|^{-8} 
\left|\Phi_1\right|^{-4} K_2^4 M_1^8 
\left( T^{\vartheta -1} K_1^{\vartheta} K_2^{2-\vartheta} M_1^{2-4\vartheta}     
+\left( 1 + {K_2^2\over M_1^5}\right)^{\!\!\vartheta} 
\!M_1^2 \right)\;. &(7.62)}$$ 
\par 
By (7.40) and (7.46), the bounds (7.61) and (7.62) imply the bounds 
for $Y$ that are stated in (7.33) and (7.34). Therefore, given the 
the conclusions of the paragraph immediately below (7.52), and bearing in mind 
the conditions 
subject to which (7.61) and (7.62) were obtained, it 
will suffice for the completion of the proof of the lemma 
that we show that, 
if it is the case that 
both the inequality $\max\{\sqrt{\Psi_0 },\sqrt{\Psi_1}\} >RS$ and 
the bound (7.12) hold, 
then it must also be the case that the bounds (7.33) and (7.34) hold.
Accordingly, we assume henceforth that the bound (7.12) holds, and 
that we have 
$$\max\left\{ \sqrt{\Psi_0 } , \sqrt{ \Psi_1}\right\} > RS\;.\eqno(7.63)$$ 
\par 
By (7.63), (7.51), (7.54), (7.55), (7.56) and (7.39), we have  
$\max\{ |\varphi|^{-2} T^{\varepsilon} M_1 ,  
|\Phi_1|^{-1} T^{\varepsilon /2} M_1^{1/2}\} \gg 
|\Phi_1\Phi_2|^{-2} M_1^2$. It therefore follows, 
given (6.1), (6.45) and (7.53), that either 
$$M_1\ll\left( {\left|\Phi_1\Phi_2\right|^2\over |\varphi|^2}\right) T^{\varepsilon} 
\leq \left( {\left|\Phi_1\Phi_2\right|^2\over 
\max\left\{\left|\Phi_1\right|^2\,,\,\left|\Phi_2\right|^2\right\}}\right) 
T^{\varepsilon} \;,\eqno(7.64)$$ 
or else 
$$M_1^{3/2}\ll\left|\Phi_1\right|\left|\Phi_2\right|^2 T^{\varepsilon /2} 
\ll \min\left\{ M_1^{1/2} \left|\Phi_2\right|^2 T^{\varepsilon}\,,\, 
\left|\Phi_1\right| M_1 T^{\varepsilon /2}\right\}\;.\eqno(7.65)$$ 
If it is the case that (7.64) holds, then it follows that one has 
$M_1 T^{-\varepsilon}\ll\min\{ |\Phi_1|^2 , |\Phi_2|^2\}$; the same conclusion 
follows if it is instead (7.65) which holds. Hence, and by (6.45), we 
are certain to have: 
$$T^{-\varepsilon} M_1\ll \left|\Phi_i\right|^2\leq |\varphi|^2\;,\qquad\quad 
\hbox{for $\,i=1,2$.}\eqno(7.66)$$ 
In view of (7.57), we find furthermore that the bounds in (7.66) imply that 
$$K_2\gg T^{(1/2)-2\varepsilon}\;.\eqno(7.67)$$ 
\par 
We remark that, in view of the bounds for $|\varphi|^2$ in (7.57), 
our conclusion in (7.66) shows that the case that we are now  
considering is extreme, in the sense that $\varphi$, which first occurs  
(within the proof of Lemma~20)  as 
a common factor of two independent variables, is almost as large in modulus as 
it can possibly be. 
\par 
Given that the bound (7.12) holds, and that we have also $\varepsilon_1 =\varepsilon$ 
and what is stated in (6.1), (7.54)-(7.56) and (7.58)-(7.60), it follows that 
$$\eqalignno{ 
Y_{\dagger}^2 &\ll_{\varepsilon,\eta}  
\,\left( T^{3\varepsilon} Q^2\right)^{\varepsilon} 
\| b\|_2^2 \| D\|_2^2 \Psi_1 P^2 S Q^2 T^{3\varepsilon} X^{2\vartheta}\log^2(T) \ll_{\varepsilon}\cr 
 &\ll_{\varepsilon} T^{7\varepsilon} \| a\|_2^4 \Psi_0  \Psi_1 P^2 S Q^2 X^{2\vartheta} \ll \cr 
 &\ll T^{7\varepsilon} \| a\|_2^4 \,e^{-(h_1+h_2+h_3)\eta} \,|\varphi|^{-8} \left|\Phi_1\right|^{-6} 
\left|\Phi_2\right|^{-2} T^{3\varepsilon -1} K_2^4 M_1^8
\left( T e^{(h_1+h_2+h_3)\eta}\right)^{\!\vartheta} \leq \cr 
 &\leq |\varphi|^{-8} \left|\Phi_1\right|^{-6} 
\left|\Phi_2\right|^{-2} T^{10\varepsilon +\vartheta -1} K_2^4 M_1^8 \| a\|_2^4\;.
 &(7.68)}$$  
By (7.68), (7.67), (7.66), (6.1) and the bound $\vartheta\leq 2/9$ of Kim and Shahidi, we obtain: 
$$\eqalign{ 
Y_{\dagger}^2 &\ll_{\varepsilon,\eta}  
\,|\varphi|^{-8} \left|\Phi_1\right|^{-4} 
\left( T^{-\varepsilon} M_1\right)^{-2} T^{10\varepsilon -1} 
\left( T^{2\varepsilon} K_2\right)^{2\vartheta} 
K_2^4 M_1^8 \| a\|_2^4 < \cr 
 &< |\varphi|^{-8} \left|\Phi_1\right|^{-4} 
T^{13\varepsilon -1} K_2^{4+2\vartheta} M_1^6 \| a\|_2^4\;.}\eqno(7.69)$$
Since we certainly have here $M_1^6\leq M_1^4 M_1^{4-3\vartheta}$, the combination 
of (7.69) with (7.40) and (7.46) yields the bound (7.33) for $Y$. 
\par 
By the hypothesis (1.20), we have furthermore  
$$\| a\|_2^4\leq 64e^{2\eta}M_1^2 \| a\|_{\infty}^4 
<128 M_1^2 \| a\|_{\infty}^4\;.\eqno(7.70)$$ 
Since we have $M_1^6 M_1^2=M_1^8\leq M_1^{10-5\vartheta}$, 
the combination of (7.69), (7.70), (7.46) and 
(7.40) yields the bound (7.34) for $Y$.  
\par 
Our work in the last four paragraphs has shown that 
if both the inequality (7.63) and 
the bound (7.12) hold, 
then so do the bounds stated in (7.33) and (7.34);  it follows that, 
for the reasons given in the paragraph immediately below (7.62),  
the proof of the lemma is now complete \quad$\square$ 

\bigskip 

\noindent{\Smallcaps The concluding steps of the proof of Theorem~2.}\quad  
Since $T$ satisfies (6.1), it follows from Lemma~20 and Lemma~21 that 
we certainly have 
$$\eqalign{\sum_{d=-\infty}^{\infty}\ \int\limits_{-\infty}^{\infty}
\left| c(d,it)\right|^2
w_0\!\left({|2d+it|^2\over\pi^2 T}\right) {\rm d}t
 &=2\pi\widehat{w_0\circ{\frak N}}(0)\,T\,\| C\|_2^2 
+(\pi /2)\left( {\cal D}^{\star}_1+{\cal D}^{\star}_2+{\cal E}^{\star}\right) \,+ \cr 
 &\quad\,+O_{\eta ,\varepsilon}\!\left( 
\left( T^{5\varepsilon + 1/2} K_1 K_2 +K_2^2\right) M_1 \| a\|_2^2\right)\qquad\quad\   
\hbox{(say),}}\eqno(7.71)$$
where 
$$\widehat{w_0\circ{\frak N}}(0) 
=\int\limits_{-\infty}^{\infty} \int\limits_{-\infty}^{\infty}  
w_0\!\left( x^2 + y^2\right) {\rm d}x\,{\rm d}y 
=\int\limits_0^{\infty} w_0\!\left( r^2\right) {\rm d}\!\left( \pi r^2\right) 
=\pi\int\limits_0^{\infty} w_0(x) {\rm d}x\;,\eqno(7.72)$$ 
while the complex numbers ${\cal D}^{\star}_1$ and ${\cal D}^{\star}_2$ are  
defined by the equations (1.23)-(1.28),  
and the complex number ${\cal E}^{\star}$ is 
given by the equations (6.6)-(6.8), (6.2) and (6.3) 
(it being assumed that $N$ is defined by (1.24)). 
Moreover, since (6.46), (6.43) and (6.2) imply that one has   
$$T^{1+\varepsilon}\left|\varphi z_2 z_3\right|^2 
\ll {|\varphi|^4 \left|\Phi_1\right|^2 T^{1+3\varepsilon}\over \Xi Z}  
= {|\varphi|^4 \left|\Phi_1\right|^2 T^{1+\varepsilon} K_1\over M_1^2 K_2}\;,$$  
it follows  by Lemma~22 and Lemma~28 that the term ${\cal E}^{\star}$ 
occurring in Equation~(7.71) satisfies both 
$${\cal E}^{\star} \ll_{\eta ,\varepsilon} 
T^{6\varepsilon} K_1 K_2^2 M_1 \| a\|_2^2 
+\left( {M_1^2\over T^{1/2}} + 
\left( {M_1^{2-(3/2)\vartheta}\over T^{(1-\vartheta)/2}}\right) 
\left( {K_2\over T^{1/2}}\right)^{\!\!\vartheta}\right) 
T^{1+9\varepsilon} K_1 K_2 \| a\|_2^2 
\eqno(7.73)$$ 
and  
$$\eqalign{ 
{\cal E}^{\star} &\ll_{\eta ,\varepsilon} 
\,T^{6\varepsilon} K_1 K_2^2 M_1 \| a\|_2^2 \ + \cr 
 &\qquad\ +\left( {M_1^2\over T^{1/2}} + 
\left(\!{M_1^2\over T^{1/2}}\right)^{\!\!1-\vartheta}
\!\left(\!\left( {K_1\over T^{1/2}}\right)^{\!\!\vartheta /2} 
\!\!\left( {K_2\over T^{1/2}}\right)^{\!\!1-(\vartheta /2)} + 
\left( {K_2\over T^{1/2} M_1^{1/2}}\right)^{\!\!\vartheta} 
\right)\!\right) 
T^{1+11\varepsilon} K_1 K_2 M_1 \| a\|_{\infty}^2\;.}\eqno(7.74)$$ 
\par 
By (6.1) and the hypotheses of the theorem ((1.18) in particular),  
we have also: 
$$\eqalign{ 
\left( T^{5\varepsilon +1/2} K_1 K_2 + K_2^2 +T^{6\varepsilon} K_1 K_2^2\right) 
M_1 \| a\|_2^2 
 &\ll T^{6\varepsilon +1/2} K_1 K_2 M_1 \| a\|_2^2 \leq \cr 
 &\leq T^{9\varepsilon +1/2} K_1 K_2 M_1^2\| a\|_2^2 = \cr 
 &=\left( M_1^2 T^{-1/2}\right)\left( T^{1+9\varepsilon} K_1 K_2 \| a\|_2^2\right) .} 
\eqno(7.75)$$ 
Moreover, as noted within the proof of Lemma~28 (in (7.70)), we have 
$$\| a\|_2^2\ll M_1 \| a\|_{\infty}^2\;.\eqno(7.76)$$ 
\par 
By (7.71)-(7.73) and (7.75), the term ${\cal E}$ determined   
by the equations (1.21)-(1.30) satisfies  
$$\eqalign{ 
{\cal E} &= O_{\eta,\varepsilon}\!\left(\left( T^{5\varepsilon +1/2} K_1 K_2 + K_2^2\right) 
M_1 \| a\|_2^2\right) + (\pi /2)\,{\cal E}^{\star} \ll_{\eta,\varepsilon}  \cr 
 &\ll_{\eta,\varepsilon}  
\left( {M_1^2\over T^{1/2}} + 
\left( {M_1^{2-(3/2)\vartheta}\over T^{(1-\vartheta)/2}}\right) 
\left( {K_2\over T^{1/2}}\right)^{\!\!\vartheta}\right) 
T^{1+9\varepsilon} K_1 K_2 \| a\|_2^2\;.}$$
We therefore have the bound stated in (1.31), which is the first result of the theorem. 
Similarly, given (1.30) and (1.22), and given that $\varepsilon$ is positive, 
it follows from (7.71), (7.72), (7.74)-(7.76) and (6.1) that the bound (1.32) 
holds\quad$\square$ 

\bigskip 

\goodbreak\proclaim{\Smallcaps Remarks~29}. In the proof of Theorem~2 just concluded, 
we have avoided making any direct use of the `Weil-Estermann bound' bound for Kloosterman sums that 
is stated in (5.29) 
(it should nevertheless be noted that we have used the bound (5.29) in 
proving the results of [45, Theorem~B and Theorem~1], 
upon which the proofs of Lemma~23 and Lemma~25 of the present paper are dependent).  
The sum of Kloosterman sums $Y$ defined in (6.47) may, of course, 
be estimated by means of a direct application of the result in (5.29). 
This approach to the estimation of $Y$ leads one,  
via Lemma~20 (in a slightly sharper revised form), Lemma~21, Lemma~22 and (7.72),  
to the upper bound 
$${\cal E}\ll_{\eta ,\varepsilon} T^{6\varepsilon} K_1 K_2^{5/2} M_1^{5/2} \| a\|_2^2\;,\eqno(7.77)$$ 
where ${\cal E}$ is the final term in the equation (1.30) of Theorem~2 
(it being assumed here that all the hypotheses of that theorem are satisfied). \hfill\break 
\indent Recall that, in addition to (7.77), we have also the bounds (1.31) and (1.32) 
for ${\cal E}$. Yet another bound for the term ${\cal E}$ may be  
by obtained by using both the inequality $xy\leq (x^2+y^2)/2$ 
and the elementary bound 
$|\{ \delta\in{\frak O} : 0<|\delta|\leq r\}|\leq 8r^2$ 
to estimate that part of the sum on the right-hand side of Equation~(6.10) 
in which one has $\xi_1\neq\xi_2$; by this method one finds that, 
subject to the hypotheses of Theorem~2, the term ${\cal E}$ in (1.30) 
must satisfy 
$${\cal E}=\cases{ 
O_{\eta,\varepsilon}\left( T^{2\varepsilon} K_1^2 K_2^2 M_1 \| a\|_2^2\right) 
 &if $\,T^{\varepsilon -1} K_1 K_2 M_1 > e^{-\eta}$, \cr 
O_{\eta,\varepsilon , A}\left( T^{-A} \| a\|_2^2\right) &otherwise,}\eqno(7.78)$$ 
where $A$ denotes an arbitrarily large positive constant. \hfill\break 
\indent Note that the only cases in which (7.77) is not weaker then (7.78) are 
those in which one has both 
$$K_1 K_2 M_1\gg T^{1-\varepsilon}\qquad{\rm and}\qquad 
K_1^2\gg T^{8\varepsilon} K_2 M_1^3\;.\eqno(7.79)$$ 
Moreover, when the conditions in both (7.79) and (1.18) are satisfied, 
the bound for ${\cal E}$ in (7.77) is weaker than that which is implied by (1.31). 
We therefore conclude that, in every case in which the 
hypotheses of Theorem~2 are satisfied, 
the combination of (1.31) with the elementary result (7.78) 
provides a bound for ${\cal E}$ that is as at least as strong 
as that in (7.77). 

\bigskip

\goodbreak 
\noindent{\bf Acknowledgements.}
\ The work of the author during 2004-6 was supported through a
fellowship associated with an EPSRC funded project, 
`The Development and Application of Mean-Value Results in Multiplicative
Number Theory' (GR/T20236/01), which was led by Glyn Harman of Royal Holloway, 
University of London. The author thanks Glyn for his 
encouragement, and for his advice on the application of Theorem~1 of 
this paper to a problem concerning the distribution of Gaussian primes. 
\par 
The author is very grateful to Roelof Bruggeman, Aleksandar Ivi\'c and Yoichi Motohashi 
for their comments on earlier draughts of this work; he also wishes to thank   
Aleksandar Ivi\'c for the gift of an offprint copy of [21]. 
\par 
The author thanks his parents for their encouragement and support. 

\vskip 5mm  

\goodbreak
\noindent{\SectionHeadingFont Bibliography}

\medskip

\parskip = 3 pt 

\item{1.} {\BibAuthorFont R. C. Baker, G. Harman, J. Pintz}, 
{\it The difference between consecutive primes, II}, Proc.
Lond. Math. Soc. (3) {\bf 83} (2001), 532-562.

\item{2.} {\BibAuthorFont R. Balasubramanian}, 
{\it On the fequency of Titchmarsh's phenomenon for $\zeta(s)$, IV},
Hardy-Ramanujan J. {\bf 9} (1986), 1-10.

\item{3.} {\BibAuthorFont V. Blomer, G. Harcos, P. Michel},  
{\it A Burgess-like subconvexity bound for twisted $L$-functions}, 
with and appendix by Z. Mao, Forum Math. {\bf 19} (2007), 61-106. 

\item{4.} {\BibAuthorFont R. W. Bruggeman, R. J. Miatello},
{\it Estimates of Kloosterman sums for groups of real rank one},
Duke Math.~J. {\bf 80} (1995), 105-137.

\item{5.} {\BibAuthorFont R. W. Bruggeman, Y. Motohashi}, 
{\it Sum formula for Kloosterman sums and fourth moment of the Dedekind
zeta-function over the Gaussian number field}, Funct. Approx. Comment. Math. 
{\bf 31} (2003), 23-92.

\item{6.} {\BibAuthorFont H. Davenport}, 
{\it Multiplicative Number Theory}, 
Graduate Texts in Mathematics, No. 74, 
2nd ed., revised by H. L. Montgomery, 
Springer-Verlag, New York, 1980.  

\item{7.} {\BibAuthorFont J.-M. Deshouillers, H. Iwaniec},
{\it Kloosterman sums and Fourier coefficients of cusp forms},
Invent. Math. {\bf 70} (1982), 219-288.

\item{8.} {\BibAuthorFont J.-M. Deshouillers, H. Iwaniec},
{\it Power mean-values of the Riemann zeta-function}, 
Mathematika {\bf 29} (1982), 202-212.

\item{9.} {\BibAuthorFont J.-M. Deshouillers, H. Iwaniec},
{\it Power mean-values for Dirichlet's polynomials and the Riemann zeta-function,~II},
Acta Arith. {\bf 43} (1984), 305-312.

\item{10.} {\BibAuthorFont W. Duke}, 
{\it Some problems in multidimensional analytic number theory},
Acta Arith. {\bf 52} (1989), 203-228.

\item{11.} {\BibAuthorFont J. Elstrodt, F. Grunewald, J. Mennicke}, 
{\it Groups Acting on Hyperbolic Space: Harmonic Analysis and Number Theory}, 
Springer Monographs in Mathematics, Springer, Berlin, 1997.

\item{12.} {\BibAuthorFont G. Harcos}, 
{\it Uniform approximate functional equation for principal $L$-functions},
Int. Math. Res. Not. (2002) No.~18, 923-932. 

\item{13.} {\BibAuthorFont G. Harcos}, 
{\it Uniform approximate functional equation for principal $L$-functions (Erratum)}, 
Int. Math. Res. Not. (2004) No.~13, 659-660.

\item{14.} {\BibAuthorFont G. H. Hardy, E. M. Wright}, 
{\it An Introduction to the Theory of Numbers}, 
5th~ed., reprinted (with corrections), Oxford University Press, 1984. 

\item{15.} {\BibAuthorFont G. Harman}, 
{\it Prime-Detecting Sieves}, 
London Mathematical Society Monographs, Vol.~33, 
Princeton University Press, Princeton, 2007. 

\item{16.} {\BibAuthorFont G. Harman, A. Kumchev, P. A. Lewis},
{\it The distribution of prime ideals of imaginary quadratic fields},
Trans. Amer. Math. Soc. {\bf 356} (2003), 599-620.

\item{17.} {\BibAuthorFont E. Hecke}, 
{\it Eine neue art von zetafunktionen und ihre beziehungen
zur verteilung der primzahlen (erste mitteilung)}, 
Math.~Z. {\bf 1} (1918), 357-376. 

\item{18.} {\BibAuthorFont E. Hecke}, 
{\it Eine neue art von zetafunktionen und ihre beziehungen
zur verteilung der primzahlen (zweite mitteilung)}, 
Math.~Z. {\bf 6} (1920), 11-51.

\item{19.} {\BibAuthorFont M. N. Huxley}, 
{\it The large sieve inequality for algebraic number fields.~II:
means of moments of Hecke zeta-functions},
Proc. Lond. Math. Soc. (3) {\bf 21} (1970), 108-128.

\item{20.} {\BibAuthorFont A. E. Ingham}, 
{\it The Distribution of Prime Numbers}, 
Cambridge Tracts in Mathematics and Mathematical Physics, No.~30, 
reissued as a volume in the Cambridge Mathematical Library series, 
with a foreward by R. C. Vaughan, Cambridge University Press, Cambridge, 1990.

\item{21.} {\BibAuthorFont A. Ivi\'{c}}, 
{\it An approximate functional equation for a class of Dirichlet series}, 
J. Anal. {\bf 3} (1995), 241-252.

\item{22.} {\BibAuthorFont A. Ivi\'{c}}, 
{\it The Riemann Zeta-Function: Theory and Applications},
reprinted with a list of errata, Dover, Mineola, 2003.

\item{23.} {\BibAuthorFont H. Iwaniec}, 
{\it On mean values for Dirichlet's polynomials and the Riemann zeta-function}, 
J. London Math. Soc.~(2) {\bf 22} (1980), 39-45.

\item{24.} {\BibAuthorFont H. Iwaniec, E. Kowalski}, 
{\it Analytic Number Theory},
American Mathematical Society Colloquium Publications, Vol.~53, 
American Mathematical Society, Providence, Rhode Island, 2004.

\item{25.} {\BibAuthorFont M. Jutila}, 
{\it A variant of the circle method}, in 
`Sieve Methods, Exponential Sums and their Applications in Number Theory: 
Proceedings of a symposium held in Cardiff, July 1995', 
G. R. H. Greaves, G. Harman, M. N. Huxley (eds.), 
London Mathematical Society Lecture Note Series, No.~237, 
Cambridge University Press, Cambridge, 1997, pp. 245-254. 

\item{26.} {\BibAuthorFont R. M. Kaufman}, 
{\it An estimate of Hecke's $L$-functions of the Gaussian field 
on the line ${\rm Re}\,s = 1/2$},
Dokl. Akad. Nauk. BSSR {\bf 22} (1978), 25-28 (in Russian).

\item{27.} {\BibAuthorFont R. M. Kaufman}, 
{\it Estimate of the Hecke $L$-functions on the critical line}, 
in `Analytic Number Theory and the Theory of Functions. Part~2', 
Zap. Nauchn. Sem. LOMI {\bf 91}, ``Nauka", Leningrad. Otdel., Leningrad, 1979, 
pp. 40-51 (in Russian, with a summary in English). 

\item{28.} {\BibAuthorFont H.H. Kim},  
{\it Functoriality for the exterior square of $GL_4$ and the symmetric fourth of $GL_2$}, 
with Appendix 1 by D.~Ramakrishnan, Appendix 2 by Kim, P.~Sarnak, 
J.~Amer. Math. Soc. {\bf 16} (2003),~139-183.  

\item{29.} {\BibAuthorFont H. H. Kim}, 
{\it On local $L$-functions and normalized intertwining operators}, 
Canad.~J. Math. {\bf 57} (2005), 535-597. 

\item{30.} {\BibAuthorFont H. H. Kim, F. Shahidi}, 
{\it Cuspidality of symmetric powers with applications}, 
Duke Math.~J. {\bf 112} (2002), 177-197. 

\item{31.} {\BibAuthorFont S. Lang}, 
{\it Real Analysis}, 2nd~ed., Addison-Wesley, Reading, Massachusetts, 1983. 

\item{32.} {\BibAuthorFont A. F. Lavrik}, 
{\it An approximate functional equation for the 
Hecke zeta function of an imaginary quadratic field}, 
Mat. Zametki {\bf 2} (1967), 475-482 (in Russian). 

\item{33.} {\BibAuthorFont P. A. Lewis}, 
{\it Finding information about Gaussian primes 
using analytic number theory sieve methods}, 
Ph.D. thesis, Cardiff University, 2002. 

\item{34.} {\BibAuthorFont H. Lokvenec-Guleska}, 
{\it Sum formula for $SL_2$ over imaginary quadratic number fields\/}, 
doctoral thesis, University of Utrecht, 2004 
(includes summaries in Dutch and Macedonian).

\item{35.} {\BibAuthorFont Y. Motohashi}, 
{\it Spectral Theory of the Riemann Zeta-Function}, 
Cambridge Tracts in Mathematics, No.~127, 
Cambridge University Press, Cambridge, 1997.

\item{36.} {\BibAuthorFont H. Rademacher}, 
{\it On the Phragm\'en-Lindel\"of theorem and some applications}, 
Math.~Z. {\bf 72} (1959), 192-204. 

\item{37.} {\BibAuthorFont K. Ramachandra}, 
{\it On the fequency of Titchmarsh's phenomenon for $\zeta(s)$,~I},
J. London Math. Soc.~(2) {\bf 8} (1974), 683-690.

\item{38.} {\BibAuthorFont K. Ramachandra}, 
{\it On the Mean-Value and Omega-Theorems for the Riemann Zeta-Function},
Tata Institute of Fundamental Research Lectures on Mathematics and Physics, No.~85, 
Nairosa Publishing House, Bombay, 1995. 

\item{39.} {\BibAuthorFont P. Sarnak}, 
{\it Fourth moments of grossencharakteren zeta functions},
Comm. Pure and Applied Math. {\bf 38} (1985), 167-178.

\item{40.} {\BibAuthorFont A. Selberg}, 
{\it On the estimation of Fourier coefficients of modular forms}, 
in `Theory of Numbers: Proceedings of the Symposium in Pure Mathematics 
Held at the California Institute of Technology, Padadena, California, 
November, 21, 22, 1963', A. L. Whiteman (editor), 
Proceedings of Symposia in  Pure Mathematics, Vol.~VIII, 
American Mathematical Society, Providence, Rhode Island, 1965, pp.~1-15. 

\item{41.} {\BibAuthorFont P. S\"{o}hne}, 
{\it An upper bound for Hecke zeta-functions with groessencharacters}, 
J. Number Theory {\bf 66} (1997), 225-250.

\item{42.} {\BibAuthorFont E. C. Titchmarsh}, 
{\it The Theory of Functions}, 2nd~ed., reprinted with corrections, 
Oxford University Press, Oxford, 1983.

\item{43.} {\BibAuthorFont N. Watt}, 
{\it Kloosterman sums and a mean value for Dirichlet polynomials}, 
J. Number Theory {\bf 53} (1995), 179-210.

\item{44.} {\BibAuthorFont N. Watt},
{\it Weighted spectral large-sieve inequalities for Hecke congruence subgroups of
${\rm SL}(2,{\Bbb Z}[i])$}, 
Funct. Approx. Comment. Math. {\bf 48} (2013), 213-376.   

\item{45.} {\BibAuthorFont N. Watt},
{\it Spectral large sieve inequalities for Hecke congruence subgroups of
${\rm SL}(2,{\Bbb Z}[i])$}, preprint, arXiv:1302.3112v1 [math.NT].

\item{46.} {\BibAuthorFont E. T. Whittaker, G. N. Watson}, 
{\it A Course of Modern Analysis}, 4th~ed., reprinted, 
Cambridge University Press, New York, 1973. 

\bye